\documentclass{daj}

\usepackage{amsmath,amsthm,amssymb,enumitem,mathtools,bbm,amsfonts,wasysym,stmaryrd}
\newtheorem{thm}{Theorem}[section]
\newtheorem{lem}[thm]{Lemma}
\newtheorem{prop}[thm]{Proposition} 

\newtheorem{conj}[thm]{Conjecture}

\newtheorem{claim}[thm]{Claim}
\newtheorem{fact}[thm]{Fact}

\newtheorem{remark}[thm]{Remark}
\newtheorem{obs}[thm]{Observation}
\newtheorem{defi}[thm]{Definition}
\theoremstyle{definition}

\theoremstyle{plain}

\newcommand{\Dirac}{\boldsymbol{\delta}}

\newcommand{\diffsymb}{\mathsf{d}} 

\newcommand{\interval}[2]{\mathsf{J}_{#1,#2}}
\newcommand{\dyainterval}[2]{\mathsf{D}_{#1,#2}}
\newcommand{\partition}[1]{\mathfrak{J}_{#1}}
\newcommand{\dyapartition}[1]{\mathfrak{D}_{#1}}
\newcommand{\alldyapartition}{\mathfrak{D}_*}
\newcommand{\eps}{\varepsilon}
\newcommand{\vphi}{\varphi}

\newcommand{\Prob}{\mathbb{P}}
\newcommand{\Expect}{\mathbb{E}}

\newcommand{\LeftConvergence}{\overset{\mathrm{left}}\rightarrow}
\newcommand{\cutn}[1]{\left\lVert #1\right\rVert_{\square}}
\newcommand{\cutnL}[1]{\left\lVert #1\right\rVert_{\XBox}}
\newcommand{\deltaL}{\delta_{\XBox}}
\newcommand{\deltaCut}{\delta_{\square}}
\newcommand{\LatinonSpace}{\mathcal{L}}
\newcommand{\AntiCompr}[1]{\nabla(#1)}
\newcommand{\inftyd}[1]{\delta_\square^{\mathbb N}(#1)}

\newcommand{\Entropy}{\mathrm{Ent}}

\dajAUTHORdetails{%
  title = {Limits of Latin Squares}, 
  author = {Frederik Garbe, Robert Hancock, Jan Hladk\'y, and Maryam Sharifzadeh},
  plaintextauthor = {Frederik Garbe, Robert Hancock, Jan Hladky, and Maryam Sharifzadeh},
  keywords = {Latin square, Latinon, limits of discrete structures, graphon},
}   

\dajEDITORdetails{%
   year={2023},
   lastpage={66},
   number={8},
   received={16 October 2020},   
   published={20 July 2023},  
   doi={10.19086/da.83253},       
}   

\begin{document}

\begin{frontmatter}[classification=text]


\author[fg]{Frederik Garbe\thanks{Supported by GA\v{C}R project 18-01472Y and RVO: 67985840. This work was done while affiliated with the Institute of Mathematics of the Czech Academy of Sciences.}}
\author[rh]{Robert Hancock\thanks{Supported by GA\v{C}R project 18-01472Y and RVO: 67985840. This work has received funding from the European Research Council (ERC) under the European Union's Horizon 2020 research and innovation programme (grant agreement No 648509) and from the MUNI Award in Science and Humanities of the Grant Agency of Masaryk University. This publication reflects only its authors' view; the European Research Council Executive Agency is not responsible for any use that may be made of the information it contains.}}
\author[jh]{Jan Hladk\'y\thanks{Supported by GA\v{C}R project 18-01472Y and RVO: 67985840. This work was done while affiliated with Institute of Mathematics of the Czech Academy of Sciences.}}
\author[ms]{Maryam Sharifzadeh
	\thanks{This work has received funding from the European Union's Horizon 2020 research and innovation programme under the Marie
		Curie Individual Fellowship agreement No 752426 and was done while affiliated with University of Warwick.}
}

\begin{abstract}
We develop a limit theory of Latin squares, paralleling the recent limit theories of dense graphs and permutations. We introduce a notion of density, an appropriate version of the cut distance, and a space of limit objects --- so-called Latinons. Key results of our theory are the compactness of the limit space and the equivalence of the topologies induced by the cut distance and the left-convergence. Last, using Keevash's recent results on combinatorial designs, we prove that each Latinon can be approximated by a finite Latin square.
\end{abstract}
\end{frontmatter}

\bigskip

An extended abstract describing this work appeared in the proceedings of EuroComb2019,~\cite{GHHS:LatinSquareEurocomb}.
\setcounter{tocdepth}{2}
\tableofcontents

\section{Introduction}\label{sec:intro}
The purpose of this paper is to lay down the basics of a combinatorial theory of limits of Latin squares. This is a further addition to a potent approach to discrete structures using analytic tools. In combinatorics, this research trend blossomed in the early 2000's with the theories of dense graph limits initiated by Borgs, Chayes, Lov\'asz, S\'os, Szegedy, and Vesztergombi,~\cite{MR2455626,MR2925382,Lovasz2006}, and sparse graph limits initiated by Benjamini and Schramm~\cite{MR1873300}. After that, many other theories followed, and Razborov provided an alternative, more syntactic limit framework, \cite{Razborov2007}.

Recall that a \emph{Latin square} is an $n\times n$ matrix filled with values of $[n]\coloneqq\{1,\dots,n\}$ in such a way that each row contains each value exactly once, and each column contains each value exactly once, as well. Let us emphasise that the rows, columns, and symbols in the matrix are ordered (from top to bottom, from left to right, and from small to large). 
In the same spirit, a \emph{permutation} is a bijection from the naturally ordered set $[n]$ into itself.

In many aspects, our theory parallels those of limits of dense graphs and of permutations. In order to make a comparison, let us quickly recall the key features of these theories. For foundations of these two respective theories, see~\cite{Lovasz2006} and~\cite{HoKo13}. So, by a `structure', we mean either a finite graph or a finite permutation.\footnote{But the same could be said about limit theories of many other discrete structures such as uniform hypergraphs.}
\begin{enumerate}[label={(F\arabic*)}]
	\item\label{en:F1} Both theories introduce a certain space of analytically defined limit objects (\emph{graphons} and \emph{permutons}).
	\item\label{en:F2} Both theories introduce a notion of \emph{densities} $t(\circ,\star)$ of a structure $\circ$ both into a structure $\star$, or into a limit object $\star$.
	\item\label{en:F3} Both theories also provide more global, `cut-like', parameters on structures and on limit objects.
\end{enumerate}
There are three key features in both theories.
\begin{enumerate}[label={(F\arabic*)},resume]
	\item\label{en:F4} Compactness theorem: Each sequence of structures contains a subsequence that converges (more precisely, `left-converges') to a suitable limit object with respect to all densities.
	\item\label{en:F5} Equivalence of left-convergence and convergence with respect to `cut-distance': A sequence of structures converges with respect to all densities if and only if it converges with respect to the above `cut-like' parameters.
	\item\label{en:F6} Sampling lemma: The densities of structures in a limit object provide a notion of random sampling of a structure $S_n$ of order $n$ from that limit object. The sampling lemma says that $S_n$ is typically very close to the limit object. More specifically, almost surely, the sequence $(S_n)_n$ converges to the original limit object. As a consequence of the sampling lemma, each limit object can be approximated by a structure with an arbitrary precision.
\end{enumerate}

In our work, we introduce a class of analytic objects, so-called \emph{Latinons}. We also introduce certain notions of densities for Latin squares and Latinons (however, with a slight difference compared to~\ref{en:F2} which we shall return to later) and `cut-like' parameters, and establish~\ref{en:F5}. Of course, the most important result is the compactness theorem~\ref{en:F4}, which we establish as well. The situation with~\ref{en:F6} is more complicated.
Given a Latinon $L$, we do not have any reasonable way of generating a random Latin square of order $n$ from $L$ (and probably no such way exists, see Remark~\ref{rem:chocolatenotsublatin}). However, we can still establish the approximability consequence of~\ref{en:F6} (which is the most important part of~\ref{en:F6}).

Theories of limits of sequences of particular discrete structures have found many applications and generated many new exciting questions, but two particular directions stand out. First and foremost, they (and the more syntactic theory of flag algebras, \cite{Razborov2007}) have led to solutions to many asymptotic problems in extremal combinatorics. As some representative examples for the class of graphs, let us mention the solution to the triangle versus edge density problem~\cite{Razborov:Triangles}, the solution to the Erd\H{o}s problem on the number of pentagons in triangle-free graphs~\cite{HHKNR:Pentagons,Grzesik:Pentagons,LidPfe:Pentagons}, or progress on the Sidorenko conjecture~\cite{Lov:Sidorenko}. For permutations, we are aware of only one application so far, but it is a remarkable one. Namely, Kr\'al' and Pikhurko~\cite{KrPi} used the theory of permutation limits to characterise sequences of quasirandom permutations. We expect the theory we develop in this article to have applications mostly in this area, that is, to answer asymptotic extremal questions about Latin squares. The set of tools we provide in this article should be complete for such purposes. Indeed, very shortly after we published this preprint at arXiv, Cooper, Kr\'al', Lamaison, Mohr~\cite{CKLM:QuasirandomLatin} used our theory to characterise sequences of quasirandom Latin squares; see Section~\ref{ssec:quasirandomness}. A second line of research where limits of discrete structures opened new possibilities is probabilistic. Let us give several examples, though the list below is by no means exhaustive. Chatterjee and Varadhan~\cite{MR2825532} developed a method of expressing certain density-based large deviation events of Erd\H{o}s--R\'enyi random graphs. This theory in particular allows us to approximately count graphs with a prescribed list of density constraints. A similar programme was carried out for permutations as well,~\cite{KeKrRaWi:Permutations}. The framework of large deviations was recently applied (with substantial technical sophistication) by Simkin~\cite{Simkin:Queens} to `queenons', an extension of the theory of permutons tailored to describe configurations of $n$ non-threatening queens on an $n\times n$ chessboard ($n\to\infty$). In particular, Simkin obtains the asymptotics for the number of such configurations thus solving a famous longstanding open problem. We do not offer such tools for Latin squares in the present article, but in Section~\ref{ssec:Entropy} we discuss a reasonable approach. Another line of research comes from studying random objects coming from~\ref{en:F6}, leading to models of inhomogeneous random graphs and random permutations. Some basic properties of these models are studied in~\cite{HlaVis:Connectivity,DHM:Cliques,IK:ChromaticGraphon} (for graphs) and~\cite{MR1334175} (for permutations). Even inhomogeneous random models that do not come directly from the sampling procedure~\ref{en:F6} can be studied using limit theories, see e.g.~\cite{MR3813988,MR4516311,MR4295570,MR4359933}. In addition to describing individual random objects, limit theories have been also useful in describing stochastic processes, see e.g. \cite{Flip1,LargeDeviationsProcesses} (for graphs) and \cite{MR4420996,MR4015358} (for permutations).

\subsection{Work on limits of related discrete structures}
The crucial feature of the structures whose limit theory we introduce is that they are ordered; that is, the rows go from top to bottom, the columns go from left to right, and the entries go from small to large. Arguably, the first class of ordered structures for which a similar limit theory was developed are posets,~\cite{Janson:PosetLimits,HlMaPaPi:Posets}. However, the most related work is a recent paper of Ben-Eliezer, Fischer, Levi, and Yoshida~\cite{BeFiLeYo18} (see also conference version~\cite{BeFiLeYo18C}). In that paper, a theory of limits of ordered graphs, that is, graphs whose vertices are put into a linear order, is developed. The limit objects they obtain (which they call `orderons') is less complex than our Latinons. More precisely, the fact that they are dealing with 0/1-matrices (adjacency matrices of ordered graphs) whereas we are dealing with more general matrices (the values are the set $[n]$) is reflected in that their limit objects are $[0,1]$-valued, while ours are distribution-valued. The necessity of such a more complex limit concept is best illustrated by examples, which we provide in Section~\ref{ssec:ExamplesOfLatinons}. Let us point out that the key difference between~\cite{BeFiLeYo18} and our paper is how we view and treat the limit structure. We believe that our approach is conceptually simpler. We give more details in Section~\ref{ssec:orderonscompactness}.

We already mentioned that our Latinons are distribution-valued, a feature that causes substantial technical complications. While indeed most limit theories lead to real-valued limit objects, the more general setting was also considered. In this direction, the most relevant paper is~\cite{KuLoSz:LimBanachDec} which introduces graphons taking values from the dual of a separable Banach space. Indeed, our Latinons take values which are probability measures on $[0,1]$, and these can be viewed as elements in the dual of the space of continuous functions on $[0,1]$ (equipped with the uniform norm). However, the notion of densities considered in~\cite{KuLoSz:LimBanachDec} is substantially different and we could not find a way to directly relate our setting to that of~\cite{KuLoSz:LimBanachDec}.

\subsection{Representation of a Latin square in a 3-dimensional matrix}\label{ssec:representationAs3Dim}
Suppose that $L$ is a Latin square of order $n$. Then we can associate $L$ with the following 3-uniform hypergraph. It has $3n$ vertices, which we label as $\{(i,\mathtt{row})\}_{i\in[n]}$, $\{(j,\mathtt{column})\}_{j\in [n]}$, $\{(k,\mathtt{value})\}_{k\in[n]}$. The hypergraph is partite, that is, only edges of the form $\mathtt{row}$-$\mathtt{column}$-$\mathtt{value}$ are present. More precisely, $(i,\mathtt{row}), (j,\mathtt{column}), (k,\mathtt{value})$ forms an edge if and only if $L_{i,j}=k$. Observe that the hypergraph has $n^2$ edges, and that
\begin{enumerate}[label=(L\arabic*)]
	\item\label{L1} for every $i,j\in [n]$ there exists a unique $k\in[n]$ such that $(i,\mathtt{row}), (j,\mathtt{column}), (k,\mathtt{value})$ forms an edge;
	\item\label{L2} for every $i,k\in [n]$ there exists a unique $j\in[n]$ such that $(i,\mathtt{row}), (j,\mathtt{column}), (k,\mathtt{value})$ forms an edge;
	\item\label{L3} for every $j,k\in [n]$ there exists a unique $i\in[n]$ such that $(i,\mathtt{row}), (j,\mathtt{column}), (k,\mathtt{value})$ forms an edge.
\end{enumerate}
Actually, this defines a one-to-one correspondence between Latin squares and 3-uniform hypergraphs with the above properties. That is, the key features~\ref{L1}--\ref{L3} could lead to a limit theory of certain 3-uniform hypergraphs of density of relative exponent $2/3$. Our understanding of general limit theories of structures of intermediate densities is limited,\footnote{See~\cite{MR3988601,MR3758733} for some progress in the case of graphs.} so such a theory may be surprising. We briefly comment on a possible limit theory stemming from interpreting~\ref{L1}--\ref{L3} symmetrically in Section~\ref{ssec:NewDiscussion}, since, crucially, our view of Latin squares in this paper is different. In particular, while the role of rows and columns is symmetric, the role of values in a Latin squares is different, because features~\ref{L2} and~\ref{L3} play a different role from that of~\ref{L1}. For example, if a sequence of Latin squares converges (in the sense of our main Definition~\ref{def:left}) then it is trivial to check that the sequence in which we swap the rows and columns in those Latin squares converges, too. However, it is not true that swapping, say, columns and values, preserves convergence. We give an example to illustrate this in Section~\ref{ssec:swapMakesDifference}. This asymmetric view arises naturally from our sampling procedure for Latin squares, but might not be desired for a general limit theory for hypergraphs.

Let us also note that a counterpart to the above construction of $3$-hypergraphs can also be made in uniformity 2, where it corresponds to permutations, and in uniformity $\ell>2$, where it leads to hypergraphs with $\Theta(n^{\ell-1})$ edges. These hypergraphs were first considered from this perspective by Linial and Luria~\cite{MR3259813}, who call them \emph{higher dimensional permutations}. Thus it would be natural to ask how much of our theory generalises to this setting. We address this in Section~\ref{ssec:higherdimPerm}.

\subsection{Organisation of the paper}

In Section~\ref{sec:basics} we introduce and summarise necessary basics from measure theory and the theories of limits of graphs and permutations. In Section~\ref{sec:latinondensitycomp} we introduce the basic concepts of our theory, i.e. the definition of the limit objects, the notion of density and state our main compactness result and approximation results, Theorem~\ref{thm:compact} and Theorem~\ref{thm:latinonsarelimits}. We further give some interesting illustrated examples. In Section~\ref{sec:cutdistcountconv} we introduce an analogue of the cut distance for Latinons and state related counting lemmas which imply the equivalence of the two topologies generated by left-convergence and the cut distance, respectively. We also state a sampling lemma for Latinons. In Section~\ref{sec:prelim} we recall some results of the theory of dense graph limits generalised to vectors of bigraphons which we need to prove our results. In Section~\ref{sec:compressions} we introduce a main concept used by our proofs, an approximation of a Latinon by a vector of bigraphons. In Section~\ref{sec:countinglemma} we use this approximation to prove the counting lemma stated in Section~\ref{sec:cutdistcountconv}. In Section~\ref{sec:compactness} we prove the compactness result stated in Section~\ref{sec:latinondensitycomp}. 
In Section~\ref{sec:sampling} we prove the sampling lemma for Latinons, which we in turn use to prove the inverse counting lemma in Section~\ref{sec:invcount}. In Section~\ref{sec:ApproximatingLatinons} we show how to construct an approximating sequence of Latin squares for a given Latinon. Finally, we discuss some open questions in Section~\ref{sec:future}. The appendix contains proofs of two auxiliary results which are pretty standard, but which we did not find in the literature.
The proof of the equivalent definition of a Latinon can be found in Appendix~\ref{app:localglobal}. In Appendix~\ref{app:nibble} we prove a result concerning a decomposition of $K_{n,n,n}$ into triangles which satisfies certain requirements on densities.

\section{Basic definitions and notation}\label{sec:basics}

\subsection{General notation}
For a set $S$ with a strict partial order $<$ and a natural number $k$ we write $S_{<}^{k}$ for the set of all $k$-tuples of $S$ which are in increasing order according to $<$, i.e.
\[S_{<}^{k}\coloneqq\{(x_1,\dots,x_k)\in S^k\mid x_i<x_{i+1} \textrm{ for all } i\in[k-1]\}\;.\]
For a set $S$ we denote by $\mathbbm{1}_S$ the indicator function, i.e. 
\[\mathbbm{1}_S(x)\coloneqq\begin{cases}1,\textrm{ if }x\in S;\\0,\textrm{ otherwise.}\end{cases}\;\]

For every $k\in\mathbb N$ we define a partition $\partition{k}$ of $[0,1]$ into $k$ parts given by $\interval{k}{i} \coloneqq [\frac{i-1}{k},\frac{i}{k})$ for $i \in [k-1]$ and $\interval{k}{k} \coloneqq [\frac{k-1}{k},1]$. 
For every $d\in\mathbb N$ we define a partition $\dyapartition{d}$ of 
$[0,1)$ into $2^d$ parts given by $\dyainterval{d}{s}\coloneqq[\frac{s-1}{2^d},\frac{s}{2^d})$
for $i\in[2^d]$. We also write 
$\alldyapartition\coloneqq\bigcup_{d\in\mathbb N}\dyapartition{d}$.

\subsection{Probability theory}
We will need the following version of the Chernoff bound, see e.g.~\cite[Corollary 2.3]{JaLuRuRandGr00}.
\begin{lem}\label{lem:chernoff}
	If $X\sim\textrm{Bi}(N,p)$, $\mu=Np$, and $\eps>0$, then
	\[\Prob[|X-\mu|\geq\eps\mu]\leq2\exp(-\eps^2\mu/3)\;.\]
\end{lem}

Further, we recall McDiarmid's concentration inequality,~\cite{McDiarmid1989}.
\begin{lem}
	\label{lem:MBoundDiff}Suppose that $r\in\mathbb{N}$,
	and that $Z$ is a random variable on a discrete product probability space $\Lambda\coloneqq \prod_{i=1}^r \Lambda_i$. Suppose that $\mathbf{c}\in \mathbb{R}^r$ is a vector with the following property.
	For each two vectors $\mathbf{x},\mathbf{x}'\in \Lambda$
	that differ on one coordinate, say the $i$-th, we have $\left|Z(\mathbf{x})-Z(\mathbf{x}')\right|\le \mathbf{c}_i$.
	Then for each $t>0$ we have that
	\[
	\Prob\left[\left|Z-\Expect[Z] \right|>t\right]\le\exp\left(-\frac{2t^{2}}{\sum_{i=1}^r \mathbf{c}_i^2}\right)\;.
	\]
\end{lem}

\subsection{Measure theory}\label{ssec:measure}
When working with measure spaces, we will suppress referring to the underlying sigma-algebra, unless it needs to be mentioned explicitly. Recall that measurable sets of zero measure are called \emph{null} and their complements are called \emph{conull}.
We denote the Lebesgue measure on $\mathbb{R}$ by $\lambda$ and by $\lambda^{\otimes n}$ the Lebesgue measure on $\mathbb{R}^n$. Throughout the paper $\Omega$ always denotes an arbitrary separable atomless probability space $(\Omega,\mu)$. The actual choice of $\Omega$ will not be important; recall that every separable atomless probability space is almost  isomorphic to $([0,1],\lambda)$. This means that there exists a conull set $\Omega'\subset\Omega$, a conull set $X\subset [0,1]$ and a measure preserving bijection $b:\Omega'\rightarrow X$ such that its inverse $b^{-1}$ is also measurable. So while we could have worked only with $([0,1],\lambda)$ in the entire paper, we prefer to work with $(\Omega,\mu)$ at the places where we want to emphasise that there is no natural linear order on that probability space. We write $\mu^{\otimes n}$ for the $n$-th power of $\mu$, which is a probability measure on $\Omega^n$.

Let $\mathcal{B}(X)$ denote the space of all Borel probability measures on a metric space $X$ and let $\mathcal{B}^0(X)$ denote the space of all Borel signed measures of total measure between $[-1,1]$ on a metric space $X$. Recall that there is a natural notion of a sigma-algebra on $\mathcal{B}(X)$ and $\mathcal{B}^0(X)$; see \cite[Section 17.E]{MR1321597}. We will not need details, except that this justifies that given a function $f:\Omega\rightarrow \mathcal{B}(X)$ we can ask whether $f$ is measurable or not. We write $\mathcal{B}\coloneqq\mathcal{B}([0,1])$ and $\mathcal{B}^0\coloneqq\mathcal{B}^0([0,1])$. For $p\in[0,1]$, we write $\Dirac_p\in \mathcal{B}$ for the Dirac measure on $p$.

\subsubsection{Construction of measures from semirings}\label{sssec:semirings}

Quite a bit of the technical work deals with constructing measures, in particular in the proof of compactness, where we construct the so-called measure representation of the limit Latinon (see Section~\ref{ssec:Cryptomorphic}). Construction of measures on a sigma-algebra is often conveniently done by extending a premeasure. We recall the necessary definitions.

\begin{defi}[Semiring]
	For a given set $\Omega$, a set $S$ of subsets of $\Omega$ is called a \emph{semiring} if
	\begin{enumerate}[label=(\roman*)]
		\item $\emptyset \in S$,
		\item $A,B\in S$ $\Rightarrow$ $A\cap B \in S$,
		\item if $A,B\in S$, then there exist disjoint sets $K_1,\dots,K_n\in S$ such that $A\setminus B=\bigcup_{i=1}^nK_i$.
	\end{enumerate}
\end{defi}

\begin{defi}[Premeasure]
	Suppose that $S$ is a semiring on $\Omega$. We say that $\mu:S\rightarrow[0,\infty)$ is a \emph{premeasure} if for all sets $A \in S$ for which there exists a countable decomposition $A=\bigcup_{i=1}^\infty A_i$ in disjoint sets $A_i\in S$, $i\in\mathbb N$, we have $\mu(A)=\sum_{i=1}^\infty\mu(A_i)$.
\end{defi}

\begin{thm}[Carath\'eodory's extension theorem]\label{thm:cara}
	Let $S$ be a semiring on $\Omega$ and let $\mu:S\rightarrow[0,1]$ be a premeasure on $S$. Then there exists a measure $\mu':\sigma(S)\rightarrow[0,1]$, where $\sigma(S)$ is the sigma-algebra generated by $S$, such that $\mu'$ is the extension of $\mu$. 
\end{thm}

\subsubsection{Weak convergence}Let us also recall the notion of weak convergence of measures. Suppose that $(X,\rho)$ is a compact metric space. We say that a sequence of Borel measures $\nu_1,\nu_2,\ldots$ on $X$ \emph{weak converges} to a Borel measure $\nu$ if for each continuous function $f:X\rightarrow \mathbb{R}$ we have that the sequence $\int f \diffsymb \nu_1,\int f \diffsymb \nu_2,\ldots$ converges to $\int f \diffsymb \nu$. It is well known that with the compactness assumption as above, if we have a sequence of Borel measures $\nu_1,\nu_2,\ldots$ that are uniformly bounded, say $\nu_n(X)\le C$ for all $n$, then there exists a subsequence $\nu_{n_1},\nu_{n_2},\ldots$ and a Borel measure $\nu$ so that $\nu_{n_1},\nu_{n_2},\ldots$ weak converges to $\nu$. 
This follows immediately from Prokhorov's theorem~\cite{prokhorov}.

\subsubsection{Disintegration theorem}
The disintegration theorem is a standard result in measure theory. Let $Y$ and $X$ be two Radon spaces. Let $\alpha\in \mathcal{B}(Y)$, let $\pi:Y\rightarrow X$ be a Borel-measurable function, and let $\beta\in \mathcal{B}(X)$ be the pushforward measure $\beta=\alpha\circ\pi^{-1}$. Then there exists a $\beta$-almost everywhere uniquely determined family of measures $\{\alpha_x\in \mathcal{B}(Y)\}_{x\in X}$ such that 
\begin{itemize}
	\item the function $x \mapsto \alpha_x$ is measurable,
	\item for $\beta$-almost all $x\in X$, $\alpha_x(Y\setminus \pi^{-1}(x))=0$,
	\item for every measurable function $s:Y\rightarrow\mathbb{R}$, we have 
	\[\int_{y\in Y} s(y) \diffsymb\alpha(y)=\int_{x\in X}\int_{y\in \pi^{-1}(x)}s(y)\diffsymb\alpha_x(y)\diffsymb\beta(x)\;.\]
\end{itemize}
The collection $\{\alpha_x\in \mathcal{B}(Y)\}_{x\in X}$ is called the \emph{disintegration of $\alpha$ according to $\pi$}. Quite often, $Y$ will be of the form $Y=X\times A$ and $\pi$ will be the natural projection. We then say that $\{\alpha_x\in \mathcal{B}(Y)\}_{x\in X}$ is the \emph{disintegration of $\alpha$ according to the first coordinate}.

\subsection{Graphons}
In this paper, we will borrow many tools from the theory of dense graph limits. A \emph{graphon} $W$ is a measurable function $W: \Omega^2 \to [0,1]$ that is symmetric, i.e., $W(x,y)=W(y,x)$. Graphons arise as limits of sequences of graphs, and could be informally understood as limits of their adjacency matrices rescaled into a square of unit area. The symmetricity property is then inherited from the fact that an adjacency matrix of a graph is symmetric. We will need a non-symmetric version of this, called bigraphons, and first introduced in \cite{MR2815610}.
A \emph{bigraphon} is a measurable function $W: \Omega^2 \to [0,1]$. We can think of a bigraphon as a limit of matrices of bipartite graphs $(G_n)_n$ on $n+n$ vertices, where the two parts are distinguished, say \emph{the first part} and \emph{the second part}. Then $W$ can be thought of as a limit of the `bipartite adjacency matrices' $B_n$, in which the entry on position $(i,j)$ is~1 if and only if the $i$-th vertex in the first part of $G_n$ is adjacent to the $j$-th vertex in the second part. The reason why bigraphons will be more relevant in this paper is that the `first part' and the `second part' will encode the rows and the columns of the Latin square respectively. 
Let $\mathcal{W}_0$ be the space of all bigraphons.
A \emph{distribution-valued bigraphon} is a measurable function $W:\Omega^2 \to \mathcal{B}$. A \emph{distribution-valued signed bigraphon} is a measurable function $W:\Omega^2 \to \mathcal{B}^0$.

\subsubsection{Cut norm and cut distance}
Given a function $W$ with domain $\Omega_1^2$ and functions $\vphi,\psi:\Omega_2 \rightarrow \Omega_1$, we define $W^{\vphi,\psi}$ to be the function with domain $\Omega_2^2$ given by $W^{\vphi,\psi}(x,y) \coloneqq W(\vphi(x),\psi(y))$ for all $x,y \in \Omega_2$. We also use the shorthand $W^{\vphi} \coloneqq W^{\vphi,\vphi}$. We use $S_{\Omega}$ to denote the set of all invertible measure preserving maps $\Omega \rightarrow \Omega$.

For $X\in L^\infty(\Omega^2)$, we define the \emph{cut norm} of $X$,
\[ \cutn{X} \coloneqq \sup_{S,T\subseteq\Omega}\left|\int_{S\times T}X(x,y) \diffsymb \mu(x) \diffsymb \mu(y) \right|\;.\] 
The \emph{cut distance} between bigraphons $U,W\in\mathcal{W}_0$ is defined as
\[\deltaCut(U,W)\coloneqq \inf_{\vphi \in S_{\Omega}} \cutn{U-W^{\vphi}}.\]

\subsubsection{Degrees in a bigraphon}
Suppose that $U:\Omega^2\rightarrow [0,1]$ is a bigraphon. Given $x\in\Omega$, we define the \emph{degree of $x$} in $U$ by $\deg_{U}(x)\coloneqq\int_\Omega U(x,y)\diffsymb y$. (So, we choose the convention of integrating over the second parameter.) We prove that the degree distribution in a cut distance convergent sequence of graphons is continuous. While this fact is not difficult, we were not able to find an explicit reference in the literature.
\begin{lem}\label{lem:degsim}
	Suppose that $U$ and $W$ are two bigraphons with $\cutn{U-W}<\eps$. Then the measure of $x\in \Omega$ for which $|\deg_U(x)-\deg_W(x)|>\sqrt{\eps}$ is at most $2\sqrt{\eps}$.
\end{lem}
\begin{proof}
	Let $A\coloneqq\{x\in\Omega\mid\deg_U(x)-\deg_W(x)>\sqrt{\eps}\}$ and $B\coloneqq\{x\in\Omega\mid\deg_W(x)-\deg_U(x)>\sqrt{\eps}\}$. We have
	\[
	\int_{A\times \Omega}U(x,y)\diffsymb (x,y)
	=\int_A \deg_U(x) \diffsymb x
	\ge \int_A \deg_W(x) \diffsymb x +\mu(A)\sqrt{\eps}
	=\int_{A\times \Omega}W(x,y)\diffsymb (x,y)+\mu(A)\sqrt{\eps}\;.
	\]
	As $\cutn{U-W}<\eps$, we in particular get $\int_{A\times \Omega}U(x,y)\diffsymb (x,y)<\int_{A\times \Omega}W(x,y)\diffsymb (x,y)+\eps$. We conclude that $\mu(A)<\sqrt{\eps}$. Similarly, we can get that $\mu(B)<\sqrt{\eps}$.
\end{proof}
From this result, we immediately obtain the following.
\begin{lem}\label{lem:degditcont}
	Suppose that $W_1,W_2,\ldots$ is a sequence of bigraphons converging to a bigraphon $W$ in the cut distance. Then for each open interval $I$, we have
	\[\lim_n \mu(\{x\in\Omega\mid\deg_{W_n}(x)\in I\})=\mu(\{x\in\Omega\mid\deg_{W}(x)\in I\})\;.\]
\end{lem}

\subsection{Permutons}\label{ssub:permutons}
We need to recall the notion of permutons, which is actually used as a building block in our definition of Latinons. For us, a \emph{permutation of size $n$} is a bijection $\pi:[n]\rightarrow[n]$. Suppose that $\rho$ is a bijection of size $k$ for some $k\le n$. Then, the \emph{density} of $\rho$ in $\pi$, denoted by $t(\rho,\pi)$ is the probability that when picking a random $k$-set $K$ of $[n]$, the \emph{restriction of $\pi$ on $K$ is $\rho$}. By that we mean that for any $a,b\in [k]$, and for the corresponding $a$-th smallest element $k_a$ of $K$ and the $b$-th smallest element $k_b$ of $K$, we have $\rho(a)<\rho(b)$ if and only if $\pi(k_a)<\pi(k_b)$. A \emph{permuton} is a probability measure on $[0,1]^2$ with uniform marginals on both axes. The reader can find more information about permutons in~\cite{KrPi} and~\cite{HoKo13}. Yet, we find it instructive to recall how each finite permutation can be represented as a permuton, and also that each sequence of permutations contains a subsequence `converging' to a permuton.\footnote{We explain later which notion of convergence is suitable for permutons.} So, suppose that $\pi:[n]\rightarrow[n]$ is a bijection 
and consider the partition $\partition{n}$.
Then we can represent $\pi$ naturally as a permuton by defining it as an $n$-multiple of the Lebesgue measure on each square 
$\interval{n}{j} \times \interval{n}{\pi(j)}$ 
and zero elsewhere. Next, suppose that $(\pi_n:[n]\rightarrow[n])_n$ is a sequence of permutations and that $(\nu_n)_n$ are the associated permutons. Then by the compactness of the weak topology, there exists a subsequence $\nu_{n_1},\nu_{n_2},\ldots$ to which we have a weak limit, say $\nu$. It is easy to check that such a measure $\nu$ must be a probability measure with uniform marginals, and hence a permuton. In this case, we say that $\nu$ is a limit of $\pi_{n_1},\pi_{n_2},\ldots$. The reason why the weak convergence is the right concept is that it preserves densities. More precisely, $\nu_{n_1},\nu_{n_2},\ldots$ weak converges to $\nu$, if and only if for each fixed permutation $\rho$, we have 
$\lim_i t(\rho,\pi_{n_i})=t(\rho,\nu)$. Here, the density of a permutation $\rho$ of size $k$ into a permuton $\nu$ is defined in the following way. Sample a $k$-tuple of independent points $(x_1,y_1),\ldots,(x_k,y_k)\in[0,1]^2$ according to the law of $\nu$. Let $t(\rho,\nu)$ be the probability that for each $a,b\in [k]$ and for the corresponding $a$-th smallest element $x_{a^*}$ of $\{x_1,\ldots,x_k\}$ and the $b$-th smallest element $x_{b^*}$ of $\{x_1,\ldots,x_k\}$, we have $\rho(a)<\rho(b)$ if and only if $y_{a^*}<y_{b^*}$.

\section{Latinons and the statements of the main results}\label{sec:latinondensitycomp}
In this section, we introduce the new key definitions.
Recall that for us, a \emph{Latin square of order $n$} is a matrix $L\in [n]^{n\times n}$ such that $L_{i,j}\neq L_{i,k}$ and $L_{j,i}\neq L_{k,i}$ for all $1\leq i\leq n$ and $1\leq j< k\leq n$.
Our first goal is to generalise the density notion of permutations to two dimensions.

Let $A\in \mathbb R^{k\times \ell}$ be a matrix and let $O_<(A)=(a_1,\dots,a_r)$ be the increasing ordering of $E(A)=\{a_{i,j}\mid 1\leq i\leq k,1\leq j\leq\ell\}$ according to the natural order $<$. We denote by $p_A:E(A)\rightarrow [r]$ the function which assigns to each $a_{i,j}$ its index in $O_<(A)$. For example for $A=\begin{pmatrix}5 & 6\\
	4 & 8
\end{pmatrix}$ we have $p_A(4)=1,p_A(5)=2,p_A(6)=3,p_A(8)=4$ and for $A'=\begin{pmatrix}5 & 5\\
	4 & 8
\end{pmatrix}$ we have $p_{A'}(4)=1,p_{A'}(5)=2,p_{A'}(8)=3$.

\begin{defi}[Structural equivalence]\label{def:strucequiv}\mbox{}
	\begin{enumerate}[label={(\roman*)}]
		\item We say that two matrices $A,B\in \mathbb R^{k\times \ell}$ are \emph{structurally equivalent}, written $A\equiv B$, if $p_A(a_{i,j})=p_B(b_{i,j})$ for all $1\leq i\leq k$ and $1\leq j\leq \ell$. Matrices of different orders are never considered structurally equivalent.
		
		\item For a matrix $A\in \mathbb R^{k\times \ell}$ and a set $S\subseteq\mathbb R$ we define
		\[\mathcal{R}^A(S)\coloneqq\{M\in S^{k\times\ell}\mid M\equiv A\}\;.\]
	\end{enumerate}
\end{defi}

\begin{remark}\label{rmk:equiv}
	Note that given $k$ and $\ell$, there are only finitely many structural equivalence classes of matrices in $\mathbb R^{k\times\ell}$. 
\end{remark}
We will be interested in matrices $A$ in $[k\ell]^{k\times\ell}$ such that $|E(A)|=k\ell$. In particular, this means that each value from $[k\ell]$ appears exactly once in $A$. We denote this set of matrices by $\mathcal{R}(k,\ell)$ and refer to these as \emph{patterns}. So, the right way of thinking of a pattern is that it is a matrix with a unique smallest entry, a unique second smallest entry etc. By structural equivalence we do not really worry about the actual values, but the relative position of those entries. Given $A \in \mathcal{R}(k,\ell)$ and $B \in [0,1]^{k \times \ell}$ such that $B \equiv A$, we call $B$ a \emph{realisation} of $A$.

Let $L$ be a Latin square of order $n$ and $I,J\subseteq[n]$. We denote by $L|_{(I,J)}$ the submatrix which results from $L$ after deleting all lines with indices $[n]\setminus I$ and all columns with indices $[n]\setminus J$. 
This notation allows us to introduce the key notion of densities in a finite Latin square, very much in parallel to the notion of densities in a finite permutation.
Denote by $\mathbf{x}\in[n]^k_<$ the set of all $\mathbf{x}\in[n]^k$ such that if $\mathbf{x}=(x_1,\dots,x_k)$ then $x_1 < x_2 < \dots < x_k$.
\begin{defi}[Densities in Latin squares]\label{def:latinsqdensities}
	Let $L$ be a Latin square of order $n$ and $A\in \mathcal{R}(k,\ell)$ for some $k,\ell\in[n]$. We define the \emph{density} of the \emph{pattern} $A$ in $L$ by
	\[t(A,L)\coloneqq\frac{1}{\binom{n}{k}\binom{n}{\ell}}\sum_{\mathbf{x}\in[n]^k_<,\mathbf{y}\in[n]^\ell_<}\mathbbm{1}_{\mathcal{R}^A([n])}(L_n|_{(\mathbf{x},\mathbf{y})})\;.\]
\end{defi}
\begin{remark}\label{rem:chocolatenotsublatin}
	Given any sequence $(L_n)_n$ of Latin squares of growing orders, observe that if $A$ is a matrix in $[k \ell]^{k \times \ell}$ such that there is a repeated entry (i.e. $A \not \in \mathcal{R}(k,\ell))$, then the sum in Definition~\ref{def:latinsqdensities} is at most $O(n^{k+\ell-1})$ and so $t(A,L_n) \rightarrow 0$ as $n \rightarrow \infty$. 
	
	So, (when $k=\ell$) restricting a large Latin square to a random $k\times k$ submatrix, we typically get a structure in which each entry appears only once, and such a structure in particular does not correspond to any Latin square. This is in contrast with limits of dense graphs and limits of permutations, where in both cases by restricting a large structure to a random sample of size $k$ we get a substructure which we can interpret in that class (i.e., a graph or a permutation). This is also related to the fact that we do not think that there is a sensible notion of an `$L$-random Latin square' (for a Latinon $L$) while there is a notion of a `$W$-random graph' (for a graphon $W$) and of a `$\nu$-random permutation' (for a permuton $\nu$, see Section~\ref{ssub:permutons}).
	
	In particular, we obtain for each fixed $k,\ell\in \mathbb{N}$ that
	\begin{equation}\label{eq:sum1}
		\lim_{n \rightarrow \infty}   \sum_{A\in \mathcal{R}(k,\ell)} t(A,L_n)=1\;.
	\end{equation}
\end{remark}

We are now ready to define our limit objects, Latinons.  
As these definitions are rather technical, in Section~\ref{ssec:ExamplesOfLatinons}, we give some basic examples of sequences of finite Latin squares and the corresponding limit object. These examples in particular explain why our Latinons are distribution-valued, and also the role of the function $f$. For this recall that we denote by $\mathcal{B}$ the space of all Borel probability measures on $[0,1]$.

\begin{defi}[Semilatinon and Latinon]\label{def:latinon}
	A \emph{semilatinon} is a pair $(W,f)$ such that $W:\Omega^2\rightarrow\mathcal{B}$ is a distribution-valued bigraphon and $f:\Omega\rightarrow[0,1]$ is measure preserving function. A \emph{Latinon} is a semilatinon which satisfies the following properties.
	\begin{enumerate}[label={(\roman*)}]
		\item\label{en:CondRow} The measure $\mu_{W,x}^1$ on $[0,1]^2$ defined by $\mu_{W,x}^1(S\times T)\coloneqq\int_{y\in f^{-1}(S)}W(x,y)(T) \diffsymb y$ for all measurable $S,T\subseteq[0,1]$ is a permuton for almost all $x\in\Omega$.\footnote{$\mu_{W,x}^1$ is defined here only on sets of the form $S\times T$. This is sufficient as such sets generate the Borel sigma-algebra on $[0,1]^2$.}
		\item\label{en:CondColumn} The measure $\mu_{W,y}^2$ on $[0,1]^2$ defined by $\mu_{W,y}^2(S\times T)\coloneqq\int_{x\in f^{-1}(S)}W(x,y)(T) \diffsymb x$ for all measurable $S,T\subseteq[0,1]$ is a permuton for almost all $y\in\Omega$.
	\end{enumerate}
	We denote the space of all Latinons by $\LatinonSpace$.
\end{defi}

Condition~\ref{en:CondRow} (respectively~\ref{en:CondColumn}) are limit counterparts to the fact that given a Latin square, we can take its slice at any given row (resp.\ column) and get a permutation.

\begin{remark}
	Note that the conditions of Definition~\ref{def:latinon} ensure for a Latinon $(W,f)$ that there is no set $X\subset \Omega$ of positive measure such that $W|_{X\times \Omega}$ or $W|_{\Omega\times X}$ is identically equal to $\Dirac_p$ for some $p\in[0,1]$.
\end{remark}

We identify Latinons if they agree everywhere except a nullset. So we are actually working with the quotient space under this relation, but we will for convenience still speak of Latinons instead of equivalence classes.

Given a function $f:\Omega\rightarrow[0,1]$ we can introduce a strict partial order $<_f$ on $\Omega$ by defining $x<_f y$ if and only if $f(x)< f(y)$. We write $\mathbf{x} \in \Omega^k_{<_f}$ to denote $\mathbf{x} = (x_1,\dots,x_k) \in \Omega^k$ such that $\mathbf{x}$ is increasingly ordered according to $<_f$, i.e. $f(x_1) < f(x_2) < \dots < f(x_k)$.
Furthermore for measures $\mu_{(1,1)},\dots,\mu_{(k,\ell)}$ on $[0,1]$ we write $\bigotimes_{(i,j)\in[k]\times[\ell]}\mu_{(i,j)}$ for the corresponding product measure on $[0,1]^{k\times\ell}$. Now we can extend our notion of density to Latinons.

\begin{defi}[Densities in Latinons]\label{def:latinondensities}
	Let $(W,f)$ be a Latinon over the ground space $\Omega$ and $A\in \mathcal R(k,\ell)$. We denote by $t(A,(W,f))$ the \emph{density} of the \emph{pattern} $A$ in $(W,f)$ and define it to be
	\[t(A,(W,f))\coloneqq k!\ell!\int_{\mathbf{x}\in\Omega^k_{<_f}}\int_{\mathbf{y}\in\Omega^\ell_{<_f}}\left(\bigotimes_{(i,j)\in[k]\times[\ell]}W(x_i,y_j)\right)(\mathcal{R}^A([0,1]))\diffsymb \mathbf{y}\diffsymb \mathbf{x}\;.\]
\end{defi}

\begin{defi}[Sampling from a Latinon]\label{def:latinonsampling}
	Let $(W,f)$ be a Latinon over the ground space $\Omega$ and $A\in \mathcal R(k,\ell)$.
	We say that $B\in [0,1]^{k\times \ell}$ is \emph{sampled} from $(W,f)$ if $B$ is a 
	matrix sampled according to the following random selection. We repeatedly choose $\mathbf{x}=(x_1,\dots,x_k)\in\Omega^k$ and $\mathbf{y}=(y_1,\dots,y_\ell)\in\Omega^\ell$ independently and uniformly at random until we find $\mathbf{x}$ and $\mathbf{y}$ which are increasingly ordered according to $<_f$. For each pair $(i,j)$ with $1\leq i\leq k$ and $1\leq j\leq \ell$ we sample a real value $b_{i,j}$ from the distribution $W(x_i,y_j)$ and set $B=(b_{i,j})$.
\end{defi}

Informally, we see that $t(A,(W,f))$ is the probability that $B\equiv A$, where $B$ is sampled from $(W,f)$. Also note that according to our earlier definitions $B$ is a realisation of the pattern $A$. The sampling view gives us the following counterpart to~\eqref{eq:sum1}.
\begin{fact}
	Suppose that $(W,f)$ is a Latinon. Then $\sum_{A\in \mathcal{R}(k,\ell)} t(A,L)=1$.
\end{fact}

The following definition will help us to associate a finite Latin square with a Latinon. The similarities to associating a permuton to a finite permutation as was done in Section~\ref{ssub:permutons} are obvious.

\begin{defi}\label{def:latinonrep}
	Let $L_n$ be a Latin square of order $n$. Given a measure preserving function $f:\Omega\rightarrow[0,1]$ we derive a Latinon $(W_{L_n},f)$ over the ground space $\Omega$ from $L_n$ by setting
	\[\left(W_{L_n}(x,y)\right)(S)\coloneqq n\cdot\lambda\left(S\cap \left[\frac{L_{nxy}-1}n,\frac{L_{nxy}}n\right]\right),\textrm{ for }x,y\in\Omega\textrm{ and }S\subseteq[0,1]\;,\]
	where $L_{nxy}\coloneqq L_n(\lceil n\cdot f(x) \rceil,\lceil n \cdot f(y) \rceil )$.
\end{defi}

Note that this indeed defines a Latinon and in fact Latinon-representations of Latin squares have almost the same densities.

\begin{prop}\label{prop:latinonrep}
	Let $L_n$ be a Latin square of order $n$ and $f:\Omega\rightarrow[0,1]$ a measure preserving function. Then for every $k,\ell \in \mathbb{N}$ and every $k \times \ell$ pattern $A$ we have
	\[|t(A,(W_{L_n},f))-t(A,L_n)|\leq\tfrac{2(k+\ell)^2}{n}.\]
\end{prop}

\begin{proof}
	Consider the probability that $B\equiv A$, where $B\in [0,1]^{k\times \ell}$ is a matrix sampled according to the following random selection. Sample as one would for $t(A,(W_{L_n},f))$ in Definition~\ref{def:latinondensities}, except rather than just sampling so that $\mathbf{x}$ and $\mathbf{y}$ are increasingly ordered according to $<_f$, but also so that for each $i,j$ we have $\lceil n\cdot f(x_i) \rceil \not= \lceil n\cdot f(x_j) \rceil$ (and similarly for the $y_i$). The crucial observation is that this probability is equal to $t(A,L_n)$ by definition of $(W_{L_n},f)$. So define $S_1$ and $S_2$ to be subsets of $[0,1]^{k+\ell}$ such that 
	\begin{align*}
		S_1&\coloneqq\left\{ \mathbf{x}\in[0,1]^k_{<_f}, \mathbf{y}\in[0,1]^{\ell}_{<_f} \mid \lceil n\cdot f(x_i) \rceil \not= \lceil n\cdot f(x_j) \rceil, \lceil n\cdot f(y_i) \rceil \not= \lceil n\cdot f(y_j) \rceil \text{ for } i\not=j  \right\}\;,\\
		S_2&\coloneqq\left\{ \mathbf{x}\in[0,1]^k_{<_f}, \mathbf{y}\in[0,1]^{\ell}_{<_f} \right\} \setminus S_1\;.
	\end{align*}
	We thus obtain
	\[t(A,(W_{L_n},f)) = k! \ell! \left( \int_{S_1} t(A,L_n) \diffsymb (\mathbf{x},\mathbf{y})+ \int_{S_2} \left(\bigotimes_{(i,j)\in[k]\times[\ell]}W(x_i,y_j)\right)(\mathcal{R}^A([0,1])) \diffsymb (\mathbf{x},\mathbf{y}) \right),\]
	and also by noting that $\lambda^{\otimes k+\ell}(S_1) + \lambda^{\otimes k+\ell}(S_2)=1/(k! \ell!)$, we have
	\begin{align*}
		|t(A,(W_{L_n},f))-t(A,L_n)| & \leq | k! \ell! \lambda^{\otimes k+\ell}(S_1) - 1 | \cdot t(A,L_n) + k! \ell! \lambda^{\otimes k+\ell}(S_2) \leq 2k! \ell! \lambda^{\otimes k+\ell}(S_2) \\ 
		& = 2 \left( 1 - \frac{(n-1) \cdots (n-k+1)}{n^{k-1}} \cdot \frac{(n-1) \cdots (n-\ell+1)}{n^{\ell-1}} \right) \\
		& \leq 2 \left( \frac{n^{k+\ell -2} - (n-k-\ell+2)^{k+\ell-2}}{n^{k+\ell-2}} \right) \leq \frac{2(k+\ell)^2}{n}.
	\end{align*} 
\end{proof}

We can now define a notion of left-convergence for Latin squares and Latinons. 
\begin{defi}[Left-convergence]\label{def:left}
	\begin{enumerate}[label={(\roman*)}]
		\item\label{ew1} Let $(L_n)_{n\in\mathbb N}$ be a sequence of Latin squares and let $(W,f)$ be a Latinon. We say that $(W,f)$ is the limit of $(L_n)_{n\in\mathbb N}$, written $L_n\LeftConvergence (W,f)$, if\\ $\lim_{n\rightarrow\infty}t(A,L_n)=t(A,(W,f))$ for every $k,\ell\in \mathbb N$ and $A\in\mathcal{R}(k,\ell)$.
		\item\label{ew2} Let $(L_n)_{n\in\mathbb N}$ be a sequence of Latinons and let $L$ be a Latinon. We say that $L$ is the limit of $(L_n)_{n\in\mathbb N}$, written $L_n\LeftConvergence L$, if $\lim_{n\rightarrow\infty}t(A,L_n)=t(A,L)$ for every $k,\ell\in \mathbb N$ and $A\in\mathcal{R}(k,\ell)$.
	\end{enumerate}
\end{defi}
Note that by Definition~\ref{def:latinonrep} and Proposition~\ref{prop:latinonrep} we have that a sequence of Latin squares (of growing orders) converges in the sense of Definition~\ref{def:left}\ref{ew1} if and only if their Latinon-representations converge in the sense of Definition~\ref{def:left}\ref{ew2}.

Note that we require the convergence only for countably many sequences of densities, but by Remark~\ref{rmk:equiv} this implies convergence for the densities of all $A\in\mathbb R^{k\times\ell}$.
Note also, that by a standard abstract compactness argument, if $(L_n)_{n\in\mathbb{N}}$ is a sequence of Latin squares of growing orders, then there exists a subsequence $(L_{n_i})_i$ so that for every $k,\ell\in \mathbb N$ and $A\in\mathcal{R}(k,\ell)$ we have that $\lim_{i\rightarrow\infty}t(A,L_{n_i})$ exists. A similar type of left-convergence is at the heart of flag algebras,~\cite{Razborov2007}. However, here we aim to obtain Latinons as a more explicit representation of the limit. One of the main results of this paper is the following compactness result. 

\begin{thm}[Compactness for Latinons]\label{thm:compact}
	Let $(L_n)_{n\in\mathbb N}$ be a sequence of Latinons over arbitrary ground spaces. There exists a subsequence $(L_{n_i})_{i\in\mathbb N}$ and a Latinon $(W,f)$ over the ground space $\Omega$ such that
	\[L_{n_i}\LeftConvergence (W,f)\;.\]
\end{thm}

It is worth recalling that in Section~\ref{ssub:permutons} we sketched an analogous compactness result for permutons, which follows immediately from the compactness of the weak topology of probability measures on $[0,1]^2$. 
Although we can represent a Latinon by a probability measure (see Definition~\ref{def:Latinonasmeasure}), the proof of Theorem~\ref{thm:compact} is much more complicated. 
The reason for this is that it is not enough to view coordinates of different values with global lenses, as our example in Section~\ref{sss:regularrandom} shows.
(That is, the probability measure needs to be on $\Omega^2 \times [0,1]$ rather than $[0,1]^3$.)

Recall that by Proposition~\ref{prop:latinonrep} this applies to sequences of Latin squares as well.
Theorem~\ref{thm:compact} can be accompanied by the following approximation result.
\begin{thm}\label{thm:latinonsarelimits}
	For every Latinon $(W,f)\in\LatinonSpace$ there exists a sequence $(L_n)_{n\in\mathbb N}$ of Latin squares such that $L_n\LeftConvergence (W,f)$.
\end{thm}

\subsection{Examples}\label{ssec:ExamplesOfLatinons}
For this section only, for convenience, we will assume that all finite $n \times n$ Latin squares have row numbers, column numbers and values in $[0,n-1]$ (rather than $[n]$).
\subsubsection{Standard cyclic Latin squares}
\label{ssec:scls}
One of the most natural examples of a sequence of Latin squares comes from the Cayley tables of cyclic groups. For each $n \in \mathbb{N}$ define the Latin square $L_n$ by $L_n(x,y)\coloneqq x+y \! \mod n$ for all $x,y \in [0,n-1]$. We define the \emph{standard cyclic Latinon} $(L,f)$ as a Latinon on $\Omega=[0,1]$ by letting $L:[0,1]^2 \to \mathcal{B}$ be defined by $L(x,y)\coloneqq \Dirac_{x+y\textrm{ mod }1}$ for all $x,y \in [0,1]$, and by setting $f:[0,1] \to [0,1]$ to be the identity. 

\begin{figure}[h]
	\begin{center}
		\begin{minipage}[t]{.25\linewidth}
			\vspace{0pt}
			\centering
			{\begin{tabular}{ c c c c c c }
					0 & 1 & 2 & 3 & 4 & 5 \\
					1 & 2 & 3 & 4 & 5 & 0 \\
					2 & 3 & 4 & 5 & 0 & 1 \\
					3 & 4 & 5 & 0 & 1 & 2 \\
					4 & 5 & 0 & 1 & 2 & 3 \\
					5 & 0 & 1 & 2 & 3 & 4 \\
			\end{tabular}}
		\end{minipage}
		\begin{minipage}[t]{.35\linewidth}
	\vspace{0pt}
	\centering
	\includegraphics[scale=0.32,clip]{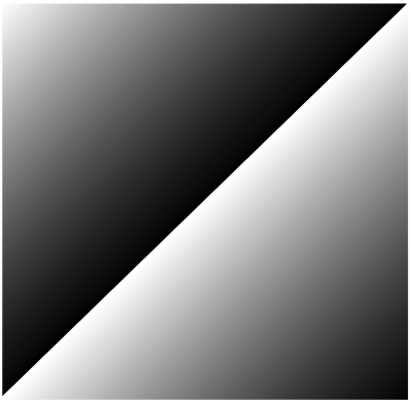}
\end{minipage}
	\end{center}
	\caption{The Latin square $L_6$ and the standard cyclic Latinon $(L,f)$ in which one should think of $x$ from $0$ to $1$ as running from left to right, $y$ from $0$ to $1$ as running from top to bottom, and within the image, the entries are scaled from white being $0$ to black being $1$.}
	\label{fig:cyclic}
\end{figure}
It can easily be shown that $(L,f)$ is a Latinon, and in Section~\ref{ssec:cyclicproof} we will show that $(L_n)_n$ converges to $(L,f)$ as $n$ tends to infinity. As mentioned above $L_n$ is a Cayley table of a cyclic group of order $n$. By thinking of $x,y$ and $x+y$ as points on the unit circle, one can also consider $(L,f)$ to be the Cayley table of the group $\mathbb{R}/\mathbb{Z}$. 

\subsubsection{Regular versus random alternation}\label{sss:regularrandom}
Consider the partition $\partition{n}$ which we defined earlier as the equipartition of $[0,1]$ into $n$ intervals.
Let $C_n:[0,1]^2\rightarrow[0,1]$ be the `chessboard graphon', that is $C_n$ is either constant $0$ or $1$ on each 
$\interval{n}{i} \times \interval{n}{j}$
depending on the parity of $i+j$. Let $R_n$ be a graphon in which a constant $0$ or $1$ on each 
$\interval{n}{i} \times \interval{n}{j}$
($i\le j$) is chosen uniformly at random. Recall that $R_n$ can be thought of as the adjacency matrix of an Erd\H{o}s--R\'enyi random graph $\mathbb{G}(n,\frac12)$ (except the values on the diagonal). Thus, basic theory of graphons tells us that $(C_n)_n$ converges to the complete balanced bipartite graphon while $(R_n)_n$ converges to the constant-one-half graphon. Our examples below can be thought of as a Latin square counterpart to this difference between `regular fifty-fifty' versus `random fifty-fifty'.

For each $n$, define an $n\times n$ matrix $H_n$,
\begin{align*}
	H_n(i,j) & \coloneqq \begin{cases}
		i+j \! \mod n & \text{if } i+j \equiv 0 \! \mod 2, \\
		-i-j \! \mod n & \text{if } i+j \equiv 1 \! \mod 2.
	\end{cases}
\end{align*}
It is easy to check that for $n$ even, $H_n$ is actually a Latin square.

Let us now consider a matrix $P_n$ which should be informally thought of in the following way:
\begin{align*}
	P_n(i,j) & \coloneqq \begin{cases}
		i+j \! \mod n & \text{with probability 1/2,} \\
		-i-j \! \mod n & \text{with probability 1/2.}
	\end{cases}
\end{align*}
One could say that $P_n$ is a Latin square `in expectation'. By that we mean that for each $k\in\{0,\ldots,n-1\}$, on each row $i$, there is exactly one column $j$ so that $k=i+j \! \mod n$ and exactly one column $j$ so that $k=-i-j \! \mod n$, indeed giving that the number of times the symbol $k$ appears is~1 in expectation. As it turns out, it is indeed possible to construct deterministic Latin squares which behave like $P_n$. We omit the details here since this construction is described in Section~\ref{sec:ApproximatingLatinons}. Now, $(P_n)_n$ exhibits the problem that $\lim_{n\rightarrow\infty} P_n(xn,yn)/n$ does not exist for $x,y\in[0,1]$. We expect the two accumulation points to be $x+y$ mod $1$ and $-x-y$ mod $1$. A solution to this problem is to let the limit object be a distribution-valued bigraphon. A candidate for the limit of $(P_n)_n$ could therefore be a distribution-valued bigraphon $P$ defined on $[0,1]^2$, where one copy of $[0,1]$ would represent rows and the other copy would represent columns, 
\[P(x,y)\coloneqq \frac12\Dirac_{(x+y\textrm{ mod }1)}+\frac12\Dirac_{(-x-y\textrm{ mod }1)}\;.\]
However, the sequence $(H_n)_n$ shows that such a definition is still not general enough. In the limit and for values $(x,y)\in[0,1]^2$ we cannot distinguish the odd from the even case. Since in a square $[(x-\eps )n,(x+\eps )n]\times[(y-\eps)n,(y+\eps)n]$ half of the time the value $H_n(i,j)$ will be $i+j \! \mod n$ and half of the time $-i-j \! \mod n$, one could guess that $H_n$ has the same limit as $P_n$. This is however not correct, as the densities of substructures in $(H_n)_n$ and $(P_n)_n$ converge to different limits; one can check that for the $2\times 3$ pattern 
\[A\coloneqq\begin{pmatrix}1 & 2 & 3\\
	4 & 5 & 6
\end{pmatrix}\]
we have that $\lim_{n\rightarrow\infty}t(A,H_n)>\lim_{n\rightarrow\infty}t(A,P_n)$. We therefore need to be able to encode `local information' in the limit object, which suggests to consider a limit object like $H:\Omega^2\rightarrow \mathcal{B}$, where in this case $\Omega=[0,1]\times\{\textrm{odd},\textrm{even}\}$ and
\[H\left((x,\eta_x),(y,\eta_y)\right)=\begin{cases}
	x+y \mod 1& \text{if $\eta_x=\eta_y$,}\\
	-x-y \mod 1 & \text{if  $\eta_x\neq\eta_y$.}
\end{cases}\]

The Latinons for the above two sequences can be described as follows. For the sequence $(H_n)$, we set $(H,h)$ to be the Latinon on $\Omega=[0,1]\times \{\text{odd, even}\}$, where $H:([0,1] \times \{\text{odd, even}\})^2 \to \mathcal{B}$ and $h:[0,1] \times \{\text{odd, even}\} \to [0,1]$ are defined by
$$h(x,a) \coloneqq x\;,
\quad
H((x,a),(y,b))\coloneqq \begin{cases} 
	x+y \! \mod 1 & \text{if } a=b, \\
	-x-y \! \mod 1 & \text{if } a\not=b.
\end{cases}
$$
For the sequence $(P_n)$, we set $(P,p)$ to be the Latinon on $\Omega=[0,1]$, where $p$ is the identity and $P:[0,1]^2 \to \mathcal{B}$ is defined by
$$
P(x,y) : = \begin{cases}
	x+y \! \mod 1 & \text{with probability 1/2,} \\
	-x-y \! \mod 1 & \text{with probability 1/2.}
\end{cases}
$$

\medskip

The sequence $(H_n)_n$, which shows that in addition to the `global position' of each row/column, some local information (`parity-like') must be recorded as well, nicely illustrates the main difference between the theory of Latinons and that of permutons. Indeed, recall that permutons record only the global position by being a measure on $[0,1]^2$. Also, recall the proof of the compactness result for permutons, as sketched in Section~\ref{ssub:permutons}. The approach using the weak topology, which is the key tool in the proof, can only measure the behaviour of the finite permutations on a global scale, that is, it tracks quantities such as `how many elements around $0.4n$ are mapped by $\pi_n$ to a value of around $0.9n$'. In other words, one cannot create any parity-like-based constructions of finite permutations that would drastically change subpermutation densities.

\subsubsection{Very Local Cyclic Latin squares}
\label{ssec:verylocal}
Let us give a sequence of Latin squares which are constructed in a similar way to standard cyclic Latin squares, but change their values in a more discontinuous fashion. For every $n=k^2$ for some $k \in \mathbb{N}$, define the Latin square $C_n$ by setting $$C_n((i-1)k+x,(j-1)k+y)\coloneqq i+j+(x+y-2)k \; \mod n$$ for every $i,j,x,y \in [k]$. Now by considering the partition of such a Latin square into $k \times k$ blocks $B_{i,j}$, $i,j \in [k]$, we see that the blocks differ very little from each other (since the role of $i,j$ is rather inferior in the definition of $C_n$) whereas within a block the numbers change much more. One possible limit representation is a Latinon $(C,f)$ on $\Omega=[0,1]\times [0,1]$ with $C((i,x),(j,y))=(x+y) \mod 1$ and $f$ being the identity. We omit a proof, which could be done using the same techniques which we use in Section~\ref{ssec:cyclicproof} to prove that $(L,f)$ is a correct limit to standard cyclic Latin squares.

\subsubsection{Columns and values swap}\label{ssec:swapMakesDifference}
Recall that swapping row numbers and column numbers, row numbers and values, and column numbers and values of a Latin square again results into a Latin square. It is not hard to see by the definition of left-convergence that if a sequence $(L_n)_n$ of Latin squares converges, so does the sequence $(M_n)_n$ with swapped rows and columns, i.e. $M_{n}(x,y)=L_{n}(y,x)$ for all $x,y\in [0,n-1]$. However, the situation presents itself differently in the case of swapping column numbers and values. For a finite $n \times n$ Latin square $L$, we denote by $L'$ the finite Latin square obtained from $L$ by swapping column numbers and values; that is for all $x,y,z \in [0,n-1]$ we have that if $L(x,y)=z$, then $L'(x,z)=y$. In this section we prove the following. 

\begin{lem}
	There exist two sequences of finite Latin squares $(J_n)_{n\in 3\mathbb{N}}$ and $(K_n)_{n\in 3\mathbb{N}}$ so that by interlacing them, we get a left-convergent sequence, yet by interlacing the respective sequences formed by column-value swaps, $(J'_n)_{n\in 3\mathbb{N}}$ and $(K'_n)_{n\in 3\mathbb{N}}$ we do not get a left-convergent sequence.
\end{lem}

\begin{proof}
	Suppose $n$ is divisible by $3$. Let $I_1\coloneqq[0,\frac{n}{3}-1]$, $I_2\coloneqq[\frac{n}{3},\frac{2n}{3}-1]$ and $I_3\coloneqq[\frac{2n}{3},n-1]$. Set $x^*\coloneqq x \! \mod (\frac{n}{3})$ and $y^*\coloneqq y \! \mod (\frac{n}{3})$. We define two Latin squares $J_n$ and $K_n$ by
	\begin{align*}
		J_n(x,y) & \coloneqq \begin{cases} 
			3(x^*+y^*) \! \mod n & \text{if $x,y \in I_1$, or $x \in I_2$ and $y \in I_3$, or $x\in I_3$ and $y \in I_2$;} \\
			3(x^*+y^*+1) \! \mod n & \text{if $x \in I_2$ and $y \in I_1$, or $x \in I_1$ and $y \in I_2$, or $x,y \in I_3$;} \\
			3(x^*+y^*+2) \! \mod n & \text{otherwise;} 
		\end{cases} \\
		K_n(x,y) & \coloneqq \begin{cases} 
			3(x^*+y^*) \! \mod n & \text{if $x,y \in I_1$, or $x,y \in I_2$, or $x,y\in I_3$;} \\
			3(x^*+y^*+1) \! \mod n & \text{if $x \in I_1$ and $y \in I_3$, or $x \in I_2$ and $y \in I_1$, or $x \in I_3$ and $y \in I_2$;} \\
			3(x^*+y^*+2) \! \mod n & \text{otherwise.} 
		\end{cases}
	\end{align*}
	Let $M_n$ be the restriction of either of $J_n$ or $K_n$ to $I_a \times I_b$ for any fixed $a,b \in [3]$. It is easy to see that $(M_n)_{n\in 3\mathbb{N}}$ converges to the standard cyclic Latinon $L$. Since this is true for any $a,b$ for both $J_n$ and $K_n$, it follows that both sequences $(J_n)_{n\in 3\mathbb{N}}$ and $(K_n)_{n\in 3\mathbb{N}}$ converge to the same Latinon. (This Latinon is formed of nine blocks of equal size in a $3 \times 3$ array, which each themselves are a copy of $L$ scaled down by three.)  
	
	Now define $J'_n$ and $K'_n$ to be the respective Latin squares obtained from $J_n$ and $K_n$ by swapping column numbers and values as described above, and let $A$ be the pattern $\begin{pmatrix} 2 & 1 \\ \end{pmatrix}$. A simple calculation shows that
	$$t(A,J'_n)=\frac{9\left(\sum_{i=0}^{\frac{n}{3}-1} \binom{\frac{n}{3}-i}{2}+\binom{i}{2} \right)+2n(\frac{n}{3})^2}{n\binom{n}{2}},$$
	whereas
	$$t(A,K'_n)=\frac{9\left(\sum_{i=0}^{\frac{n}{3}-1} \binom{\frac{n}{3}-i}{2}+\binom{i}{2} \right)+n(\frac{n}{3})^2}{n\binom{n}{2}}.$$ 
	From this we observe that 
	$$\lim_{n \to \infty} |t(A,J'_n) - t(A,K'_n)| = \lim_{n \to \infty} \frac{n(\frac{n}{3})^2}{n\binom{n}{2}}=\frac{2}{9},$$
	and thus $(J'_n)_{n \in 3\mathbb{N}}$ and $(K'_n)_{n \in 3\mathbb{N}}$ do not converge to the same Latinon. 
\end{proof}


\subsection{Two cryptomorphic definitions of Latinons}\label{ssec:Cryptomorphic}
We give two alternative definitions of a Latinon. The first one is as a measure on $\Omega^2\times [0,1]$. In the second one, we do not index the rows (and columns) of a Latinon by $\Omega$ together with the measure preserving function $f:\Omega\rightarrow[0,1]$ but rather by $[0,1]\times \Gamma$ (and without any measure preserving function), where $\Gamma$ is a standard probability space.
\begin{defi}[Latinons as measures]\label{def:Latinonasmeasure}
	A \emph{Latinon in measure representation} is a pair $(\nu,f)$ where $f:\Omega\rightarrow [0,1]$ is a measure preserving function and $\nu$ is a Borel measure on $\Omega^2\times [0,1]$ which has uniform marginals, that is, for all subsets $X,Y\subseteq\Omega$ and $Z\subseteq[0,1]$ we have $\nu(X\times Y\times [0,1])=\mu(X)\mu(Y)$, $\nu(\Omega\times Y\times Z)=\mu(Y)\lambda(Z)$, and $\nu(X\times \Omega\times Z)=\mu(X)\lambda(Z)$.
\end{defi}
\begin{defi}[Latinons indexed by global and local coordinates]\label{def:Latinongloballocal}
	Suppose that $\Gamma$ is a standard measure space with a probability measure $\gamma$ and $U$ is a function $U:([0,1]\times \Gamma)^2\rightarrow \mathcal{B}$. $U$ is said to be a \emph{Latinon indexed by global and local coordinates} if for almost every $(a,x), (b,y) \in [0,1]\times \Gamma$ the measures $\mu^1_{a,x}$ on $[0,1]^2$ defined by $\mu^1_{a,x}(S\times T)\coloneqq\int_{b\in S}\int_{y\in\Gamma}U\big((a,x),(b,y)\big)(T)\diffsymb y\diffsymb b$ and $\mu^2_{b,y}$ on $[0,1]^2$ defined by $\mu^2_{b,y}(S\times T)\coloneqq\int_{a\in S}\int_{x\in\Gamma}U\big((a,x),(b,y)\big)(T)\diffsymb x \diffsymb a$ are permutons.
\end{defi}

We shall see in the sections below that one can indeed transform a Latinon from one representation to another. It is important that we could translate the notion of density to any of these alternative definitions. To this end, we only need to be able to sample random rows and columns.
\begin{itemize}
	\item In case of a representation as a measure, we use the measure $\mu$.
	\item In case of a representation indexed by global and local coordinates, we use the product measure $\lambda\times\gamma$.
\end{itemize}

\subsubsection{The role of $[0,1]\times \Gamma$ in the global-local representation} We saw in Section~\ref{ssec:ExamplesOfLatinons} that sometimes we need to distinguish different types of rows/columns in a matrix even when these rows/columns are very close together. Because of this necessity $f$ (in Definition~\ref{def:latinon}) need not be injective. Definition~\ref{def:Latinongloballocal} provides an alternative way of recording different types of local behaviour. That is, the component $[0,1]$ records the global position of rows/columns of some large Latin square while the coordinate $\Gamma$ is used to store different types of rows/columns with a similar global position. Sometimes the local coordinate is irrelevant, e.g. in the standard cyclic Latinon from Section~\ref{ssec:scls}, whereas for $(H,h)$ from Section~\ref{sss:regularrandom} the space $\Gamma$ would be partitioned into two pieces of measure half corresponding to odd and even.\footnote{Note that we do not have a universal partition of $\Gamma$ into the odd and even part; the partition may change as we vary the global coordinate.} The example from Section~\ref{ssec:verylocal} shows a limit graphon whose values need not depend on the global coordinate but are highly nonconstant on the local coordinate.

\subsubsection{From a standard representation to a representation by a measure} 

Suppose that $(W,f)$ is a standard representation of a Latinon. Then define $\nu$ by setting $\nu(X\times Y\times T)\coloneqq\int_{x\in X}\int_{y\in Y} W(x,y)(T) \diffsymb y \diffsymb x$ for each $X,Y\subset \Omega$ and $T\subset [0,1]$. Clearly, the uniform marginal condition is satisfied.

\subsubsection{From a representation by a measure to a standard representation} Suppose that $(\nu,f)$ is a Latinon as in Definition~\ref{def:Latinonasmeasure}. Now, it is enough to define $\{W(x,y)\}_{x,y\in\Omega}$ to be the disintegration of $\nu$ according to $\Omega^2$. It is straightforward to check that Latinons are transformed correctly in this construction.

\subsubsection{From a standard representation to a representation with global and local coordinates}\label{ssec:fromStandardToLocGlob}
The key to the whole transformation is Proposition~\ref{prop:transformIntoLocalGlobal} below. Proposition~\ref{prop:transformIntoLocalGlobal} looks like a standard result in real analysis. However, we were not able to find it anywhere. We thank Jan Greb\'ik for suggesting its proof, and Martin Dole\v zal for discussions about its subtleties. The proof is included in the appendix.
\begin{prop}\label{prop:transformIntoLocalGlobal}
	Suppose that $(\Omega,\mu)$ is a separable atomless probability measure space, and that $f:\Omega\rightarrow [0,1]$ is a measure preserving function. Then there exists a measure preserving function $h:[0,1]\times [0,1]\rightarrow \Omega$ such that for almost every $(a,x)\in [0,1]\times [0,1]$ we have $f(h(a,x))=a$.
\end{prop}
Given Proposition~\ref{prop:transformIntoLocalGlobal}, it is easy to perform the desired transformation of the domain $\Omega^2$ to $([0,1]^2)^2$. Indeed, given a Latinon $(M,f)$ in standard representation, we take the function $h:[0,1]^2\rightarrow \Omega$ from Proposition~\ref{prop:transformIntoLocalGlobal} and define $U:([0,1]^2)^2\rightarrow \mathcal{B}$ by $U((a,x),(b,y))\coloneqq M(h(a,x),h(b,y))$.

\subsubsection{From a representation with global and local coordinates to a standard representation}
Suppose that $U:([0,1]\times \Gamma)^2\rightarrow \mathcal{B}$ is a Latinon as in Definition~\ref{def:Latinongloballocal}. Define a function $g:[0,1]\times \Gamma\rightarrow [0,1]$ by $(a,x)\mapsto a$. Clearly, this function is measure preserving. If $\Omega=[0,1]\times \Gamma$ (as a measure space), then $(U,g)$ is obviously a Latinon as in Definition~\ref{def:latinon}. In the general case, we first fix a measure preserving bijection $\phi:\Omega\rightarrow [0,1]\times \Gamma$, and set $f:\Omega\rightarrow [0,1]$, $f(z)\coloneqq g(\phi(z))$ and $W:\Omega^2\rightarrow\mathcal{B}$, $W(u,v)\coloneqq U(\phi(u),\phi(v))$. It is straightforward to check that Latinons are transformed correctly in this construction.

\subsection{A symmetric limit theory of Latin squares}\label{ssec:NewDiscussion}
As we explained in Section~\ref{ssec:representationAs3Dim}, each Latin square of order $n$ can be represented using a $3$-uniform tripartite hypergraph on $n+n+n$ vertices satisfying~\ref{L1}--\ref{L3}, and that, in fact, this representation is a bijection. The way we introduce convergence in Definition~\ref{def:left} is \emph{a priori} asymmetric in that the rows and the columns play a different role to the values. Indeed, we saw in Section~\ref{ssec:swapMakesDifference} that swapping columns and values can turn a convergent sequence of Latin square into a divergent one. 

In this section, we briefly discuss a potential limit theory motivated purely by the representation of a Latin square as a $3$-uniform tripartite hypergraph satisfying~\ref{L1}--\ref{L3}. While the extra symmetries may make such a theory seem more appealing in the context of hypergraphs we do not see how to introduce any meaningful notion of densities to derive combinatorial applications. Such a limit theory is also briefly mentioned (in the context of quasirandomness) in~\cite[Section~3.1]{CorRaz:NaturalQuasiramness}. 

So, suppose that $L$ is a Latin square encoded by a 3-uniform hypergraph $H$ as described in Section~\ref{ssec:representationAs3Dim}. For $i,j,k\in[n]$ we write $(i,j,k)\in H$ if $(i,\mathtt{row}), (j,\mathtt{column}), (k,\mathtt{value})$ forms an edge. We can represent $L$ by a Borel measure $\kappa_L$ on $[0,1]^3$ by setting
\[
\kappa_L(B):=n\cdot\lambda^{\otimes 3}\left(B\cap\bigcup_{(i,j,k)\in H}[\tfrac{i-1}{n},\tfrac{i}{n}]\times [\tfrac{j-1}{n},\tfrac{j}{n}]\times [\tfrac{k-1}{n},\tfrac{k}{n}]\right)
\]
for each Borel set $B\subset [0,1]^3$. In particular, $\kappa_L$ is a Borel probability measure and properties~\ref{L1}--\ref{L3} tell us that the marginals on each of the three coordinates are permutons. We call all probability measures on $[0,1]^3$ with the above marginal condition \emph{symmetric-Latinons}. Recall that the space of Borel probability measures on $[0,1]^3$ equipped with the weak topology is compact, and hence so is the space of all symmetric-Latinons. That is, if we take the weak topology as `the cut distance'\footnote{metrized for example by the Lévy--Prokhorov metric} for the general limit framework~\ref{en:F1}--\ref{en:F6} then compactness is guaranteed. To complete this program, we would need a notion of densities (as in~\ref{en:F2} and~\ref{en:F5}) and a sampling lemma (as in~\ref{en:F6}). This could be an opening for future research, but quite possibly no combinatorially meaningful interpretation exists.
%

Note that relying solely on the weak topology is similar to that in the theory of permutation limits,~\cite{HoKo13}.

\section{A cut distance, counting lemmas, and equivalence of convergence}\label{sec:cutdistcountconv}
In this section we introduce a cut distance for Latinons $\deltaL(\cdot,\cdot)$, which is a counterpart to the cut distance for graphons. The definition is given in Definition~\ref{def:cutdistlati}. Before giving preliminary definitions that are needed to this end, let us explain on a high level, the features we want this cut distance to have. Suppose that we have Latinons $L_1=(W,f)$ and $L_2=(U,g)$ over the ground space $\Omega$. Firstly, suppose that for each interval $(a,b)\subset [0,1]$ we have that the law of $W$ on $f^{-1}((a,b))$ is the same as the law of $U$ on $g^{-1}((a,b))$.\footnote{Note that these laws are on $\mathcal{B}$-valued functions.} Informally, this means that for almost every $x\in[0,1]$, $W$ on the fibre $f^{-1}(x)$ can be obtained by permuting $U$ on the fibre $g^{-1}(x)$ in a measure preserving fashion. In such a scenario, we clearly want to have $\deltaL(L_1,L_2)=0$. We want $\deltaL(L_1,L_2)$ to be small if and only if we are in a situation which is only a slightly perturbed version of the above scenario, in that
\begin{enumerate}[label={[cd\arabic*]}]
	\item\label{cd1} the laws of $W$ and $U$ on the individual fibres need not be equal, but it suffices to be approximately equal for most fibres,
	\item\label{cd2} instead of pairing up a fibre $f^{-1}(x)$ with $g^{-1}(x)$, we may pair it up with $g^{-1}(x')$ if $x\approx x'$ (for most choices of $x$).
\end{enumerate}
Definition~\ref{def:cutnormlati} quantifies~\ref{cd1}. In Definition~\ref{def:orderfix} we then introduce `order fixing bigraphons'. These order fixing bigraphons are then used in the terms $\cutn{O^f-O^{g\circ\vphi}}$ and $\cutn{O^f-O^{g\circ\psi }}$ in the main formula~\eqref{eq:defdeltaL} to quantify~\ref{cd2}.

\begin{defi}[Cut norm for distribution-valued signed bigraphons]\label{def:cutnormlati}
	Let $W:\Omega^2 \to \mathcal{B}^0$. We define
	\[\cutnL{W} \coloneqq \sup_{\substack{R,C\subseteq\Omega,\\V\subseteq[0,1]\textrm{ interval}}} \left| \int_{x \in R} \int_{y \in C} W(x,y)(V) \diffsymb x \diffsymb y \right|. \]
\end{defi}

\begin{defi}[Order fixing bigraphon]\label{def:orderfix}
	Define $O:[0,1]^2\rightarrow[0,1]$ to be the bigraphon such that $$O(x,y)=\begin{cases}1,\; \textrm{ if }x<y;\\0,\textrm{ otherwise.}\end{cases}$$
\end{defi}
\begin{defi}[Cut distance for Latinons]\label{def:cutdistlati}
	Let $L_1=(W,f)$ and $L_2=(U,g)$ be Latinons over the ground space $\Omega$. We define 
	\begin{equation}\label{eq:defdeltaL}
		\deltaL(L_1,L_2)\coloneqq\inf_{\vphi,\psi\in S_{\Omega}}\left(\cutn{O^f-O^{g\circ\vphi}}+\cutn{O^f-O^{g\circ\psi }}+\cutnL{W-U^{\vphi,\psi}}\right)\;.
	\end{equation}
	(Here, $W-U^{\vphi,\psi}$ is a difference between two distribution-valued bigraphons, and hence a distribution-valued signed bigraphon.)
\end{defi}
It is easy to verify that $\deltaL(\cdot,\cdot)$ is a pseudometric, i.e., it is symmetric and obeys the triangle inequality.

\begin{remark}\label{rem:cutdistsemi}
	Definition~\ref{def:cutdistlati} can also be extended to semilatinons, in which case it also obeys the triangle inequality.
\end{remark}

The two key lemmas below state that the topologies given by left-convergence and the cut distance are equivalent. We will prove Lemma~\ref{lem:count} in Section~\ref{sec:countinglemma} and Lemma~\ref{lem:invcounting} in Section~\ref{sec:invcount}.

\begin{lem}[Counting lemma for Latinons]\label{lem:count}
	Let $k,\ell\in\mathbb N$. Then there exists a constant $c_{k,\ell}$ such that for every $A\in\mathcal{R}(k,\ell)$ and Latinons $L_1,L_2$ over the ground space $\Omega$ we have
	\[|t(A,L_1)-t(A,L_2)|<c_{k,\ell} \cdot \deltaL(L_1,L_2)^{1/(2k\ell)}\;.\]
\end{lem}

\begin{lem}[Inverse counting lemma for Latinons]\label{lem:invcounting}
	Let $L_1$ and $L_2$ be Latinons over the ground space $\Omega$ such that $\deltaL(L_1,L_2)\geq d>0$. Then for every $k>2^{(120/d)^{10}}$ there exists a pattern $A\in\mathcal{R}(k,k)$ such that $|t(A,L_1)-t(A,L_2)|\geq \tfrac{1}{2(k^2!)}$.
\end{lem}

By combining Lemma~\ref{lem:count} and Lemma~\ref{lem:invcounting}, we immediately get the following.
\begin{thm}[Equivalence of convergence]\label{thm:equivalenceofconverge}
	Let $(L_n)_{n\in\mathbb N}$ be a sequence of Latinons and $L$ be a Latinon. The following are equivalent.
	\begin{enumerate}[label={(\roman*)}]
		\item $(L_n)_{n\in\mathbb N}\LeftConvergence L$.
		\item $(L_n)_{n\in\mathbb N}\overset{\deltaL}{\rightarrow}L$.
	\end{enumerate}
\end{thm}

\subsubsection{Latinons on different ground spaces}\label{sec:diffgroundspaces}
We said that throughout the paper, we will work with a fixed probability space $(\Omega,\mu)$ on which all the Latinons will live. There could be a situation where we need to work with an additional probability space, say $(\Omega',\mu')$. Even in this situation, we can define the cut distance between a Latinon $L$ over the ground space $\Omega$ and a Latinon $K$ over the ground space $\Omega'$, as follows,
$$\deltaL^*(L,K)\coloneqq\inf_{\psi:\Omega\rightarrow \Omega'\text{ measure preserving bijection}} \deltaL(L,K^{\psi})\;.$$ 
Here, we use the notation that if $K=(U,g)$ is a Latinon, then $K^\psi\coloneqq (U^\psi,g\circ \psi)$.

Firstly, note that transforming $K$ into $K^{\psi}$ does not change the densities.

Secondly, it is easy to check that $\deltaL^*(\cdot,\cdot)$ is an extension of $\deltaL(\cdot,\cdot)$, that is, if $\Omega'=\Omega$, then $\deltaL^*(L,K)=\deltaL(L,K)$.

Thirdly, it is straightforward to check that $\deltaL^*(\cdot,\cdot)$ is symmetric and it satisfies the triangle inequality, i.e., if in addition $H$ is a Latinon over a ground space $\Omega''$ then we have $\deltaL^*(L,H)\le \deltaL^*(L,K)+\deltaL^*(K,H)$. 

In view of Theorem~\ref{thm:equivalenceofconverge} and Proposition~\ref{prop:latinonrep} it is also reasonable to define a cut distance between a Latinon $L=(W,f)$ over the ground space $\Omega$ and a Latin square $L_n$ via its Latinon-representation $L' \coloneqq (W_{L_n},f)$: we define $\deltaL(L,L_n) \coloneqq \deltaL(L,L')$.

\subsection{A sampling lemma}
The key result for proving the inverse counting lemma is a sampling lemma.\footnote{This approach parallels that for graphons; see Section~10.6 of~\cite{Lovasz2012}.}
The sampling lemma, stated below as Lemma~\ref{lem:sampling}, tells us that 
a pattern sampled from a Latinon (as in Definition~\ref{def:latinonsampling})
is close in the cut distance.

Recall that we have a way of creating a Latinon from a Latin square (Definition~\ref{def:latinonrep}). We also require an analogy for general matrices; since the object we start with is not a Latin square, we cannot hope to produce a Latinon, so we instead produce a semilatinon. The key result is that the semilatinon that will be associated with a sampled matrix from a Latinon is close to the Latinon in the cut distance with high probability.

\begin{defi}\label{def:assocsemilatinons}Suppose that $k\in\mathbb{N}$.
	\begin{enumerate}[label={(\roman*)}]
		\item For a distribution-valued matrix $A\in\mathcal{B}^{k\times k}$ together with a measure preserving map $f:\Omega\rightarrow [0,1]$ we define an associated semilatinon $L^A=(W^A,f)$ in the following way. 
		Given the partition $\partition{k}$,
		set $W^A(x,y)\coloneqq A_{i,j}$ for $x\in f^{-1}(
		\interval{k}{i}
		)$ and $y\in f^{-1}(
		\interval{k}{j}
		)$. 
		\item For a matrix $B\in[0,1]^{k\times k}$ we define the associated semilatinon $L^B=L^{A}$, where $A\in \mathcal{B}^{k\times k}$ is the distribution-valued matrix defined by $A_{i,j}=\Dirac_{B_{i,j}}$ for $i,j\in[k]$.
		\item For a pattern $C \in \mathcal{R}(k,k)$, define the matrix $B\in[0,1]^{k \times k}$ by $B_{i,j}\coloneqq C_{i,j}/k^2$. Then the associated semilatinon $L^C$ of the pattern $C$ is given by $L^{B}$ defined in~(ii).
	\end{enumerate}
\end{defi}

We can now state the sampling lemma. Note that by Remark~\ref{rem:cutdistsemi}, the cut distance can also be applied to semilatinons.
\begin{lem}[Sampling lemma]\label{lem:sampling}
	Suppose that $L=(W,f)$ is a Latinon and $k\in\mathbb{N}$ is arbitrary. If $A\in[0,1]^{k\times k}$ is sampled from $L$
	according to Definition~\ref{def:latinonsampling},
	then with probability at least $1-30\exp(-k^{1/8}/10)$ we have
	\[\deltaL\left(L,(W^A,f)\right)\leq \frac{30}{\log(k)^{1/8}}\;.\]
\end{lem}
We prove the sampling lemma in Section~\ref{sec:sampling}.

\subsection{Return to the standard cyclic Latinon example}\label{ssec:cyclicproof}
Recall the standard cyclic Latinon $(L,f)$ and the Latin squares $L_n$ from Section~\ref{ssec:ExamplesOfLatinons}. We can now proceed with the proof that $L_n \LeftConvergence (L,f)$. First let $L'_n\coloneqq (W_{L_n},f)$ denote the Latinon-representation of $L_n$. By Proposition~\ref{prop:latinonrep} and Lemma~\ref{lem:count}, it suffices to show that $L'_n \overset{\deltaL}{\rightarrow} (L,f)$.

By taking $\phi$ and $\psi$ to be the identity, we obtain 
\begin{align}\label{standard1}
	\deltaL(L'_n,(L,f)) \leq \sup_{\substack{R,C,V\subseteq[0,1]\\ V\textrm{interval}}}\int_{x\in R}\int_{y\in C}W_{L_n}(x,y)(V)-L(x,y)(V)\diffsymb y \diffsymb x.
\end{align}
Next by the uniform marginals property of Latinons, we have for any interval $V=[z_1,z_2]$,
\begin{align}\label{standard2}
	\int_{x \in [0,1]} \int_{y \in [0,1]} W_{L_n}(x,y)(V)-L(x,y)(V) \diffsymb y \diffsymb x \leq z_2-z_1.
\end{align}
We have one of the following three conditions:
\begin{itemize}
	\item{There exists an integer $a \in [0,n-1]$ such that $\frac{a}{n} \leq z_1 < z_2 \leq \frac{a+1}{n}$;}
	\item{There exists an integer $a \in [1,n-1]$ such that $\frac{a-1}{n} \leq z_1 \leq \frac{a}{n} \leq z_1 \leq \frac{a+1}{n}$;}
	\item{We can write $V$ as $[z_1,\frac{a}{n}] \cup [\frac{a}{n},\frac{b}{n}] \cup [\frac{b}{n},z_2]$ for integers $a,b \in [0,n-1]$, $a<b$.}
\end{itemize}
Thus via (\ref{standard1}) and (\ref{standard2}) we obtain
\begin{align}\label{standardsup}
	\deltaL(L'_n,(L,f)) \leq \frac{2}{n} + \sup_{\substack{R,C \subseteq [0,1] \\ a,b \in [0,n-1]}} \int_{x\in R}\int_{y\in C}W_{L_n}(x,y)\left(\left[\frac{a}{n},\frac{b}{n}\right]\right)-L(x,y)\left(\left[\frac{a}{n},\frac{b}{n}\right]\right)\diffsymb y \diffsymb x,
\end{align}
where $a,b \in [0,n-1]$ and $a<b$.
For each $r,c \in [0,n-1]$ define 
\begin{align*}
	A_{r,c} & \coloneqq \left\{(x,y)\mid x\in \interval{n}{r}, y \in \interval{n}{c}, x+y \leq \frac{(r+c \! \mod n)+1}{n}\right\}; \\
	B_{r,c} & \coloneqq \left\{(x,y)\mid x\in \interval{n}{r}, y \in \interval{n}{c}, x+y > \frac{(r+c \! \mod n)+1}{n}\right\}.
\end{align*}

We have
\begin{align*}
	\int_{(x,y) \in A_{r,c}} W_{L_n}(x,y) \left( \interval{n}{v} \right) \diffsymb y \diffsymb x= & \begin{cases} \frac{1}{2n^2} & \text{if } v=r+c \! \mod n; \\ 0 & \text{otherwise;} \end{cases} \\
	\int_{(x,y) \in B_{r,c}} W_{L_n}(x,y) \left( \interval{n}{v} \right) \diffsymb y \diffsymb x= &\begin{cases} \frac{1}{2n^2} & \text{if } v=r+c \! \mod n; \\ 0 & \text{otherwise;} \end{cases} \\
	\int_{(x,y) \in A_{r,c}} L(x,y) \left( \interval{n}{v} \right) \diffsymb y \diffsymb x = &\begin{cases} \frac{1}{2n^2} & \text{if } v=r+c \! \mod n; \\ 0 & \text{otherwise;} \end{cases} \\
	\int_{(x,y) \in B_{r,c}} L(x,y) \left( \interval{n}{v} \right) \diffsymb y \diffsymb x = &\begin{cases} \frac{1}{2n^2} & \text{if } v=r+c-1 \! \mod n; \\ 0 & \text{otherwise.} \end{cases} \\
\end{align*}
Combining the above we have
\begin{align*}
	& \int_{(x,y) \in D_{r,c}} W_{L_n}(x,y) \left( \interval{n}{v} \right) -  L(x,y) \left( \interval{n}{v} \right) \diffsymb y \diffsymb x \\
	= & \begin{cases} \frac{1}{2n^2} & \text{if } D=B \text{ and } v=r+c \! \mod n; \\ -\frac{1}{2n^2} & \text{if } D=B \text{ and } v=r+c-1 \! \mod n; \\ 0 & \text{otherwise.} \end{cases}
\end{align*}

Now by a telescoping sum we obtain
\begin{align*}
	& \int_{(x,y) \in D_{r,c}} W_{L_n}(x,y) \left( \left[ \frac{a}{n},\frac{b}{n} \right] \right) -  L(x,y) \left( \left[ \frac{a}{n},\frac{b}{n} \right] \right) \diffsymb y \diffsymb x \\
	= & \begin{cases} \frac{1}{2n^2} & \text{if } D=B \text{ and } b=r+c-1 \! \mod n; \\ -\frac{1}{2n^2} & \text{if } D=B \text{ and } a=r+c-1 \! \mod n; \\ 0 & \text{otherwise.} \end{cases}
\end{align*}
From this it is clear that the supremum in (\ref{standardsup}) is bounded by $\frac{1}{2n}$ (since at most there are $n$ areas of weight $\frac{1}{2n^2}$ in any choice of $R$ and $C$). Thus we have $\lim_{n \to \infty} \deltaL(L'_n,(L,f)) \to 0$ as required.

\section{Vectors of bigraphons}\label{sec:prelim}

In our proofs we will often reduce the assertion to a claim about a vector of bigraphons. We therefore recall several results of this theory which we will use later on. Recall that a bigraphon is a measurable function $W: \Omega^2 \to [0,1]$ which is not necessarily symmetric. We denote the space of bigraphons by $\mathcal{W}_0$.
We can define a metric $\delta_\square^{\mathbb N}$ on $\mathcal{W}_0^{\mathbb N}$ by setting 
\begin{equation}\label{eq:definftyd}
	\inftyd{
		(U_n)_{n\in\mathbb N},(W_n)_{n\in\mathbb N}
	}
	\coloneqq
	\inf_{\vphi \in S_{\Omega}}\sum_{n=1}^\infty\frac{1}{2^n}\cdot \cutn{U_n-W_n^{\vphi}}\;.
\end{equation}
Note that this is not the same as $\sum_{n=1}^\infty\tfrac{1}{2^n}\delta_{\square}(U_n,W_n)$ and therefore compactness of $(\mathcal{W}_0^{\mathbb N},\delta_\square^{\mathbb N})$ does not follow immediately from Tychonoff's theorem applied to $(\mathcal{W}_0,\delta_\square)$.

\begin{thm}[Generalised compactness for bigraphons]\label{thm:genLoSz}
	$(\mathcal{W}_0^{\mathbb N},\delta_\square^{\mathbb N})$ is compact.
\end{thm}

Theorem~\ref{thm:genLoSz} is not really new. In fact, most proofs of the Lov\'asz--Szegedy compactness theorem generalise. Nevertheless, we provide a sketch of the generalisation of the original proof of Lov\'asz and Szegedy. To this end, we need the following `multicolour' version of the (weak) regularity lemma. Given a measurable partition $\mathcal{P}$ of $\Omega$ and a bigraphon $W$, we write $W^{\Join\mathcal{P}}$ for the \emph{stepping} of $W$ according to $\mathcal{P}$. For each $A,B\in\mathcal{P}$, $W^{\Join\mathcal{P}}$ on $A\times B$ is defined to be a constant $\frac{1}{\mu(A)\mu(B)}\cdot\int_{A\times B}W(x,y) \diffsymb(x,y)$.
Recall that $\mathcal{P}$ is an \emph{equipartition} if each cell has measure $\frac{1}{|\mathcal{P}|}$.

Next, let us state a weak regularity lemma for tuples of bigraphons.
\begin{lem}\label{lem:weak_reg_tup}
	Suppose that $r,m\in\mathbb{N}$ are arbitrary. Let $(W_1,\dots,W_m)$ be an $m$-tuple of bigraphons on $\Omega^2$. Let $\mathcal{P}^*$ be an arbitrary partition of $\Omega$. Then,
	\begin{enumerate}[label={(\roman*)}]
		\item\label{rl:part1}
		there exists a partition $\mathcal{P}$ of $\Omega$ into at most $r\cdot |\mathcal{P}^*|$ classes that refines $\mathcal{P}^*$ and such that
		$\cutn{W_i-(W_i)^{\Join\mathcal{P}}}
		<
		\sqrt{\tfrac{2m}{\log(r)}}$
		for each $i\in[m]$, and
		\item\label{rl:part2} 
		there exists an equipartition $\mathcal{K}$ of $\Omega$ into at most $r^2\cdot |\mathcal{P}^*|$ classes that refines $\mathcal{P}^*$ and such that
		$\cutn{W_i-(W_i)^{\Join\mathcal{K}}}
		<
		\sqrt{\tfrac{2m}{\log(r)}}+\frac{2}{r}$
		for each $i\in[m]$.
	\end{enumerate}
\end{lem}

We consider Lemma~\ref{lem:weak_reg_tup} standard enough to omit giving a full proof, but still provide a summary. Let us start with Part~\ref{rl:part1}. Like in the usual proofs of the weak regularity lemma (for a single graph/graphon), the central quantity is that of an index of a partition. That is, we start with a trivial partition $\mathcal{P}\coloneqq \mathcal{P}^*$, and record $m$ indices of that partition  with respect to all bigraphons $(W_1,\dots,W_m)$. Now, if $\cutn{W_i-(W_i)^{\Join\mathcal{P}}}\ge \sqrt{\tfrac{2m}{\log(r)}}$ for some $i$, then there exists a refinement of $\mathcal{P}$ in which each cluster is split into at most~2, and such that $i$-th index goes up by at least $\tfrac{m}{\log(r)}$. Also, a refinement of a partition never decreases an index, which in particular applies to the remaining $m-1$ indices. Since the sum of the indices is between $0$ and $m$, we conclude that we end up with a regular partition in at most $m/{\tfrac{m}{\log(r)}}=\log(r)$  steps. So the number, say $N$, of clusters at the end of this process will be at most $|\mathcal{P}^*| \cdot 2^{\log(r)}=r\cdot |\mathcal{P}^*|$. 

For the proof of Part~\ref{rl:part2} we use the proof for Part~\ref{rl:part1}, but need to turn the partition into an equipartition. We cut each cell $X$ of $\mathcal{P}$ arbitrarily into $\lfloor \frac{\mu(X)}{rN}\rfloor$ cells, of individual measure $\frac{1}{rN}$ each; of course, there is a possible leftover of measure less than $\frac{1}{rN}$. We aggregate all the leftover sets; this aggregate set $A$ has measure less than $N\cdot \frac{1}{rN}=\frac{1}{r}$. We cut $A$ arbitrarily into $\frac{\mu(A)}{rN}$ cells, of individual measure $\frac{1}{rN}$ each (note that $\mu(A)$ is an integer multiple of $\frac{1}{rN}$, so this is really possible). Clearly, we now have an equipartition $\mathcal{K}$ with exactly $rN\le r^2\cdot |\mathcal{P}^*|$ parts. Furthermore, it is a property of the cut norm distance that it does not increase due to refinements. So, the only increase in the cut norm distance could have been caused by the leftover $A$, that is, $\cutn{W_i-(W_i)^{\Join\mathcal{K}}}
\le 
\cutn{W_i-(W_i)^{\Join\mathcal{P}}}+2\mu(A)< \sqrt{\tfrac{2m}{\log(r)}}+\frac{2}{r}$, as needed.

\begin{proof}[Sketch of proof of Theorem~\ref{thm:genLoSz}]
	Let us consider an arbitrary sequence $\mathbf{W}^1,\mathbf{W}^2,\ldots$ of countable tuples of bigraphons, $\mathbf{W}^i=(W_1^i,W_2^i,\ldots)$. For a given $i$ and $m\in\mathbb{N}$, let us write $\mathbf{W}_{\restriction m}^i$ for the $m$-tuple $\mathbf{W}_{\restriction m}^i\coloneqq (W_1^i,W_2^i,\ldots,W_{m-1}^i,W_m^i)$. Let us take a sequence $(\eps_\ell)_{\ell\in\mathbb{N}}$ of positive numbers that tend to~$0$. Set $\mathbf{j}_0\coloneqq (1,2,\ldots)$ and $\mathcal{P}_{0,q}\coloneqq \{\Omega\}$ for each $q\in\mathbb{N}$. Set $r_0\coloneqq 1$ and $r_{i+1}\coloneqq r_i\cdot 2^{2/\eps_i^3}$.
	
	Now, for $\ell=1,2,\ldots$ we do the following. Set $m\coloneqq \lfloor 1/\eps_\ell\rfloor$. For each $q=1,2,\ldots$, we apply Lemma~\ref{lem:weak_reg_tup} with a prepartition $\mathcal{P}_{\ell-1,q}$, $m$-tuple of bigraphons that is $q$-th with respect to $\mathbf{j}_{\ell-1}$, that is, $\mathbf{W}_{\restriction m}^{\mathbf{j}_{\ell-1}(q)}$, with prepartition $\mathcal{P}_{\ell-1,q}$. Lemma~\ref{lem:weak_reg_tup} yields a partition $\mathcal{Q}_{\ell,q}$, $|\mathcal{Q}_{\ell,q}|\le r_i\cdot |\mathcal{P}_{\ell-1,q}|\le r_{\ell}$, which refines $\mathcal{P}_{\ell-1,q}$ and with the property that $\delta_{\square}\left(W^{\mathbf{j}_{\ell-1}(q)}_i,(W^{\mathbf{j}_{\ell-1}(q)}_i)^{\Join \mathcal{Q}_{\ell,q}}\right)<\eps_\ell$ for each $i\in [\lceil 1/\eps_\ell\rceil]$. We now prepare ourselves for the step $\ell+1$. We pass to an infinite subsequence $\mathbf{j}_\ell\subset \mathbf{j}_{\ell-1}$ so that in this subsequence all the partitions $\mathcal{Q}_{\ell,q}$ have the same size, say $R$, and so that we have the following.
	\begin{itemize}
		\item Order the cells of $\mathcal{Q}_{\ell,q}$ arbitrarily as $Q_{\ell,q;1},\ldots,Q_{\ell,q;R}$. Then the vectors $(\mu(Q_{\ell,q;1}),\ldots,\mu(Q_{\ell,q;R}))$ converge, as $q\rightarrow \infty$, 
		\item For each $i\in [\lceil 1/\eps_\ell\rceil]$, the bigraphons $(W^{\mathbf{j}_{\ell}(q)}_i)^{\Join \mathcal{Q}_{\ell,q}}$, when viewed as $R\times R$ matrices ordered as above, converge, as $q\rightarrow \infty$.
	\end{itemize}
	Note that such a subsequence exists, since the above convergence notions refer to the (compact) spaces $[0,1]^R$ and $[0,1]^{R^2}$, respectively. We can proceed with the step $\ell+1$.
	
	Now, the diagonalisation is done as in the Lov\'asz--Szegedy proof.
\end{proof}

We will also need to consider homomorphism densities of oriented graphs in $m$-tuples of bigraphons. This extension of the concept of a homomorphism to tuples plays a prominent role for example in~\cite{Zh15}. Recall that in an oriented graph, each edge has an initial and a terminal vertex; parallel or counterparalled edges and self-loops are not allowed.

\begin{defi}[Densities in vectors of bigraphons]
	Let $\mathbf{W}=(W_1,\dots,W_m)$ be an $m$-tuple of bigraphons on $\Omega^2$, let $F$ be an oriented graph and let $\alpha:E(F)\rightarrow[m]$ be any map. We then set
	\[t_\alpha(F,\mathbf{W})\coloneqq\int_{\mathbf{x}\in\Omega^{V(F)}}\prod_{(x_i,x_j)\in E(F)}W_{\alpha(x_i,x_j)}(x_i,x_j) \diffsymb \mathbf{x}\;.\] 
\end{defi}

We will make use of a counting lemma which is the following slightly generalised version of a lemma from \cite[Lemma 5.3]{Zh15}.

\begin{lem}[Counting lemma for bigraphons]\label{lem:counting}
	Let $\mathbf{U}=(U_1,\dots,U_m)$ and $\mathbf{W}=(W_1,\dots,W_m)$ be two $m$-tuples of bigraphons. Suppose that $\cutn{W_i-U_i}\leq\eps$ for each $i\in[m]$. Then for any oriented graph $F$ and any map $\alpha:E(F)\rightarrow [m]$ we have
	\[|t_\alpha(F,\mathbf{U})-t_\alpha(F,\mathbf{W})|\leq \eps \cdot e(F)\;,\]
	where $e(F)$ is the number of edges in $F$.
\end{lem}

Furthermore, we will make use of the following `First sampling lemma' for bigraphons in a form analogous to that stated in~\cite[Lemma 10.6]{Lovasz2012} for graphons. Note that considering bigraphons does not affect the argument of the proof, so the analogous statement holds. Given a bigraphon $U$ on $[0,1]^2$, $k\in\mathbb N$, and $S,T\in[0,1]^k$, we define $U[S,T]$ to be the bigraphon on $[0,1]^2$ defined in the following way. For every $i,j \in [k]$ and 
$x \in \interval{k}{i}$, $y \in \interval{k}{j}$
we set $U[S,T](x,y)\coloneqq U(S_i,T_j)$.
\begin{lem}[First sampling lemma for bigraphons]\label{lem:firstsampling}
	Let $k\geq 1$, and let $U$ be a bigraphon on $[0,1]^2$. If $S,T$ are chosen independently uniformly at random from $[0,1]^k$, then with probability at least $1-4\exp(-\sqrt{k}/10)$,
	\[ \left| \cutn{U[S,T]}-\cutn{U} \right| \leq \tfrac{8}{k^{1/4}}\;.\]
\end{lem}

We will need an ordered version which fits better to Definition~\ref{def:assocsemilatinons}. It follows immediately from Lemma~\ref{lem:firstsampling} by noting that $\cutn{U^{\vphi,\psi}}=\cutn{U}$ for every $\vphi,\psi\in S_{[0,1]}$.

\begin{lem}[First sampling lemma for bigraphons (ordered version)]\label{lem:firstsamplingordered}
	Let $k\geq 1$, and let $U$ be a bigraphon on $[0,1]^2$. If $S,T$ are chosen independently uniformly at random from $[0,1]^k_<$, then with probability at least $1-4\exp(-\sqrt{k}/10)$,
	\[ \left| \cutn{U[S,T]}-\cutn{U} \right| \leq \tfrac{8}{k^{1/4}}\;.\]
\end{lem}

\section{Compressions of Latinons}\label{sec:compressions}

Several of our proofs are based on the idea that we can approximate a Latinon by a vector of bigraphons, where the length of the vector depends on the precision of the approximation. We will also consider vectors of infinite length which contain all of the finite approximations. For the following recall that by $\dyapartition{d}$ we denote the partition of $[0,1)$ into $2^d$ parts given by $\dyainterval{d}{s}\coloneqq[\frac{s-1}{2^d},\frac{s}{2^d})$ for $i\in[2^d]$ and that $\alldyapartition=\bigcup_{d\in\mathbb N}\dyapartition{d}$.

\begin{defi}[Compressions of Latinons]\label{def:compression}
	Let $d\in\mathbb N$. For a distribution-valued bigraphon $W$ on $\Omega^2$ we define the \emph{compression of depth} $d$ of $W$ to be the $2^d$-tuple of bigraphons on $\Omega^2$
	\[\mathbf{W}^d=(W_{d,1},\dots,W_{d,2^d})\;,\]
	where for every $s\in[2^d]$ we have
	\[W_{d,s}(x,y)=W(x,y)(\dyainterval{d}{s})\;
	\text{ for all } x,y\in \Omega\;.\]
	We define the \emph{infinite compression} of $W$ to be the vector $\mathbf{W}^{\mathbb N}\in\mathcal{W}_0^{\mathbb N}$ with
	\[\mathbf{W}^{\mathbb N}=(\mathbf{W}^1,\mathbf{W}^2,\dots)\;.\]
	For a Latinon $L=(W,f)$ we define the \emph{compression of depth} $d$ of $L$ to be the $2^d+1$-tuple of bigraphons on $\Omega^2$
	\[\mathbf{L}^d=(O^f,\mathbf{W}^d)\;.\]
	We define the \emph{infinite compression} of $L$ to be the vector $\mathbf{L}^{\mathbb N}\in\mathcal{W}_0^{\mathbb N}$ with
	\[\mathbf{L}^{\mathbb N}=(O^f,\mathbf{W}^{\mathbb N})\;.\]
\end{defi}

The next observation states that most times when working with a Latinon, it does not make a difference whether we work with open, closed or half open intervals in the definition of a compression. 

\begin{obs}\label{obs:intervals}
	Suppose that $L=(W,f)$ is a Latinon. Then there exists a nullset $S\subseteq\Omega^2$ such that for every $(x,y)\in\Omega^2\setminus S$ we have that $W(x,y)((a,b))=W(x,y)([a,b])$ for all $a,b\in\mathbb Q\cap[0,1]$ with $a<b$.
\end{obs}

\begin{proof}
	It is enough to prove that for each $t\in [0,1]$, the set of $(x,y)\in\Omega^2$ for which $W(x,y)(\{t\})>0$ is null. Indeed, once this is settled, it is enough to take $S$ to be the union of these exceptional nullsets, over all $t\in\mathbb Q\cap[0,1]$. But the former easily follows from the uniform marginals property (Definition~\ref{def:latinon}),
	$$\int_{x\in\Omega}\int_{y\in\Omega} W(x,y)(\{t\})\diffsymb x\diffsymb y=\lambda(\{t\})=0\;.$$
\end{proof}

Thus, with the help of Observation~\ref{obs:intervals}, we have the following nestedness property.
\begin{obs}\label{obs:nested}
	Suppose that $L=(W,f)$ is a Latinon. Then there exists a nullset $S\subseteq\Omega^2$ such that for each $d\in\mathbb{N}$ and each $(x,y)\in \Omega^2\setminus S$ we have for its compressions of depth $d$ and $d+1$ and each $i\in [2^d]$ that $W_{d,i}(x,y)=W_{d+1,2i-1}(x,y)+W_{d+1,2i}(x,y)$.
	
	In particular, for $(x,y)\in \Omega^2\setminus S$ we have $\sum_{i=1}^{2^d}W_{d,i}(x,y)=1$.
\end{obs}

Next, we show how to proceed in the other direction, from compressions to Latinons.
\begin{defi}
	Suppose that $d\in\mathbb{N}$. Given a $2^d$-tuple of bigraphons $\mathbf{W}=(W_1,\dots,W_{2^d})$ on $\Omega^2$ with 
	\begin{equation}\label{eq:roky}
		\sum_{i=1}^{2^d}W_i(x,y)=1
	\end{equation} for each $(x,y)\in\Omega^2$, we can define a distribution-valued bigraphon, denoted by $\AntiCompr{\mathbf{W}}$, for any $(x,y)\in\Omega^2$ and measurable set $X\subset [0,1]$, by
	\[
	\AntiCompr{\mathbf{W}}(x,y)(X)\coloneqq \sum_{i=1}^{2^d}W_i(x,y)\cdot 2^d \cdot \lambda(X\cap\dyainterval{d}{i})\;.
	\]
	Note that~\eqref{eq:roky} ensures that the values of $\AntiCompr{\mathbf{W}}$ are indeed probability distributions.
\end{defi}
The following observations are immediate.
\begin{fact}\label{fact:anticompressions}
	Suppose that $(W,f)$ is a Latinon and $d\in\mathbb{N}$. Let $(O^f,W_1,\dots,W_{2^d})$ be the compression of depth $d$ of $W$. Then $(\AntiCompr{W_1,\dots,W_{2^d}},f)$ is a Latinon. Further, its compression of depth $d$ is $(O^f,W_1,\dots,W_{2^d})$.
\end{fact}
Fact~\ref{fact:anticompressions} allows to extend the symbol $\AntiCompr{\cdot}$ even to compressions of Latinons.

\begin{lem}[Latinon-distance of compressions]\label{lem:cutdistcomp} Let $d\in\mathbb N$.
	\begin{enumerate}[label={(\roman*)}]
		\item\label{en:caj1} If $W$ is a distribution-valued bigraphon with uniform marginals and $\mathbf{W}^d$ its compression of depth $d$, then
		\[\cutnL{W-\AntiCompr{\mathbf{W}^d}}\leq\frac{1}{2^{d-1}}\;.\]
		\item\label{en:caj2} If $L$ is a Latinon and $\mathbf{L}^d$ its compression of depth $d$, then
		\[\deltaL\left(L,\AntiCompr{\mathbf{L}^d}\right)\leq\frac{1}{2^{d-1}}\;.\]
	\end{enumerate}
\end{lem}

\begin{proof}
	For~\ref{en:caj1} let $S,T\subseteq\Omega$ and $V\subseteq [0,1]$ be an interval with its smallest value in $\dyainterval{d}{s}$ and its largest in $\dyainterval{d}{t}$. We then have
	\begin{align*}
		&\int_{S}\int_TW(x,y)(V)\diffsymb y \diffsymb x\\
		=&\int_{S}\int_T \sum_{i=s+1}^{t-1} W_{d,i}(x,y)\diffsymb y \diffsymb x + \int_{S}\int_TW(x,y)\left(V\cap(\dyainterval{d}{s}\cup \dyainterval{d}{t})\right)\diffsymb y \diffsymb x\\
		=&\int_S\int_T \mathbf{W}^d(x,y)(V)\diffsymb y \diffsymb x\pm \tfrac{2}{2^d}\;,
	\end{align*}
	where we use that $W$ has uniform marginals for the last equality.
	
	Next, observe that~\ref{en:caj2} follows immediately from~\ref{en:caj1}.
\end{proof}

The next fact helps to compare compressions of Latinons by comparing the individual bigraphons via cut distance.

\begin{prop}\label{prop:boundbycompr}
	Let $d\in\mathbb{N}$, let $W_1$, $W_2$ be distribution-valued bigraphons with uniform marginals, and let $\mathbf{W_1}^d=(W_1^{d,1},\dots,W_1^{d,2^d})$ and $\mathbf{W_2}^d=(W_2^{d,1},\dots,W_2^{d,2^d})$ be their compressions of depth $d$. Then
	\[\cutnL{\mathbf{W_1}^d-\mathbf{W_2}^d}\leq\sum_{i=1}^{2^d} \cutn{W_1^{d,i}-W_2^{d,i}}\;.\]
\end{prop}

\begin{proof}
	Set $\mathbf{U}^d\coloneqq \mathbf{W_1}^d-\mathbf{W_2}^d$ and $U^{d,i}\coloneqq W_1^{d,i}-W_2^{d,i}$.
	For $S,T\subseteq\Omega$ and an interval $V\subseteq [0,1]$ we then have
	\begin{align*}
		\left|\int_{S}\int_T\mathbf{U}^d(x,y)(V)\diffsymb x \diffsymb y\right|=&\left|\int_{S}\int_T \sum_{i=1}^{2^d} \mathbf{U}^{d}(x,y)\left(V\cap \dyainterval{d}{i}\right)\diffsymb x \diffsymb y\right|\\
		=&\left|\sum_{i=1}^{2^d} 2^d \cdot \lambda\left(V\cap \dyainterval{d}{i}\right)\cdot \int_{S}\int_T  U^{d,i}(x,y)\diffsymb x \diffsymb y\right|\\
		\leq&\sum_{i=1}^{2^d}\cutn{W_1^{d,i}-W_2^{d,i}}\;.
	\end{align*}
	
\end{proof}

We conclude the section by showing that the infinity-distance of infinite compressions bounds the Latinon-distance of the finite compressions.

\begin{prop}\label{prop:cutdistinftydist}
	Let $d\in\mathbb N$. If $L_1$ and $L_2$ are Latinons, then
	\[\deltaL(\mathbf{L_1}^d,\mathbf{L_2}^d)<2^{2^{d+1}}\cdot\inftyd{\mathbf{L_1}^{\mathbb N},\mathbf{L_2}^{\mathbb N}}.\]
\end{prop}

\begin{proof}
	We write $\mathbf{L_1}^{\mathbb N}=(O^f,W_{1,1},\dots)$, $\mathbf{L_2}^{\mathbb N}=(O^g,U_{1,1},\dots)$ and $\delta=\inftyd{\mathbf{L_1}^{\mathbb N},\mathbf{L_2}^{\mathbb N}}$. Now there exists $\vphi\in S_{\Omega}$ such that
	\[\tfrac{1}{2}\cutn{O^{f\circ\vphi}-O^g}+\sum_{d=1}^\infty\sum_{k=1}^{2^d}\tfrac{1}{2^{2^{d}+k-1}}\cutn{W_{d,k}^\vphi-U_{d,k}}<2\delta\;,\]
	and therefore by Proposition~\ref{prop:boundbycompr},
	\[\deltaL(\mathbf{L_1}^d,\mathbf{L_2}^d)\leq2\cutn{O^{f\circ\vphi}-O^g}+\sum_{k=1}^{2^d}\cutn{W_{d,k}^\vphi-U_{d,k}}<2^{2^{d+1}}\delta\;.\]
\end{proof}

\section{Proof of the counting lemma}\label{sec:countinglemma}
In this section, we prove the counting lemma, Lemma~\ref{lem:count}. We make use of the compressions of Latinons defined in Section~\ref{sec:compressions}. For this recall Definition~\ref{def:strucequiv}, in particular that $\mathcal{R}^A(S)$ is the set of matrices with entries from a set $S$ which are structurally equivalent to a matrix $A$.
\begin{proof}[Proof of Lemma~\ref{lem:count}]
	Let $L_1=(W,f)$ and $L_2=(U,g)$ be Latinons, $A\in\mathcal{R}(k,\ell)$ and $d\in\mathbb N$. We write $\mathbf{L_1}^d=(O^f,W_{d,1},\dots,W_{d,2^d})$ and $\mathbf{L_2}^d=(O^g,U_{d,1},\dots,U_{d,2^d})$ for the corresponding compressions of depth $d$. The calculations below are valid for each $d\in \mathbb{N}$, which we will choose only in the last step of the proof.
	
	Note that by Definition~\ref{def:cutdistlati} there exist $\vphi,\psi\in S_{\Omega}$ such that
	\begin{align*}
		\deltaL(L_1,L_2)&\geq\cutn{O^f-O^{g\circ\vphi}}\;,
		\\
		\deltaL(L_1,L_2)&\geq\cutn{O^f-O^{g\circ\psi}}\;,
	\end{align*}
	and
	\begin{align*}
		\deltaL(L_1,L_2) \geq&\cutnL{W-U^{\vphi,\psi}}\\
		\geq&\sup_{R,C\subseteq\Omega}\left|\int_{x\in R}\int_{y\in C}W(x,y)
		(\dyainterval{d}{i})
		-U^{\vphi,\psi}(x,y)
		(\dyainterval{d}{i})
		\diffsymb x \diffsymb y\right|\\
		=&\sup_{R,C\subseteq\Omega}\left|\int_{x\in R}\int_{y\in C}W_{d,i}(x,y)-U_{d,i}(\vphi(x),\psi(y))\diffsymb x \diffsymb y\right|= \cutn{W_{d,i}-U_{d,i}^{\vphi,\psi}}\;,
	\end{align*}
	for every $i\in[2^d]$. We now set $\mathbf{W}=(O^{f},O^{f},W_{d,1},\dots,W_{d,2^d})$ and $\mathbf{U}=(O^{g\circ\vphi},O^{g\circ\psi},U_{d,1}^{\vphi,\psi},\dots,U_{d,2^d}^{\vphi,\psi})$. We also define an oriented graph $F$ with vertex set 
	$$\{a_i\mid i\in[k]\}\cup\{b_i\mid i\in[\ell]\}$$
	and edge set
	$$\{(a_i,b_j)\mid i\in[k],j\in[\ell]\}\cup\{(a_i,a_j)\mid 1\leq i<j\leq k\}\cup\{(b_i,b_j)\mid 1\leq i<j\leq \ell\}\;.$$
	Furthermore, given a matrix $P=(p_{i,j})_{i\in[k],j\in[\ell]}\in\mathcal{R}^A([2^d])$, we define a map $\alpha_P:E(F)\rightarrow[2+2^d]$ by $(a_i,a_j)\mapsto 1$, $(b_i,b_j)\mapsto 2$, and $(a_i,b_j)\mapsto p_{i,j}+2$ for $i\in[k],j\in[\ell]$. We then get by Lemma~\ref{lem:counting} and the above inequalities that
	\begin{equation}
		|t_{\alpha_P}(F,\mathbf{W})-t_{\alpha_P}(F,\mathbf{U})|\leq|F|\cdot\deltaL(L_1,L_2)= (k\ell+\tbinom{k}{2}+\tbinom{\ell}{2})\cdot\deltaL(L_1,L_2)\;.\label{eq:vectorbound}
	\end{equation}

	We will make use of the following notation. Let $S$ be a set and $(m_{i,j})=M\in S^{k\times \ell}$. Furthermore given a partition $\mathcal{Q}=(Q_1,\dots,Q_d)$ of a collection of subsets of $S$, we write $M\perp\mathcal{Q}$ if there are no $(i_1,j_1),(i_2,j_2)\in[k]\times[\ell]$ and $s\in[d]$ such that $(i_1,j_1)\neq (i_2,j_2)$ and $m_{i_1,j_1},m_{i_2,j_2}\in Q_s$. 
	
	\begin{claim}\label{clm:ILikeLabels1}
		We have
		$$\sum_{P\in\mathcal{R}^A([2^d])}t_{\alpha_P}(F,\mathbf{W})
		=
		\int_{\mathbf{x}\in\Omega^k_{<_f}}\int_{\mathbf{y}\in\Omega^\ell_{<_f}}\left(\bigotimes_{(i,j)\in[k]\times[\ell]}W(x_i,y_j)\right)(\{M\in\mathcal{R}^A([0,1])\mid M\perp \dyapartition{d}\})\diffsymb\mathbf{y}\diffsymb\mathbf{x}\;.
		$$  
	\end{claim}
	\begin{proof}[Proof of Claim~\ref{clm:ILikeLabels1}]
		Indeed,
		\begin{align*}
			\sum_{P\in\mathcal{R}^A([2^d])}t_{\alpha_P}(F,\mathbf{W})
			=&\sum_{P\in\mathcal{R}^A([2^d])}\int_{\mathbf{x}\in\Omega^k}\int_{\mathbf{y}\in\Omega^\ell}\prod_{i=1}^k\prod_{j=1}^\ell
			W_{d,p_{i,j}}(x_i,y_j)\mathbbm{1}_{\Omega^k_{<_f}}(\mathbf{x})\mathbbm{1}_{\Omega^\ell_{<_f}}(\mathbf{y})\diffsymb\mathbf{y}\diffsymb\mathbf{x}\\
			=&\sum_{P\in\mathcal{R}^A([2^d])}\int_{\mathbf{x}\in\Omega^k_{<_f}}\int_{\mathbf{y}\in\Omega^\ell_{<_f}}\prod_{i=1}^k\prod_{j=1}^\ell W_{d,p_{i,j}}(x_i,y_j)\diffsymb\mathbf{y}\diffsymb\mathbf{x}\\
			=&\sum_{P\in\mathcal{R}^A([2^d])}\int_{\mathbf{x}\in\Omega^k_{<_f}}\int_{\mathbf{y}\in\Omega^\ell_{<_f}}\prod_{i=1}^k\prod_{j=1}^\ell W(x_i,y_j)(\dyainterval{d}{p_{i,j}})
			\diffsymb\mathbf{y}\diffsymb\mathbf{x}\\
			=&\int_{\mathbf{x}\in\Omega^k_{<_f}}\int_{\mathbf{y}\in\Omega^\ell_{<_f}}\left(\bigotimes_{(i,j)\in[k]\times[\ell]}W(x_i,y_j)\right)(\bigcup_{P\in\mathcal{R}^A([2^d])}(\dyainterval{d}{p_{i,j}})_{(i,j)\in[k]\times[\ell]})\diffsymb\mathbf{y}\diffsymb\mathbf{x}\\
			=&\int_{\mathbf{x}\in\Omega^k_{<_f}}\int_{\mathbf{y}\in\Omega^\ell_{<_f}}\left(\bigotimes_{(i,j)\in[k]\times[\ell]}W(x_i,y_j)\right)(\{M\in\mathcal{R}^A([0,1])\mid M\perp
			\dyapartition{d}
			\})\diffsymb\mathbf{y}\diffsymb\mathbf{x}\;.
		\end{align*}
	\end{proof}
	\begin{claim}\label{clm:ILikeLabels2}
		We have
		$$ 
		\sum_{P\in\mathcal{R}^A([2^d])}t_{\alpha_P}(F,\mathbf{U})
		=
		\int_{\mathbf{x}\in\Omega^k_{<_g}}\int_{\mathbf{y}\in\Omega^\ell_{<_g}}\left(\bigotimes_{(i,j)\in[k]\times[\ell]}U(x_i,y_j)\right)(\{M\in\mathcal{R}^A([0,1])\mid M\perp\dyapartition{d}\})\diffsymb\mathbf{y}\diffsymb\mathbf{x}\;.
		$$
	\end{claim}
	\begin{proof}[Proof of Claim~\ref{clm:ILikeLabels2}]
		The calculations are similar as in Claim~\ref{clm:ILikeLabels1}, but in addition use that $\vphi$ and $\psi$ are measure preserving bijections.
		\begin{align*}
			\sum_{P\in\mathcal{R}^A([2^d])}t_{\alpha_P}(F,\mathbf{U})=&\sum_{P\in\mathcal{R}^A([2^d])}\int_{\mathbf{x}\in\Omega^k}\int_{\mathbf{y}\in\Omega^\ell}\prod_{i=1}^k\prod_{j=1}^\ell U_{d,p_{i,j}}(\vphi(x_i),\psi(y_j))\mathbbm{1}_{\Omega^k_{<_g}}(\vphi(\mathbf{x}))\mathbbm{1}_{\Omega^\ell_{<_g}}(\psi(\mathbf{y}))\diffsymb\mathbf{y}\diffsymb\mathbf{x}\\
			=&\sum_{P\in\mathcal{R}^A([2^d])}\int_{\mathbf{x}\in\Omega^k}\int_{\mathbf{y}\in\Omega^\ell}\prod_{i=1}^k\prod_{j=1}^\ell U_{d,p_{i,j}}(x_i,y_j)\mathbbm{1}_{\Omega^k_{<_g}}(\mathbf{x})\mathbbm{1}_{\Omega^\ell_{<_g}}(\mathbf{y})\diffsymb\mathbf{y}\diffsymb\mathbf{x}\\
			=&\int_{\mathbf{x}\in\Omega^k_{<_g}}\int_{\mathbf{y}\in\Omega^\ell_{<_g}}\left(\bigotimes_{(i,j)\in[k]\times[\ell]}U(x_i,y_j)\right)(\{M\in\mathcal{R}^A([0,1])\mid M\perp\dyapartition{d}\})\diffsymb\mathbf{y}\diffsymb\mathbf{x}\;.
		\end{align*}
	\end{proof}
	
	\begin{claim}\label{clm:line}
		We have
		\[\int_{x_1\in\Omega}\int_{y_1\in\Omega}\int_{y_2\in\Omega}
		W(x_{1},y_{1})
		(\dyainterval{d}{r})
		W(x_{1},y_{2})
		(\dyainterval{d}{r})
		\diffsymb y_2\diffsymb y_1\diffsymb x_1=\left(\tfrac{1}{2^d}\right)^2\]
		and
		\[\int_{x_1\in\Omega}\int_{x_2\in\Omega}\int_{y_1\in\Omega}
		W(x_{1},y_{1})
		(\dyainterval{d}{r})
		W(x_{2},y_{1})
		(\dyainterval{d}{r})
		\diffsymb y_1\diffsymb x_2\diffsymb x_1=\left(\tfrac{1}{2^d}\right)^2\]
	\end{claim}
	
	\begin{proof}[Proof of Claim~\ref{clm:line}]
		It is enough to prove the first equality, the second equality follows by symmetry of the following argument.
		\begin{align*}
			&\int_{x_1\in\Omega}\int_{y_1\in\Omega}\int_{y_2\in\Omega}
			W(x_{1},y_{1})
			(\dyainterval{d}{r})
			W(x_{1},y_{2})
			(\dyainterval{d}{r})
			\diffsymb y_2\diffsymb y_1\diffsymb x_1\\
			=&\int_{x_1\in\Omega}\int_{y_1\in\Omega}W(x_{1},y_{1})
			(\dyainterval{d}{r})
			\int_{y_2\in\Omega}
			W(x_{1},y_{2})
			(\dyainterval{d}{r})
			\diffsymb y_2\diffsymb y_1\diffsymb x_1\\
			=&\int_{x_1\in\Omega}\int_{y_1\in\Omega}W(x_{1},y_{1})
			(\dyainterval{d}{r})
			\mu_{W,x_1}^1(\Omega\times 
			\dyainterval{d}{r}
			)\diffsymb y_1\diffsymb x_1\\
			=&\int_{x_1\in\Omega}\mu_{W,x_1}^1(\Omega\times 
			\dyainterval{d}{r}
			)\int_{y_1\in\Omega}W(x_{1},y_{1})
			(\dyainterval{d}{r})
			\diffsymb y_1\diffsymb x_1\\
			=&\int_{x_1\in\Omega}\mu_{W,x_1}^1(\Omega\times 
			\dyainterval{d}{r}
			)\mu_{W,x_1}^1(\Omega\times 
			\dyainterval{d}{r}
			)\diffsymb x_1\\
			=&\int_{x_1\in\Omega}\tfrac{1}{2^d}\tfrac{1}{2^d}\diffsymb x_1\\
			=&\left(\tfrac{1}{2^d}\right)^2\;,
		\end{align*}
		where we use the fact that the 2-dimensional measure $\mu_{W,x_1}^1$ has uniform marginals for almost every $x_1\in\Omega$, as $(W,f)$ is a Latinon. 
	\end{proof}

	\begin{claim}\label{clm:diagonal}
		We have
		\[
		\int_{x_1\in\Omega}\int_{x_2\in\Omega}\int_{y_1\in\Omega}\int_{y_2\in\Omega}
		W(x_{1},y_{1})
		(\dyainterval{d}{r})
		W(x_{2},y_{2})
		(\dyainterval{d}{r})
		\diffsymb y_2\diffsymb y_1\diffsymb x_2\diffsymb x_1=\left(\tfrac{1}{2^d}\right)^2
		\;.\]
	\end{claim}
	
	\begin{proof}[Proof of Claim~\ref{clm:diagonal}]
		\begin{align*}
			&\int_{x_1\in\Omega}\int_{x_2\in\Omega}\int_{y_1\in\Omega}\int_{y_2\in\Omega}
			W(x_{1},y_{1})
			(\dyainterval{d}{r})
			W(x_{2},y_{2})
			(\dyainterval{d}{r})
			\diffsymb y_2\diffsymb y_1\diffsymb x_2\diffsymb x_1\\
			=&\int_{x_1\in\Omega}\int_{y_1\in\Omega}W(x_{1},y_{1})
			(\dyainterval{d}{r})
			\int_{x_2\in\Omega}\int_{y_2\in\Omega}
			W(x_{2},y_{2})
			(\dyainterval{d}{r})
			\diffsymb y_2\diffsymb x_2\diffsymb y_1\diffsymb x_1\\
			=&\int_{x_1\in\Omega}\int_{y_1\in\Omega}W(x_{1},y_{1})
			(\dyainterval{d}{r})
			\int_{x_2\in\Omega}\mu_{W,x_2}^1
			(\dyainterval{d}{r})
			\diffsymb x_2\diffsymb y_1\diffsymb x_1\\
			=&\tfrac{1}{2^d}\int_{x_1\in\Omega}\int_{y_1\in\Omega}W(x_{1},y_{1})
			(\dyainterval{d}{r})
			\diffsymb y_1\diffsymb x_1\\
			=&\tfrac{1}{2^d}\int_{x_1\in\Omega}\mu_{W,x_1}^1
			(\dyainterval{d}{r})
			\diffsymb x_1\\
			=&\left(\tfrac{1}{2^d}\right)^2\;,
		\end{align*}
		where we use the fact that the 2-dimensional measure $\mu_{W,x}^1$ has uniform marginals for almost every $x\in\Omega$, as $(W,f)$ is a Latinon.
	\end{proof}
	
	We therefore get
	\begin{align}
		\begin{split}\label{eq:split1}
			&\left|t(A,L_1)-k!\ell!\sum_{P\in\mathcal{R}^A([2^d])}t_{\alpha_P}(F,\mathbf{W})\right|\\
			\overset{C\ref{clm:ILikeLabels1}}=&k!\ell!\int_{\mathbf{x}\in\Omega^k_{<_f}}\int_{\mathbf{y}\in\Omega^\ell_{<_f}}\left(\bigotimes_{(i,j)\in[k]\times[\ell]}W(x_i,y_j)\right)(\{M\in\mathcal{R}^A([0,1])\mid M\not\perp\dyapartition{d}\})\diffsymb\mathbf{y}\diffsymb\mathbf{x}\\
			\leq& k!\ell!\sum_{r\in[2^d]}\sum_{\substack{(i_1,j_1)\\ \neq(i_2,j_2)\\ \in[k]\times[\ell]}}\int_{\substack{\mathbf{x}\in\Omega^k\\\mathbf{y}\in\Omega^\ell}}W(x_{i_1},y_{j_1})
			(\dyainterval{d}{r})
			W(x_{i_2},y_{j_2})
			(\dyainterval{d}{r})
			\prod_{\substack{(i,j)\in[k]\times[\ell]\\ \neq(i_1,j_1),(i_2,j_2)}}W(x_i,y_j)([0,1])\diffsymb\mathbf{y}\diffsymb\mathbf{x}\\
			=&k!\ell!\sum_{r\in[2^d]}\sum_{\substack{i_1,i_2\in[k],j_1,j_2\in[\ell]\\i_1\neq i_2, j_1\neq j_2}}\int_{\substack{(x_{i_1},x_{i_2})\in\Omega^2\\(y_{j_1},y_{j_2})\in\Omega^2}}W(x_{i_1},y_{j_1})
			(\dyainterval{d}{r})
			W(x_{i_2},y_{j_2})
			(\dyainterval{d}{r})
			\diffsymb(y_{j_1},y_{j_2})\diffsymb(x_{i_1},x_{i_2})\\
			+&k!\ell!\sum_{r\in[2^d]}\sum_{\substack{i_1,i_2\in[k],j_1,j_2\in[\ell]\\i_1= i_2, j_1\neq j_2}}\int_{\substack{x_{i_1}\in\Omega\\(y_{j_1},y_{j_2})\in\Omega^2}}W(x_{i_1},y_{j_1})
			(\dyainterval{d}{r})
			W(x_{i_1},y_{j_2})
			(\dyainterval{d}{r})
			\diffsymb(y_{j_1},y_{j_2})\diffsymb(x_{i_1})\\
			+&k!\ell!\sum_{r\in[2^d]}\sum_{\substack{i_1,i_2\in[k],j_1,j_2\in[\ell]\\i_1\neq i_2, j_1= j_2}}\int_{\substack{(x_{i_1},x_{i_2})\in\Omega^2\\y_{j_1}\in\Omega}}W(x_{i_1},y_{j_1})
			(\dyainterval{d}{r})
			W(x_{i_2},y_{j_1})
			(\dyainterval{d}{r})
			\diffsymb y_{j_1}\diffsymb(x_{i_1},x_{i_2})\\
			\overset{C\ref{clm:line},\ref{clm:diagonal}}=&k!\ell!\cdot 2^d\cdot k\ell(k\ell-1)\cdot(\tfrac{1}{2^d})^2=k!\ell!\cdot k\ell(k\ell-1)\cdot\tfrac{1}{2^d}\;.
		\end{split}
	\end{align}
	
	A similar calculation to~\eqref{eq:split1} which uses Claim~\ref{clm:ILikeLabels2} instead of Claim~\ref{clm:ILikeLabels1} yields
	\begin{equation}\label{eq:split2}
		\left| t(A,L_2)-k!\ell!\sum_{P\in\mathcal{R}^A([2^d])}t_{\alpha_P}(F,\mathbf{U})\right|
		\leq k!\ell!\cdot k\ell(k\ell-1)\cdot\tfrac{1}{2^d}\;.
	\end{equation}
	
	We are now in a position to bound $|t(A,L_1)-t(A,L_2)|$. Hence, take $d$ such that $$\tfrac{1}{\deltaL(L_1,L_2)^{1/(2k\ell)}}\leq 2^d\leq\tfrac{2}{\deltaL(L_1,L_2)^{1/(2k\ell)}}$$ and set $c_{k,\ell}\coloneqq2k!\ell!k\ell(k\ell-1)+2^{k\ell}k!\ell!(k\ell+\tbinom{k}{2}+\tbinom{\ell}{2})$.
	Starting with the triangle inequality, we get
	\begin{align*}
		&|t(A,L_1)-t(A,L_2)|\\
		\leq &\left|t(A,L_1)-k!\ell!\sum_{P\in\mathcal{R}^A([2^d])}t_{\alpha_P}(F,\mathbf{W})\right|\\
		&+k!\ell!\sum_{P\in\mathcal{R}^A([2^d])}|t_{\alpha_P}(F,\mathbf{W})-t_{\alpha_P}(F,\mathbf{U})|\\
		&+\left|k!\ell!\sum_{P\in\mathcal{R}^A([2^d])}t_{\alpha_P}(F,\mathbf{U})-t(A,L_2)\right|\\
		\overset{C\ref{clm:ILikeLabels1},~\eqref{eq:vectorbound},C\ref{clm:ILikeLabels2}}\leq &2k!\ell!k\ell(k\ell-1)/2^d+|\mathcal{R}^A([2^d])|k!\ell!(k\ell+\tbinom{k}{2}+\tbinom{\ell}{2})\deltaL(L_1,L_2)\\
		\leq &2k!\ell!k\ell(k\ell-1)/2^d+2^{dk\ell}k!\ell!(k\ell+\tbinom{k}{2}+\tbinom{\ell}{2})\deltaL(L_1,L_2)\\
		\leq &2k!\ell!k\ell(k\ell-1)\deltaL(L_1,L_2)^{1/(2k\ell)}+2^{k\ell}k!\ell!(k\ell+\tbinom{k}{2}+\tbinom{\ell}{2})\sqrt{\deltaL(L_1,L_2)}\\
		\leq &c_{k,\ell}\deltaL(L_1,L_2)^{1/(2k\ell)}\;.
	\end{align*}
\end{proof}

\section{Proof of compactness}\label{sec:compactness}

In this section, we prove Theorem~\ref{thm:compact}. In fact,  Theorem~\ref{thm:compact} follows from Lemma~\ref{lem:count}, remarks in Section~\ref{sec:diffgroundspaces} and the proposition below, which we prove instead.
\begin{prop}\label{prop:Latinoncompact}
	The space $(\LatinonSpace,\deltaL)$ is compact.
\end{prop}
Proposition~\ref{prop:Latinoncompact} follows immediately from Lemma~\ref{lem:compconthot1} and Lemma~\ref{lem:compconthot2} below. To state these lemmata, we introduce a function $\iota$ which maps a Latinon to its infinite compression (see Definition~\ref{def:compression}), i.e.
\[\iota:(\LatinonSpace,\deltaL)\rightarrow(\mathcal{W}_0^{\mathbb N},\delta_\square^{\mathbb N}),\;L\mapsto \mathbf{L}^{\mathbb N}\;.\]
As we show in the lemma below, the function $\iota$ is injective.
\begin{lem}\label{lem:iotainjective}
	The function $\iota$ is injective.     
\end{lem}
\begin{proof}
	Suppose that $L=(W,f),L'=(W',f')\in\LatinonSpace$ are such that $\inftyd{\iota(L),\iota(L')}=0$. We need to show that $\deltaL(L,L')<\alpha$ for an arbitrary given $\alpha>0$. Let us write $h=\lceil 4/\alpha\rceil$. Let us write $\iota(L)=(O^f,W_{1,1},W_{1,2},\dots)$ and $\iota(L')=(O^{f'},W'_{1,1},W'_{1,2},\dots)$. Let $\vphi\in S_{\Omega}$ be as in~\eqref{eq:definftyd} such that 
	\begin{equation}\label{eq:halohalo}
		\frac12\cutn{O^f-O^{f'\circ\vphi}}+\sum_{\ell=1}^\infty \sum_{q=1}^{2^\ell} \frac{1}{2^{2^{\ell}+q-1}}\cdot \cutn{W_{\ell,q}-(W'_{\ell,q})^{\vphi}}<\frac{\alpha}{2^{3+2^{h+1}}}\;.
	\end{equation}
	Taking $\vphi=\psi$ in~\eqref{eq:defdeltaL}, we see (by \eqref{eq:halohalo}) that the sum of the first and the second term is at most $\frac{4\alpha}{2^{3+2^{h+1}}}<\alpha/4$. To conclude that $\deltaL(L,L')<\alpha$, it thus remains to bound $\cutnL{W-(W')^{\vphi}}$. So, let $R,C\subseteq\Omega$ be arbitrary, and let $V=[v_1,v_2]\subseteq[0,1]$ be an interval. Let $V'=[w_1,w_2]$ be an interval whose endpoints $w_1,w_2$ are dyadic rationals with dyadic denominators at most $2^h$ and with $|w_i-v_i|<\alpha/4$ (clearly, such a dyadic approximation exists). By the uniform marginals property, we have 
	$$\int_{x \in R}\int_{y \in C} (W-(W')^{\vphi})(x,y)(V)\diffsymb y\diffsymb x=\int_{x \in R}\int_{y \in C} (W-(W')^{\vphi})(x,y)(V')\diffsymb y\diffsymb x \;\pm\; (|w_1-v_1|+|w_2-v_2|)\;,$$
	and so it suffices to prove that $|\int_{x \in R}\int_{y \in C} (W-(W')^{\vphi})(x,y)(V')\diffsymb y\diffsymb x|<\alpha/4$, which we do now. Let us write $Q\subset [2^h]$ for the indices of dyadic intervals at depth $h$ encoding $V'$, i.e., $Q\coloneqq\{i\in[2^h]\mid
	\dyainterval{h}{i}
	\subset V'\}$. Then we have
	\begin{align*}
		\left|\int_{x \in R}\int_{y \in C} (W-(W')^{\vphi})(x,y)(V') \diffsymb y\diffsymb x \right|
		&=
		\left|\int_{x \in R}\int_{y \in C} \sum_{i\in Q}(W_{h,i}-(W'_{h,i})^{\vphi})(x,y) \diffsymb y\diffsymb x \right|\\
		&\le 
		\sum_{i\in Q}\left|\int_{x \in R}\int_{y \in C} (W_{h,i}-(W'_{h,i})^{\vphi})(x,y) \diffsymb y\diffsymb x\right|\\
		&\le 
		\sum_{i\in Q}\cutn{W_{h,i}-(W'_{h,i})^{\vphi}}\overset{\eqref{eq:halohalo}}{<}\alpha/4\;,
	\end{align*}
	proving the statement.
\end{proof}
\begin{lem}\label{lem:compconthot1}
	The space $(\iota(\LatinonSpace),\delta_\square^{\mathbb N})$ is compact.
\end{lem} 
\begin{proof}
	Let $(L_n)_{n\in\mathbb N}$ be a sequence of Latinons. We have to show that $(\iota(L_n))_{n\in\mathbb N}$ contains a convergent subsequence with respect to $\delta_\square^{\mathbb N}$. For each $n\in\mathbb N$ we write $\iota(L_n)=(O^{n},W_{1,1}^n,W_{1,2}^n,\dots)$. By Theorem~\ref{thm:genLoSz} there exist bigraphons $\widetilde{O},\widetilde{W}_{d,k}$, $d\in\mathbb N,k\in[2^d]$, such that, after passing to a subsequence, we have that
	\begin{equation}
		(O^{n},W^{n}_{1,1},W^{n}_{1,2},W^{n}_{2,1},W^{n}_{2,2},W^{n}_{2,3},\dots)\overset{\delta_\square^{\mathbb N}}{\rightarrow}\mathbf{\widetilde{W}}=(\widetilde{O},\widetilde{W}_{1,1},\widetilde{W}_{1,2},\widetilde{W}_{2,1},\widetilde{W}_{2,2},\widetilde{W}_{2,3},\dots)\;.\label{eq:conv}
	\end{equation}
	Again by passing to a subsequence we can additionally assume that for every $n\in\mathbb N$,
	\begin{equation}
		\inftyd{\iota(L_n),\mathbf{\widetilde{W}}}<\tfrac{1}{2^{n+2}}\;,\label{eq:subsdist}
	\end{equation}
	and in particular we can fix $\vphi_n\in S_{\Omega}$ for every $n\in\mathbb N$ such that
	\begin{equation}
		\tfrac{1}{2}\cutn{O^{n,\vphi_n}-\widetilde{O}}+\sum_{d=1}^\infty\sum_{k=1}^{2^d}\tfrac{1}{2^{2^{d}+k-1}}\cutn{W_{d,k}^{n,\vphi_n}-\widetilde{W}_{d,k}}<\tfrac{1}{2^{n+1}}\;.\label{eq:subsnorm}
	\end{equation}
	
	We have to show that there exists a Latinon $L=(\nu,f)$ such that $\iota(L)=\mathbf{\widetilde{W}}$ $=(\widetilde{O},\widetilde{W}_{1,1},\widetilde{W}_{1,2},\widetilde{W}_{2,1},\dots)$.
	
	First, we derive a measure preserving map $f:\Omega\rightarrow [0,1]$ from $\widetilde{O}$ by setting $f(x)\coloneqq\deg_{\widetilde{O}}(x)$. Note that then $O^f=\widetilde{O}$ up to a null set. Since the degree distribution function of each $O^{n,\vphi_n}$ is measure preserving, Lemma~\ref{lem:degditcont} implies that this property is inherited to the limit, i.e. $f$ is indeed measure preserving.
	
	For the definition of the Latinon-measure $\nu$ note that the dyadic intervals $\alldyapartition$ together with $\emptyset$ form a semiring which we call $R$, and that $R$ generates the Borel sigma-algebra on $[0,1)$.\footnote{See Section~\ref{sssec:semirings} for constructions of measures from semirings.} Consequently, 
	\[
	R'\coloneqq\{S\times T\times J\mid S,T\subseteq\Omega\text{ Borel measurable},J\in \alldyapartition\}\cup\{\emptyset\}
	\]
	is a semiring which generates the Borel sigma-algebra on $\Omega^2\times[0,1)$. We define 
	$$\nu^*(S\times T\times \dyainterval{d}{s} )\coloneqq \int_S\int_T\widetilde{W}_{d,s}(x,y)\diffsymb x\diffsymb y$$
	for all $S,T\subset \Omega$, $d\in\mathbb N$ and $s\in[2^d]$. We want to show that $\nu^*$ defines a premeasure on $\Omega^2\times[0,1)$. This is shown by the following claim.
	
	\begin{claim}\label{clm:subadditivity}
		For every partition $S\times T\times 
		\dyainterval{d'}{k'}
		=\bigcup_{i=1}^\infty S_i\times T_i\times 
		\dyainterval{d_i}{k_i}
		$, where $S,T,S_i,T_i\subseteq\Omega$ and $\dyainterval{d'}{k'},\dyainterval{d_i}{k_i}
		\in\alldyapartition$, we have
		\[\nu^*(S\times T\times 
		\dyainterval{d'}{k'}
		)=\sum_{i=1}^\infty\nu^*(S_i\times T_i\times 
		\dyainterval{d_i}{k_i}
		)\;.\]
	\end{claim}
	
	\begin{proof}[Proof of Claim~\ref{clm:subadditivity}]
		Let $S\times T\times 
		\dyainterval{d'}{k'}
		=\bigcup_{i=1}^\infty S_i\times T_i\times 
		\dyainterval{d_i}{k_i}
		$ be as in the assumption. We need to show that for an arbitrary $\eps>0$ (which we consider fixed from now on) there exists $N_0\in\mathbb{N}$ such that for every $N>N_0$ we have
		\begin{equation}\label{eq:MilosZeman}
			\left|\nu^*(S\times T\times 
			\dyainterval{d'}{k'}
			)-\sum_{i=1}^N \nu^*(S_i\times T_i\times 
			\dyainterval{d_i}{k_i}
			)\right|<\eps\;.
		\end{equation}
		We define a map $g:
		\dyainterval{d'}{k'}
		\rightarrow\mathbb N$ by
		\begin{equation}\label{eq:g(p)def}
			g(p)\coloneqq\min \left \{r\in\mathbb N \mid \mu^{\otimes 2}\left(S\times T\setminus\bigcup_{i\in[r]\textrm{ s.t. }p\in 
				\dyainterval{d_i}{k_i}
			}S_i\times T_i\right)<\eps/4 \right\}. 
		\end{equation}
		Note that $g$ is measurable and finite everywhere and that $\sum_{i=1}^\infty \lambda(g^{-1}(i))=\lambda(\bigcup_{i\in\mathbb N}g^{-1}(i))= \lambda
		(\dyainterval{d'}{k'})
		\le 1$ is convergent. Therefore there exists $N_0\in\mathbb N$ such that
		\begin{equation}\label{eq:ginversesmall}
			\sum_{i=N_0+1}^\infty \lambda(g^{-1}(i))<\eps/4\;.
		\end{equation}
		We set $E\coloneqq\bigcup_{i=N_0+1}^\infty g^{-1}(i)$.
		We now prove the following claim.
		\begin{claim}\label{clm:marg}
			Suppose that $(U,h)$ is an arbitrary Latinon. Then
			\[\sum_{i=N_0+1}^\infty\int_{S_i}\int_{T_i}U(x,y)
			(\dyainterval{d_i}{k_i})
			\diffsymb x\diffsymb y <\eps/2\;.\]
		\end{claim}
		
		\begin{proof}
			$U$ induces a measure on $\Omega^2\times[0,1]$ with uniform marginals. Let $\{\xi_p\}_{p\in[0,1]}$ be the family on $\Omega^2$ given by applying the disintegration theorem to this measure with respect to the third coordinate. It follows from the uniform marginals property that almost every measure $\{\xi_p\}_{p\in[0,1]}$ is a probability measure. Then
			\begin{align*}
				&\sum_{i=N_0+1}^\infty\int_{S_i}\int_{T_i} U(x,y)
				(\dyainterval{d_i}{k_i})
				\diffsymb x\diffsymb y\\
				&=\sum_{i=N_0+1}^\infty\int_{[0,1]} \mathbbm{1}_{
					\dyainterval{d_i}{k_i}
				}(p)\cdot\xi_p(S_i\times T_i)\diffsymb p
				=\sum_{i=N_0+1}^\infty\int_{
					\dyainterval{d'}{k'}
				} \mathbbm{1}_{
					\dyainterval{d_i}{k_i}
				}(p)\cdot\xi_p(S_i\times T_i)\diffsymb p\\
				&=\int_{
					\dyainterval{d'}{k'}
				}\xi_p\left(\bigcup_{i>N_0\textrm{ s.t. }p\in 
					\dyainterval{d_i}{k_i}
				}S_i\times T_i\right)\diffsymb p\\
				&=\int_{E}\xi_p\left(\bigcup_{i>N_0\textrm{ s.t. }p\in 
					\dyainterval{d_i}{k_i}
				}S_i\times T_i\right)\diffsymb p+\int_{
					\dyainterval{d'}{k'}
					\setminus E}\xi_p\left(\bigcup_{i>N_0\textrm{ s.t. }p\in 
					\dyainterval{d_i}{k_i}
				}S_i\times T_i\right)\diffsymb p\\
				&\le\int_{E}\xi_p\left(\Omega^2\right)\diffsymb p+\int_{
					\dyainterval{d'}{k'}
					\setminus E}\xi_p\left(\bigcup_{i>N_0\textrm{ s.t. }p\in 
					\dyainterval{d_i}{k_i}
				}S_i\times T_i\right)\diffsymb p
				\overset{\eqref{eq:g(p)def},\eqref{eq:ginversesmall}}{<}\eps/4+\eps/4=\eps/2\;.
			\end{align*}
		\end{proof}

		Let us now return to proving~\eqref{eq:MilosZeman}. Suppose that $N>N_0$ is arbitrary. Choose $M\in\mathbb N$ such that $\tfrac{N+1}{M}<\eps/2$ and $d_{max}=\max(\{d_i\mid i\in [N]\}\cup\{d'\})$. By~\eqref{eq:subsnorm} there exists $n_0\in\mathbb N$ such that for all $n\geq n_0$, $d \leq d_{max}$ and $k \in [2^d]$,
		\begin{equation}\label{eq:DonaltDrump}
			\left|\int_S\int_T W_{d,k}^{n,\vphi_{n}}(x,y)\diffsymb x\diffsymb y-\int_S\int_T\widetilde{W}_{d,k}(x,y)\diffsymb x\diffsymb y\right|<\tfrac{1}{M}\;.
		\end{equation}
		Recalling the definition of $\nu^*$, we then get for $n\geq n_0$ that
		\begin{align*}
			&\;\left|\nu^*(S\times T\times 
			\dyainterval{d'}{k'}
			)-\sum_{i=1}^N\nu^*(S_i\times T_i\times 
			\dyainterval{d_i}{k_i}
			)\right|\\
			\leq&\;\left|\nu^*(S\times T\times 
			\dyainterval{d'}{k'}
			-\int_S\int_T W^{n,\vphi_{n}}_{d',k'}(x,y)\diffsymb x\diffsymb y\right|\\
			+&\;\left|\int_S\int_T W^{n,\vphi_{n}}_{d',k'}(x,y)\diffsymb x\diffsymb y-\sum_{i=1}^N\int_{S_i}\int_{T_i} W^{n,\vphi_{n}}_{d_i,k_i}(x,y)\diffsymb x\diffsymb y\right|\\
			+&\;\sum_{i=1}^N\left| \int_{S_i}\int_{T_i}W^{n,\vphi_{n}}_{d_i,k_i}(x,y)\diffsymb x\diffsymb y-\nu^*(S_i\times T_i\times 
			\dyainterval{d_i}{k_i}
			)\right|\\
			\overset{2\times\eqref{eq:DonaltDrump}}{\leq}&\;\tfrac{1}{M}+\left|\sum_{i=N+1}^\infty\int_{S_i}\int_{T_i} W^{n,\vphi_{n}}_{d_i,k_i}(x,y)\diffsymb x\diffsymb y\right|+\tfrac{N}{M}\overset{C\ref{clm:marg}}{<}\eps\;,
		\end{align*}
		indeed proving~\eqref{eq:MilosZeman}.
	\end{proof}
	
	By Claim~\ref{clm:subadditivity} we now have that $\nu^*$ defines a premeasure on $R'$. Therefore by Carath\'eodory's extension theorem (Theorem~\ref{thm:cara}) there exists a unique measure $\nu$ on $\Omega^2\times[0,1]$ such that $\nu(S\times T\times 
	\dyainterval{d}{s}
	)=\int_S\int_T\widetilde{W}_{d,s}(x,y)(
	\dyainterval{d}{s}
	)\diffsymb x\diffsymb y$ for all $d\in\mathbb N$ and $s\in[2^d]$. 
	Lastly, we need to check that $\nu$ has uniform marginals;
	this will also verify that the total mass of $\nu$ is~1.
	Note that it is sufficient to restrict ourselves to sets of the semiring $R'$. Let $S,T\subseteq\Omega$ and $J\in\alldyapartition$. By making use of the uniform marginals of $W^n$ we get
	\begin{align*}
		\nu(S\times T\times[0,1])&=\lim_{n\rightarrow\infty}\int_S\int_T W^{n,\vphi_n}(x,y)([0,1])\diffsymb x\diffsymb y=\lim_{n\rightarrow\infty}\mu(S)\mu(T)=\mu(S)\mu(T)\;,\\
		\nu(\Omega\times T\times J)&=\lim_{n\rightarrow\infty}\int_\Omega\int_T W^{n,\vphi_n}(x,y)(J)\diffsymb x\diffsymb y=\lim_{n\rightarrow\infty}\mu(T)\lambda(J)=\mu(T)\lambda(J)\;,\\
		\nu(S\times \Omega\times J)&=\lim_{n\rightarrow\infty}\int_S\int_\Omega W^{n,\vphi_n}(x,y)(J)\diffsymb x\diffsymb y=\lim_{n\rightarrow\infty}\mu(S)\lambda(J)=\mu(S)\lambda(J)\;.
	\end{align*}
	Hence $\nu$ is a Latinon according to Definition~\ref{def:Latinonasmeasure}.
\end{proof}

To state Lemma~\ref{lem:compconthot2}, we need to recall that $\iota^{-1}$ is well-defined by Lemma~\ref{lem:iotainjective}.
\begin{lem}\label{lem:compconthot2}
	The map $\iota^{-1}:(\iota(\LatinonSpace),\delta_\square^{\mathbb N})\rightarrow (\LatinonSpace,\deltaL)$ is continuous.
\end{lem}
\begin{proof}
	Let $\eps>0$. We fix $d\in\mathbb N$ large enough such that $\tfrac{1}{2^{d-2}}<\eps/2$ and then choose $\delta>0$ small enough such that $2^{2^{d+1}}\cdot\delta<\eps/2$. Let $L_1$ and $L_2$ be a pair of Latinons such that $\delta_\square^{\mathbb N}(\iota(L_1),\iota(L_2))<\delta$. By Lemma~\ref{lem:cutdistcomp} and Proposition~\ref{prop:cutdistinftydist} we get that
	\[\deltaL(L_1,L_2)\leq \deltaL(L_1,\mathbf{L_1}^d)+\deltaL(\mathbf{L_1}^d,\mathbf{L_2}^d)+\deltaL(\mathbf{L_2}^d,L_2)\leq\tfrac{1}{2^{d-2}}+2^{2^{d+1}}\delta<\eps\;.\] 
\end{proof}

\begin{proof}[Proof of Proposition~\ref{prop:Latinoncompact}]
	By Lemma~\ref{lem:compconthot1}, the space $(\iota(\LatinonSpace),\delta_\square^{\mathbb N})$ is compact. By Lemma~\ref{lem:compconthot2}, $\iota^{-1}$ is continuous. Since a continuous image of a compact space is compact, we conclude that $(\LatinonSpace,\deltaL)$ is compact.
\end{proof}

\subsection{An additional application of the method: compactness of orderons}\label{ssec:orderonscompactness}
Here, we show how a simplified version of the above approach can be used to treat orderons. An \emph{orderon}, as defined~\cite{BeFiLeYo18}, is a symmetric function $W:([0,1]\times[0,1])^2\rightarrow[0,1]$. Here, the first component $[0,1]$ encodes the global position and the second component $[0,1]$ is used to encode local information. Definition~2.7 of~\cite{BeFiLeYo18} gives a definition of `CS-distance' between two orderons $U$ and $W$. The philosophy behind the `CS-distance' is to modify the usual cut distance --- in which one would consider all measure preserving bijections $\psi$ from $[0,1]\times [0,1]$ to itself and take the infimum of $\cutn{U-W^\psi}$ --- by penalising those $\psi$ that shuffle the global coordinates substantially. So, instead of repeating the definition from~\cite{BeFiLeYo18}, let us give an alternative but equivalent one:\footnote{`Equivalent' meaning that it generates the same topology as Definition~2.7 of~\cite{BeFiLeYo18}. This is routine to check.} define 
$$\delta_\triangle(U,W)\coloneqq\inf_{\psi} \cutn{O^{\pi_1}-O^{\pi_1\circ \psi}}+\cutn{U-W^\psi}\;,$$
where $\psi$ is as above, $\pi_1:[0,1]\times [0,1]\rightarrow [0,1]$ is the projection on the first coordinate and $O$ is the order fixing bigraphon.

Let us show that this alternative view together with the strategy we used to prove Proposition~\ref{prop:Latinoncompact} allows us to easily establish compactness of orderons (with respect to $\delta_\triangle$); the original proof occupies Section~4 of~\cite{BeFiLeYo18}. Indeed, suppose that we are given a sequence of orderons $W_1,W_2,W_3,\ldots$. We represent these orderons as pairs $(O^{\pi_1},W_1), (O^{\pi_1},W_2), (O^{\pi_1},W_3),\ldots$. Applying Theorem~\ref{thm:genLoSz} (on pairs rather than on infinite vectors), we get that there exists a sequence of indices $i_1<i_2<\ldots$, a bigraphon $\widetilde{O}$, a graphon $\widetilde{W}$, and a sequence of measure preserving bijections $\psi_{i_1},\psi_{i_2},\ldots$ such that
\begin{align*}
	\cutn{O^{\pi_1\circ \psi_{i_n}}-\widetilde{O}}&\rightarrow 0\;\mbox{, and}\\
	\cutn{W_{i_n}^{\psi_{i_n}}-\widetilde{W}}&\rightarrow 0\;.\\
\end{align*}
We might stop here and say that the pair $(\widetilde{O},\widetilde{W})$ is the sought after accumulation point. But let us do one more step and transform $(\widetilde{O},\widetilde{W})$ back into the global-local representation which is used in~\cite{BeFiLeYo18}. We proceed as described in Section~\ref{ssec:fromStandardToLocGlob}. That is, we appeal to Proposition~\ref{prop:transformIntoLocalGlobal} with a function $f:[0,1]\times [0,1]\rightarrow [0,1]$, $f((a,x))\coloneqq\deg_{\widetilde{O}}((a,x))$ and get a function $h:[0,1]\times[0,1]\rightarrow[0,1]\times[0,1]$. Then $W^*((a,x),(b,y))\coloneqq\widetilde{W}(h(a,x),h(b,y))$ is the accumulation point in the global-local representation.

\section{Proof of the sampling lemma for Latinons}\label{sec:sampling}

In this section we prove the sampling lemma, Lemma~\ref{lem:sampling}. The proof splits into two parts. In the first part, we show that we can approximately mimic sampling from a Latinon $L$ on $\Omega^2$ by sampling from a suitably chosen Latinon $(U,id)$ on $[0,1]^2$. In the second part, we show that a pattern sampled from a Latinon $(U,id)$ on $[0,1]^2$ is close to $(U,id)$ in the cut distance for Latinons. We start with the first part contained in the following lemma.

\begin{lem}\label{lem:couplsampl}
	Suppose that $L=(W,f)$ is a Latinon on $\Omega^2$, $k\in\mathbb N$, and $\eps>0$. Then there exists a Latinon $L'=(U,id)$ on $[0,1]^2$ such that 
	\begin{enumerate}[label={(\roman*)}]
		\item $\deltaL^*(L,L')<\eps$ and
		\item there exists a coupling $\mathcal{C}$ of sampling a $k\times k$ pattern $A$ from $L$ and sampling a $k\times k$ pattern $B$ from $L'$ such that $\Prob_{\mathcal{C}}[A\neq B]<\eps k^2$.
	\end{enumerate}
\end{lem}

\begin{proof}
	Choose $m\in\mathbb N$ such that $2/m<\eps$. 
	Given the partition $\partition{m}$, for each $i\in[m]$
	we choose an arbitrary measure preserving bijection $\vphi_i:f^{-1}(
	\interval{m}{i}
	)\rightarrow 
	\interval{m}{i}
	$ (where the domain is equipped with the restriction of the measure $\mu$ and the codomain is equipped with the Lebesgue measure). Define $\vphi:\Omega\rightarrow [0,1]$ by $\vphi(x)=\vphi_i(x)$, if $x\in 
	\interval{m}{i}
	$. We define $U$ to be the distribution-valued bigraphon on $[0,1]^2$ with $U(x,y)=W(\vphi^{-1}(x),\vphi^{-1}(y))$. Note that then for $L'=(U,id)$ we have
	\[\deltaL^*(L,L')\leq 2\cutn{O^{id}-O^{f\circ\vphi^{-1} }}+\cutnL{U-W^{\vphi^{-1}}}\leq 2/m+0<\eps\;,\]
	as $f\circ\vphi^{-1}$ only reorders elements inside each $
	\interval{m}{i}
	$. 
	We define the coupling for the sampling process in the following way. We sample $k$-sets $R,C\subseteq \Omega$ and order them according to $<_f$, written $(r_1,\dots,r_k)$ and $(c_1,\dots,c_k)$. We then consider the sets $\vphi(R)$ and $\vphi(C)$ and reorder the $k$-tuples according to $<$, written $(\vphi(r_{\sigma(1)}),\dots,\vphi(r_{\sigma(k)}))$ and $(\vphi(c_{\tau(1)}),\dots,\vphi(c_{\tau(k)}))$. We then sample values $a_{i,j}=b_{\sigma(i),\tau(j)}$ from $W(r_i,c_j)=U(\vphi(r_i),\vphi(c_j))$. Thus the event $A\neq B$ can only happen if $\tau\neq id$ or $\sigma\neq id$, which in turn can only happen if two elements of $R$ (or $C$) lie in the same set $f^{-1}(
	\interval{m}{i}
	)$ for some $i\in[k]$. Hence
	\[\Prob_{\mathcal{C}}[A\neq B]\leq 2\cdot \binom{k}{2}\cdot m\cdot(1/m)^2\leq k^2/m<\eps k^2\,.\]
\end{proof}

We can now move to the second part of the argument in which show that a pattern sampled from a Latinon $(U,id)$ on $[0,1]^2$ is close to $(U,id)$ in the cut distance. To this end we first need the following auxiliary lemma.

\begin{lem}\label{lem:SampPaternFromDistr}
	Suppose that $k\in\mathbb N$.
	Suppose that $L=(W,id)$ is a semilatinon on $[0,1]^2$ which is constant on every set 
	$\interval{k}{i} \times \interval{k}{j}$ for each $i,j \in [k].$
	For each $i,j\in[k]$ sample from 
	$W(\interval{k}{i},\interval{k}{j})$
	a value $a_{i,j}$ and call the resulting $k\times k$ matrix $A$. Then with probability at least $1-\exp(-5k)$ we have
	\[\cutnL{W-W^A}\leq 6/k^{1/4}\;.\]
\end{lem}

\begin{proof}
	We choose $d=\log(k^{1/4})$ and $\eps=2/\sqrt{k}$. For every $s\in[2^d]$ and $i,j$ let $X_{i,j}^s\in\{0,1\}$ be the random variable $X_{i,j}^s=\mathbbm{1}_{a_{i,j}\in 
		\dyainterval{d}{s}
	}$.
	Let $S,T\subseteq[k]$. For every $S,T\subseteq[k]$ and $s\in[2^d]$ we have by Lemma~\ref{lem:MBoundDiff} (with $Z=\sum_{i\in S}\sum_{j\in T}X_{i,j}^s$ and $t=\eps k^2$) that with probability at least $1-\exp(-2\eps^2k^2)$
	\begin{equation}\label{eq:IlseStoebe}
		\left|\sum_{i\in S}\sum_{j\in T}(X_{i,j}^s-\Expect[X_{i,j}^s])\right|\leq\eps k^2\;.
	\end{equation}
	Hence we have with probability at least $1-2^{2k}2^d\exp(-2\eps^2k^2)>1-\exp(-5k)$ that for every $s\in[2^d]$ and $S,T\subseteq[k]$
	\[\left|\int_{\bigcup_{i\in S}
		\interval{k}{i}
	}\int_{\bigcup_{j\in T}
		\interval{k}{j}
	}W^A_{d,s}(x,y)-W_{d,s}(x,y)\diffsymb x\diffsymb y\right|=\left|\tfrac{1}{k^2}\sum_{i\in S}\sum_{j\in T}(X_{i,j}^s-\Expect[X_{i,j}^s])\right|<\eps\;.\]
	Note that, since $W_{d,s}^A-W_{d,s}$ is a step-function, the maximum $\cutn{W^A_{d,s}-W_{d,s}}$ is attained on $\bigcup_{i\in S}
	\interval{k}{i}
	\times\bigcup_{j\in T}
	\interval{k}{j}
	$ for some $S,T\subseteq[k]$ and therefore
	\[\cutn{W^A_{d,s}-W_{d,s}}<\eps\;.\]
	By Lemma~\ref{lem:cutdistcomp} and Lemma~\ref{prop:boundbycompr} we therefore get that with probability at least $1-\exp(-5k)$
	\begin{align*}
		\cutnL{W^A-W}\leq\;&\cutnL{W^A-\bf{(W^A)}^d}+\cutnL{\bf{(W^A)}^d-\bf{W}^d}+\cutnL{\bf{W}^d-W}\\
		\leq\;&\tfrac{4}{2^d}+2^d\eps\leq 6/k^{1/4}\;.
	\end{align*}
\end{proof}

We also need the following definition of sampling a distribution-valued step-bigraphon, which can be seen as an intermediate step of sampling a pattern.

\begin{defi}\label{def:sampleDistr}
	Let $U$ be a distribution-valued bigraphon on $[0,1]^2$ and $R,C\in [0,1]^n_{<}$ be $k$-tuples. We write $R=(r_1,\dots,r_k)$, $C=(c_1,\dots,c_k)$. 
	We define $U[R,C]$ to be the distribution-valued bigraphon defined by $U[R,C](x,y)\coloneqq U(r_i,c_j)$, 
	if $x \in \interval{k}{i}$ and $y \in \interval{k}{j}$.
\end{defi}

Now we can prove the second part of the overall argument. Recall that by $W^{\Join\mathcal{P}}$ we denote the step bigraphon derived from a bigraphon $W$ with respect to the partition $\mathcal{P}$.

\begin{lem}\label{lem:sample01}
	Suppose that $(U,id)$ is a Latinon on $[0,1]^2$ and $k\in\mathbb N$. If $A$ is a $k\times k$ pattern sampled from $(U,id)$, then with probability at least $1-29\exp(-k^{1/8}/10)$ we have
	$$\deltaL((U,id),(W^A,id))<29/\sqrt{\log(k)}\;.$$
\end{lem}

\begin{proof}
	Set $m,r=k^{1/8}$ and $d=\log(\log(k)^{1/4})$ (we omit rounding symbols and assume that $m$, $k$ and $d$ are integers). Let $L^d=(O,U_{d,1},\dots,U_{d,2^d})$ be the compression of $L$ of depth $d$. 
	Consider the partition $\partition{r}$.
	By Lemma~\ref{lem:weak_reg_tup}\ref{rl:part2} we find an equipartition $\mathcal{P}$ of $[0,1]$ refining $\partition{r}$, i.e. 
	$\interval{r}{j}
	=P_{j,1}\cup\cdots\cup P_{j,m}$,
	and such that, for every $1\leq i\leq 2^d$,
	\begin{equation}
		\cutn{U_{d,i}-(U_{d,i})^{\Join\mathcal{P}}}\leq\sqrt{\tfrac{2^{d+1}}{\log(m)}}+\tfrac{2}{m}\;.\label{eq:regularity}
	\end{equation}
	This defines a step Latinon $L^d_{\mathcal{P}}=\left(U_1,id\right)$ where $U_1=\AntiCompr{((U_{d,i})^{\Join\mathcal{P}})_{i\in[2^d]}}$. Now, let $R,C\in[0,1]^n_<$ be $k$-sets chosen uniformly at random. We write $L[R,C]$ for the semilatinon $(U[R,C],id)$, $L^d[R,C]$ for the semilatinon $(\mathbf{U}^d[R,C],id)$ and $L^d_{\mathcal{P}}[R,C]$ for the semilatinon $(U_2,id)$, where $U_2=U_1[R,C]$. Furthermore, let $A$ be the pattern generated by sampling $A_{i,j}$ from $U(r_i,c_j)$ for each $i,j\in[k]$. We then have with probability at least $$1-(2^d\cdot4\exp(-\sqrt{k}/10)+20\exp(-k^{1/8}/10)+\exp(-5k))>1-29\exp(-k^{1/8}/10)$$ that  Claims~\ref{clm:IBoughtASweatshirtToday}, \ref{clm:partdsmpl2}, \ref{clm:partdsmpl1} (stated and proven below) and Lemma~\ref{lem:SampPaternFromDistr} hold. Hence with probability at least $1-29\exp(-k^{1/8}/10)$ we have
	\begin{align*}
		\deltaL((U,id),(W^A,id))\leq&\;\deltaL(L,L^d)&(\textrm{Lemma}~\ref{lem:cutdistcomp})\\
		+&\;\deltaL(L^d,L^d_{\mathcal{P}})&(\textrm{Claim}~\ref{clm:IBoughtASweatshirtToday})\\
		+&\;\deltaL(L^d_{\mathcal{P}},L^d_{\mathcal{P}}[R,C])&(\textrm{Claim}~\ref{clm:partdsmpl1})\\
		+&\;\deltaL(L^d_{\mathcal{P}}[R,C],L^d[R,C])&(\textrm{Claim}~\ref{clm:partdsmpl2})\\
		+&\;\deltaL(L^d[R,C],L[R,C])&(\textrm{Lemma}~\ref{lem:cutdistcomp})\\
		+&\;\deltaL(L[R,C],(W^A,id))&(\textrm{Lemma}~\ref{lem:SampPaternFromDistr})\\
		\leq&\; \tfrac{1}{2^{d-2}} + 2^{d+1}\left(\sqrt{\tfrac{2^{d+1}}{\log(m)}}+\tfrac{2}{m}\right)+\tfrac{4}{k^{1/16}} + \tfrac{8\cdot2^d}{k^{1/4}} +\tfrac{6}{k^{1/4}}&\\
		\leq&\;\tfrac{29}{\log(k)^{1/8}}\;.&
	\end{align*}
	
	\begin{claim}\label{clm:IBoughtASweatshirtToday} We have
		\[\deltaL(L^d,L^d_{\mathcal{P}}) \le 2^d\cdot\left(\sqrt{\tfrac{2^{d+1}}{\log(m)}}+\tfrac{2}{m}\right)\;.\]
	\end{claim}
	\begin{proof}
		The Latinons $L^d$ and $L^d_{\mathcal{P}}$ share the same order-fixing bigraphon $O^f$. Thus, their distance can be bounded using Proposition~\ref{prop:boundbycompr} with~\eqref{eq:regularity} for each $i\in [2^d]$.
	\end{proof}
	
	\begin{claim}\label{clm:partdsmpl2} We have with probability at least $1-2^d\cdot4\exp(-\sqrt{k}/10)$ that
		\[\deltaL(L^d_{\mathcal{P}}[R,C],L^d[R,C])\leq\tfrac{8\cdot 2^d}{k^{1/4}}+2^d\left(\sqrt{\tfrac{2^{d+1}}{\log(m)}}+\tfrac{2}{m}\right)\;.\]
	\end{claim}
	
	\begin{proof}
		Let us fix an $i\in[2^d]$. By Lemma~\ref{lem:firstsamplingordered} we have with probability at least $1-4\exp(-\sqrt{k}/10)$ that
		\[\left|\cutn{U_{d,i}[R,C]-(U_{d,i})^{\Join\mathcal{P}}[R,C]}-\cutn{U_{d,i}-(U_{d,i})^{\Join\mathcal{P}}}\right|\leq\tfrac{8}{k^{1/4}}\;,\]
		and hence together with~\eqref{eq:regularity} that
		\begin{align*}
			&\cutn{U_{d,i}[R,C]-(U_{d,i})^{\Join\mathcal{P}}[R,C]}\\
			\leq&\left|\cutn{U_{d,i}[R,C]-(U_{d,i})^{\Join\mathcal{P}}[R,C]}-\cutn{U_{d,i}-(U_{d,i})^{\Join\mathcal{P}}}\right|+\cutn{U_{d,i}-(U_{d,i})^{\Join\mathcal{P}}}\\
			\leq&\tfrac{8}{k^{1/4}}+\sqrt{\tfrac{2^{d+1}}{\log(m)}}+\tfrac{2}{m}\;.
		\end{align*}
		Since this holds for any $i\in[2^d]$ we get by the union bound and Proposition~\ref{prop:boundbycompr} that with probability at least $1-2^d\cdot 4\exp(-\sqrt{k}/10)$ we have
		\[\deltaL(L^d_{\mathcal{P}}[R,C],L^d[R,C])\leq\sum_{i=1}^{2^d}\cutn{U_{d,i}[R,C]-(U_{d,i})^{\Join\mathcal{P}}[R,C]}\leq\tfrac{8\cdot 2^d}{k^{1/4}}+2^d\left(\sqrt{\tfrac{2^{d+1}}{\log(m)}}+\tfrac{2}{m}\right)\;.\]
	\end{proof}
	
	\begin{claim}\label{clm:partdsmpl1}We have with probability at least $1-20\exp(-k^{1/8}/10)$ that
		\[\deltaL(L^d_{\mathcal{P}},L^d_{\mathcal{P}}[R,C])\leq\tfrac{4}{k^{1/16}}\;.\]
	\end{claim}
	\begin{proof}
		Recall that $L^d_{\mathcal{P}}=(U_1,id)$ and $L^d_{\mathcal{P}}[R,C]=(U_2,id)$. Furthermore we write $R=(r_1,\dots,r_k)$ and $C=(c_1,\dots,c_k)$. 
		For $r_i\in R$ and $c_i\in C$ we assign an interval of measure $1/k$ by setting $\interval{k}{r_i}\coloneqq \interval{k}{c_i}\coloneqq \interval{k}{i}$. 
		Also we want to count how many elements of $R$ fall into the intervals of $\partition{r}$ and their refinement $\mathcal{P}$ by setting $R_{i}\coloneqq R\cap 
		\interval{r}{i}
		$ and $R_{i,j}\coloneqq R\cap P_{i,j}$. Similarly, we write $C_{i}\coloneqq C\cap 
		\interval{r}{i}
		$ and $C_{i,j}\coloneqq C\cap P_{i,j}$. Finally we write $R_{i,j}'\coloneqq \bigcup_{x\in R_{i,j}}
		\interval{k}{x}
		$ and $C_{i,j}'\coloneqq \bigcup_{x\in C_{i,j}}
		\interval{k}{x}$. The key observation is that by the definition of the sampling process the two distribution-valued bigraphons $U_1$ and $U_2$ are very similar. Both are step Latinons and $U_1$ takes the same value on $P_{i,j}\times P_{s,t}$ as $U_2$ on $R_{i,j}'\times C_{s,t}'$. The only difference lies in the size and the position of the steps. It will turn out that with high probability most of $R_{i,j}'$ and $C_{i,j}'$ actually lies in $\interval{r}{i}$ and is of similar measure as $P_{i,j}$. Hence a suitable reordering inside each $\interval{r}{i}$ matches $U_1$ to $U_2$ up to a small error set of $[0,1]^2$.
		
		We first observe that with high probability the chosen elements of $R$ distribute among the partite sets $P_{i,j}$ as expected. By Lemma~\ref{lem:chernoff} with $\eps=1/k^{5/16}$ and $\mu=\lambda(P_{i,j})\cdot k=k/(mr)=k^{3/4}$ we have with probability at least $(1-rm\cdot\exp(-k^{1/8}/3))>1-10\exp(-k^{1/8}/10)$ that
		\begin{equation}
			|R_{i,j}|=\lambda(P_{i,j})\cdot k\pm k^{11/16}\lambda(P_{i,j})=k^{3/4}\pm k^{7/16}\textrm{ for every }i\in[r],j\in[m]\;.\label{eq:ClaraZetkin}
		\end{equation}
		Now, we can observe that most of $R_{i,j}'$ is contained in 
		$\interval{r}{i}$
		. In the worst case each of the at most $rm$ intervals before $R_{i,j}'$ are too short (or too long) by the maximum size, by \eqref{eq:ClaraZetkin} i.e. $rm\cdot k^{7/16}\cdot 1/k$, and hence
		\[\lambda(R_{i,j}'\cap 
		\interval{r}{i}
		)=\lambda\left(\left(\bigcup_{x\in R_{i,j}}
		\interval{k}{x}
		\right)\cap 
		\interval{r}{i}
		\right)\geq\lambda(P_{i,j})- rm/ k^{9/16}=1/k^{1/4}-1/k^{5/16}\;.\]
		Thus for each $i,j$ we can choose sets $A_{i,j}\subseteq R_{i,j}'\cap 
		\interval{r}{i}
		$ and $B_{i,j}\subseteq P_{i,j}\subseteq 
		\interval{r}{i}
		$ of measure $1/k^{1/4}-1/k^{5/16}$ and a bijection $\vphi_{i,j}:A_{i,j}\rightarrow B_{i,j}$. Let $\vphi:[0,1]\rightarrow[0,1]$ be defined by 
		\begin{equation*}
			x\mapsto\begin{cases}\vphi_{i,j}(x),\textrm{ if }x\in A_{i,j},\\x,\textrm{ otherwise.}\end{cases}
		\end{equation*}
		Note that then, since $\vphi$ only permutes elements inside each 
		$\interval{r}{i}$
		, that 
		\[\cutn{O-O^{\vphi}}\leq r\cdot 1/r^2=1/r=1/k^{1/8}\;.\]
		With probability at least $1-10\exp(-k^{1/8}/10)$ we can repeat the same procedure for $C$ to get sets $A_{i,j}'\subseteq C_{i,j}'\cap 
		\interval{r}{i}
		$ and $B_{i,j}'\subseteq P_{i,j}\subseteq 
		\interval{r}{i}
		$ of measure $1/k^{1/4}-1/k^{5/16}$ and a bijection $\psi$ defined in an analogous way such that
		\[\cutn{O-O^{\psi}}\leq 1/k^{1/8}\;.\]
		Note that then $U_1^{\vphi,\psi}(x,y)=U_2(x,y)$ for every $x\in \bigcup_{i\in[r],j\in[m]}A_{i,j}$ and $y\in \bigcup_{i\in[r],j\in[m]}A_{i,j}'$. Thus $U_1^{\vphi,\psi}$ and $U_2$ agree everywhere on $[0,1]^2$ except on a set of measure at most $2\cdot rm\cdot 1/k^{5/16}=2/k^{1/16}$ and we have
		\[\cutn{U_1^{\vphi,\psi}-U_2}\leq 2/k^{1/16}\;.\]
		Hence in total with probability at least $1-20\exp(-k^{1/8}/10)$
		\[\deltaL((U_1,id),(U_2,id))\leq \tfrac{2}{k^{1/16}}+\tfrac{2}{k^{1/8}}<\tfrac{4}{k^{1/16}}\;.\]
	\end{proof}
	
\end{proof}

We can now prove Lemma~\ref{lem:sampling}.

\begin{proof}[Proof of Lemma~\ref{lem:sampling}]
	Let $k\in\mathbb N$. By Lemma~\ref{lem:couplsampl} with $\eps=\exp(-k)/k^2$ there exists a Latinon $L'=(U,id)$ on $[0,1]^2$ such that
	\[\deltaL^*(L,L')<\exp(-k)/k^2\;,\]
	and a coupling $\mathcal{C}$ of sampling a $k\times k$ pattern $A$ from $L$ and a $k\times k$ pattern $B$ from $L'$ such that
	\[\Prob_{\mathcal{C}}[A\neq B]<\eps k^2=\exp(-k)\;.\]
	By Lemma~\ref{lem:sample01} we then have with probability at least $$1-29\exp(-k^{1/8}/10)-\exp(-k)>1-30\exp(-k^{1/8}/10)$$ that
	\[\deltaL^*(L,(W^A,id))\leq\deltaL^*(L,L')+\deltaL^*(L',(W^B,id))\leq \exp(-k)+\tfrac{29}{\log(k)^{1/8}}<\tfrac{30}{\log(k)^{1/8}}\;.\]
\end{proof}

\section{Proof of the inverse counting lemma}\label{sec:invcount}
In this section we prove the inverse counting lemma, Lemma~\ref{lem:invcounting}. To this end, we need to draw a connection between sampled realisations and patterns. The following definition will be useful.

\begin{defi}\label{def:spread}
	Given $\eps>0$, we say that a finite multiset $S\subset[0,1]$ is \emph{$\eps$}-spread if for each interval $I\subset [0,1]$ we have that $|I\cap S|=(\lambda(I)\pm\eps)|S|$.
\end{defi}

Our first observation is that sampled realisations are almost always spread.

\begin{lem}\label{lem:samplingspread}
	Let $L$ be a Latinon on ground space $\Omega$, and $k\in\mathbb{N}$. Let $S$ be the (random) multiset of $k^2$ values obtained from $L$ by $k\times k$-sampling. Then we have with probability at least $1-\exp(-k^{0.1})$ that $S$ is $(4k^{-0.4})$-spread.
\end{lem}

\begin{proof}
	Let us fix an arbitrary set $Y\subset [0,1]$. Let $(i,j)\in [k]^2$ be arbitrary. Now, consider the sampling from $L$ for the value $v$~corresponding to the $i$-th sampled row and $j$-th sampled column.\footnote{Notice that we are really talking about the row that was sampled as the $i$-th, not about the row that has the $i$-th smallest global position. Same for the column.} By the uniform marginals property we have 
	\[\Prob[v\in Y]=\lambda(Y)\;. \]
	Thus, by the linearity of expectation we have for the multiset of all $k^2$ sampled values that
	\[\Expect[|S\cap Y|]=\lambda(Y)k^2\;. \]
	In the rest of the proof, we argue that this value is concentrated. To this end we want to use McDiarmid's inequality, Lemma~\ref{lem:MBoundDiff}. We need to set up a suitable product probability space. There is an option which reflects our sampling. The first $k$ coordinates will be copies of $\Omega$ representing the selection of rows, the next $k$ coordinates will be copies of $\Omega$ representing the selection of columns, and the last $k^2$ coordinates should model the choice of the individual entries of the sampled matrix $M$. We have to be careful because these last probability spaces need to be fixed from the beginning in the setting of McDiarmid's inequality. Let $\Lambda=\Omega^k\times\Omega^k\times [0,1]^{k\times k}$. For $(i,j)\in [k]^2$, and $\mathbf{w}\in\Lambda$, let $q_{ij}$ be defined as follows. Take $\mathcal D_{ij}\coloneqq L(\mathbf{w}_i,\mathbf{w}_{k+j})$. That is, $\mathcal D_{ij}$ is the distribution on $[0,1]$ at the $i$-th sampled row and the $j$-th sampled column. Let $F_{\mathcal D_{ij}}$ be the cumulative distribution function of $\mathcal D_{ij}$. Now, let $q_{ij}$ be the infimum of $t$'s for which $F_{\mathcal D_{ij}}(t)> \mathbf{w}_{k+i+kj}$. This construction ensures that the matrix $(q_{ij})_{i,j}$ has the same distribution as the matrix sampled from $L$ (before reordering of rows and columns).
	
	We can now define a map $Z:\Lambda\rightarrow \mathbb{R}$, $Z(\mathbf{w})\coloneqq \sum_{i,j} \mathbbm{1}_{q_{ij}\in Y}$. The above justifies that $Z$ satisfies the requirements of Lemma~\ref{lem:MBoundDiff} with $\mathbf{c}$ being $k$ on the first $2k$ coordinates and being~$1$ on the remaining $k^2$ coordinates. Thus, Lemma~\ref{lem:MBoundDiff} tells us that for each $t>0$ we have
	\begin{equation}\label{eq:OrangeBlossomSpecial}
		\Prob\left[\left|Z-\lambda(Y)k^2\right|>t\right]\le\exp\left(-\frac{2t^{2}}{2k\cdot k^2+k^2\cdot 1^2}\right)\le \exp\left(-\frac{2t^{2}}{3k^3}\right)\;.
	\end{equation}
	
	Now, let us take $t\coloneqq \lceil k^{1.6}\rceil$, and let us apply the above for intervals $Y_{\ell_1,\ell_2}\coloneqq \left[\frac{\ell_1}t,\frac{\ell_2}t\right]$, $0\le \ell_1<\ell_2\le t$, $\ell_1,\ell_2\in\mathbb{N}$. 
	The union bound over~\eqref{eq:OrangeBlossomSpecial} tells us that
	\begin{align}
		\begin{split}\label{eq:Iamsogoodevent}
			\Prob\left[\mbox{for each $\ell_1$ and $\ell_2$ we have }||S\cap Y_{\ell_1,\ell_2}|-\lambda(Y_{\ell_1,\ell_2})k^2|\le \lceil k^{1.6}\rceil\right] &\ge 1-t^2\exp\left(-\frac{2t^{2}}{3k^3}\right)\\
			&\ge 1-\exp(-k^{0.1})\;.
		\end{split}
	\end{align}
	Now, we claim that if the good event from~\eqref{eq:Iamsogoodevent} is satisfied, then $S$ is $(40k^{-0.5})$-spread. Indeed, let $I\subset [0,1]$ be an arbitrary interval. We now choose a minimal interval $Y_{\ell_1^+,\ell_2^+}$ which contains $I$. Note that $\lambda(Y_{\ell_1^+,\ell_2^+})\le \lambda(I)+\frac{2}{\lceil k^{1.6}\rceil}$. Then we have 
	\begin{align*}
		|S\cap I|&\le |S\cap Y_{\ell_1^+,\ell_2^+}|
		\le \lambda(Y_{\ell_1^+,\ell_2^+})k^2+\lceil k^{1.6}\rceil
		\le \left(\lambda(I)+\frac{2}{\lceil k^{1.6}\rceil}\right)k^2+\lceil k^{1.6}\rceil \\
		&\le \left(\lambda(I)+4k^{-0.4}\right)k^2\;.
	\end{align*}
	Similarly, taking a maximal interval $Y_{\ell_1^-,\ell_2^-}$ which is contained in $I$, we can obtain that $|S\cap I|\ge\left(\lambda(I)-4k^{-0.4}\right)k^2$. This concludes the proof.
\end{proof}

Our second observation is that realisations which are spread and belong to the same pattern are close in the cut norm. 

\begin{lem}\label{lem:spreaddist}
	Let $A\in\mathcal{R}(k,k)$ be a $k\times k$ pattern and $A^*\in\mathcal{R}^A([0,1])$ be a realisation. Let $L^A=(W^A,id)$ and $L^{A^*}=(W^{A^*},id)$ be the associated semilatinons with ground space $[0,1]$. If $A^*$ is $\eps$-spread, then
	\begin{equation}\label{eq:LALAs}
		\cutnL{W^A-W^{A^*}}\leq 4\eps+\tfrac{4}{k^2}\;.
	\end{equation}
\end{lem}

\begin{proof}
	Observe that $A_{i,j}=|\{A^*\}\cap[0,A^*_{i,j}]|=k^2(A^*_{i,j}\pm\eps)$ for each fixed $i,j\in [k]$, where the first equality holds because of $A^*\in\mathcal{R}^A([0,1])$ and the second, as $A^*$ is $\eps$-spread. 
	
	In order to verify~\eqref{eq:LALAs} using Definition~\ref{def:cutnormlati}, let $R,C\subseteq[0,1]$ be arbitrary and $V\coloneqq[a,b]\subseteq [0,1]$ an arbitrary interval. Set $Y\coloneqq\{(i,j)\in[k]^2\mid \;k^{-2}A_{i,j}\in V,A_{i,j}^*\notin V\}$ and $Y^*\coloneqq\{(i,j)\in[k]^2\mid A_{i,j}^*\in V,\;k^{-2}A_{i,j}\notin V\}$. 
	Since $A_{i,j}^*= k^{-2}A_{i,j} \pm \eps$, the set $Y'\coloneqq\{(i,j)\in[k]^2\mid \;k^{-2}A_{i,j} \geq a,A_{i,j}^*<a\}$ has size at most $\eps k^2+1$. Similarly the set $Y''\coloneqq\{(i,j)\in[k]^2\mid \;k^{-2}A_{i,j} \leq b,A_{i,j}^*>b\}$ has size at most $\eps k^2+1$. Then $|Y| \leq |Y'|+|Y''| \leq 2\eps k^2+2$. We similarly obtain $|Y^*| \leq 2\eps k^2+2$.
	
	Therefore
	\[\left|\int_{x\in R}\int_{y\in C}W^A(x,y)[V]-W^{A^*}(x,y)[V]\diffsymb x\diffsymb y\right|\leq \tfrac{1}{k^2}(|Y|+|Y^*|) \leq \tfrac{1}{k^2}(4\eps k^2+4) \leq 4\eps+\tfrac{4}{k^2}\;.\]
\end{proof}

We now have all the tools to prove Lemma~\ref{lem:invcounting}.
\begin{proof}[Proof of Lemma~\ref{lem:invcounting}]
	We set $\eps\coloneqq 4k^{-0.4}$. For patterns $A\in\mathcal{R}(k,k)$ and realisations $A^*\in[0,1]^{k\times k}$ we will work with the associated semilatinons $L^A=(W^A,id)$ and $L^{A^*}=(W^{A^*},id)$ over the ground space $[0,1]$. Furthermore we set
	\begin{align*}
		\mathcal{A}_1&\coloneqq \{A\in\mathcal{R}(k,k)\mid \deltaL^*(L^A,L_1)<d/2\}\;,\\
		\mathcal{A}_2&\coloneqq \{A\in\mathcal{R}(k,k)\mid \deltaL^*(L^A,L_2)<d/2\}\;,\\
		\mathcal{A}_1^{*,\eps}&\coloneqq \{A^*\in[0,1]^{k\times k}\mid A^*\textrm{ is $\eps$-spread and $\deltaL^*(L^{A^*},L_1)<d/4$}\}\;,\\
		\mathcal{A}_2^{*,\eps}&\coloneqq \{A^*\in[0,1]^{k\times k}\mid A^*\textrm{ is $\eps$-spread and $\deltaL^*(L^{A^*},L_2)<d/4$}\}\;.
	\end{align*}
	
	Note that the triangle inequality for $\deltaL^*$ implies that $\mathcal{A}_1\cap\mathcal{A}_2=\emptyset$. Furthermore, as the choice of $\eps$ and the condition on the size of $k$ guarantee that $4\eps+4/k^2<d/4$, we have by Lemma~\ref{lem:spreaddist} that the following holds. If $A\in\mathcal{R}(k,k)$ and $A^*\in\mathcal{R}^A([0,1])$ is $\eps$-spread and such that $\deltaL^*(L^{A^*},L_1)<d/4$, then $\deltaL^*(L^{A},L_1)< 4\eps+4/k^2+d/4\leq d/2$. Therefore we have that for every $A^*\in\mathcal{A}_1^{*,\eps}$ there exists $A\in\mathcal{A}_1$ such that $A^*\equiv A$. Hence we have for uniformly chosen $<_f$-ordered $k$-sets $S,T\subseteq \Omega$ that
	\[\sum_{A\in\mathcal{A}_1}t(A,L_1)\geq \mathbb{P}[L[S,T]\in\mathcal{A}_1^{*,\eps}]\;,\]
	and, as the size of $k$ guarantees $\frac{30}{\log(k)^{1/10}}<d/4$, it follows from Lemma~\ref{lem:sampling} and Lemma~\ref{lem:samplingspread} that
	\[\mathbb{P}[L[S,T]\in\mathcal{A}_1^{*,\eps}]\geq 1-\exp(-k^{0.1})-\exp(-4k/\log(k)^{0.25})\geq \tfrac{3}{4}\;.\]
	Similarly we have
	\[\sum_{A\in\mathcal{A}_2}t(A,L_2)\geq \mathbb{P}[L[S,T]\in\mathcal{A}_2^{*,\eps}]\geq \tfrac{3}{4}\;,\]
	and therefore
	\[\sum_{A\in\mathcal{A}_1}t(A,L_2)\leq 1-\sum_{A\in\mathcal{A}_2}t(A,L_2)\leq \tfrac{1}{4}\;.\]
	Thus
	\[\sum_{A\in\mathcal{A}_1}(t(A,L_1)-t(A,L_2))\geq\tfrac{1}{2}\]
	and by averaging, and since $|\mathcal{A}_1|\leq k^2!$, there exists a pattern $A\in\mathcal{A}_1$ such that
	\[|t(A,L_1)-t(A,L_2)|\geq \left(\sum_{A\in\mathcal{A}_1}t(A,L_1)-\sum_{A\in\mathcal{A}_1}t(A,L_2)\right)/(k^2!)\geq\tfrac{1}{2(k^2!)}\;.\]
\end{proof}

\section{Approximating Latinons}\label{sec:ApproximatingLatinons}

In this section, we prove Theorem~\ref{thm:latinonsarelimits}. In fact,  Theorem~\ref{thm:latinonsarelimits} is a direct corollary of the lemma below, which we prove instead.
\begin{lem}\label{lem:approx}
	Let $L$ be a Latinon on $\Omega$. For every $\eps>0$ and $n_0\in\mathbb N$ there exists $n\ge n_0$ and a Latin square $L_n$ of order $n$ with the property that
	\[\deltaL^*(L,L_n)<\eps\;.\]
\end{lem}

To find such a Latin square we first approximate the given Latinon $L$ by a step Latinon.

\begin{defi}
	A step Latinon $(U,f)$ on $\Omega$ is a Latinon such that there exists a partition $\mathcal{P}=(P_1,\dots,P_m)$ of $\Omega$ and a partition $\mathcal{Q}=(Q_1,\dots,Q_d)$ of $[0,1]$ consisting of consecutive intervals such that there exist $\alpha_{i,j,k}\in\mathbb [0,1]$, $i,j\in[m],k\in[d]$, such that for every $I\subseteq Q_k$,
	\[U(x,y)[I]=\alpha_{i,j,k}\frac{\lambda(I)}{\lambda(Q_k)}\;\textrm{ for every }(x,y)\in P_i\times P_j\;.\]
	In other words, the distribution of $U$ is constant on each $P_i\times P_j$, and consists of a linear combination (with coefficients between $0$ and $1$) of Lebesgue measures on $\{Q_k\}_{k\in [d]}$.
\end{defi}

Our first technical lemma says that each Latinon can be approximated by a step Latinon. For technical reasons, we wish to switch in this approximation from the ground space $\Omega$ to ground space $[0,1]$. We leave the proof of the lemma for later.
\begin{lem}\label{lem:stepapprox}
	Let $L$ be a Latinon on $\Omega$. For every $\eps>0$ there exists a step Latinon $(U,id)$ on $[0,1]$ such that its corresponding partitions $\mathcal{P}$ and $\mathcal{Q}$ consist of intervals and
	\[\deltaL^*(L,(U,id))<\eps\;.\]
\end{lem}

Then we find a Latin square which is close to the previously constructed step Latinon.

\begin{lem}\label{lem:LSapprox}
	Let $\eps>0$ and $n_0\in\mathbb N$. If $(U,id)$ is a step Latinon on $[0,1]$ whose corresponding partitions consist of intervals, then there exists $n\ge n_0$ and a Latin square $L_n$ of order $n$ with the property that
	\[\deltaL((U,id),L_n)<\eps\;.\]
\end{lem}

Using both lemmas the proof of Lemma~\ref{lem:approx} then follows immediately.

\begin{proof}[Proof of Lemma~\ref{lem:approx}]
	Given $\eps>0$, we use Lemma~\ref{lem:stepapprox} with $\eps/2$ in place of $\eps$ to get a step Latinon $(U,id)$ on $[0,1]$ with $\deltaL(L,(U,id))<\eps/2$. We then use Lemma~\ref{lem:LSapprox} again with $\eps/2$ in place of $\eps$ to get a number $n\ge n_0$ and a Latin square $L_n$ of order $n$ such that $\deltaL((U,id),L_n)<\eps/2$. Together we therefore get
	\[\deltaL^*(L,L_n)\leq\deltaL^*(L,(U,id))+\deltaL((U,id),L_n)<\eps/2+\eps/2=\eps\;.\]
\end{proof}

We first show the proof of Lemma~\ref{lem:stepapprox} and defer the proof of Lemma~\ref{lem:LSapprox} to Section~\ref{sec:constrlatinsquares}. For this we introduce the phrase \emph{for $\alpha$ sufficiently smaller than $\beta$}, written $\alpha\ll\beta$, meaning that for any $\beta>0$ there exists $\alpha_0>0$ such that for any $\alpha\leq\alpha_0$ the subsequent statement holds.

\begin{proof}[Proof of Lemma~\ref{lem:stepapprox}]
	First we choose constants $d,r,m\in\mathbb N$ such that
	\[\tfrac{1}{r}\ll\tfrac{1}{d},\tfrac{1}{m}\ll\eps\;.\]
	
	Let $L=(W,f)$ be the given Latinon. Let $L^d=(O^f,W_{d,1},\dots,W_{d,2^d})$ be the compression of $L$ of depth $d$. 
	Consider the partition $\partition{m}$.
	By Lemma~\ref{lem:weak_reg_tup}\ref{rl:part1} we find a partition $\mathcal{P}=\{P_{j,t}\mid j\in[m],t\in[r]\}$ of $\Omega$ refining $f^{-1}(
	\partition{m}
	)$, i.e. $f^{-1}(
	\interval{m,j}
	)=P_{j,1}\cup\cdots\cup P_{j,r}$, $j\in[m]$, and such that, for every $1\leq i\leq 2^d$,
	\begin{equation}
		\cutn{W_{d,i}-(W_{d,i})^{\Join\mathcal{P}}}\leq\sqrt{\tfrac{2^{d+1}}{\log(r)}}\;.\label{eq:regular}
	\end{equation}
	Note that this defines a step Latinon $L^d_{\mathcal{P}}=(W',f)$ with $W'=\AntiCompr{(W_{d,1})^{\Join\mathcal{P}},\dots,(W_{d,2^d})^{\Join\mathcal{P}}}$ on $\Omega$ and we get by \textrm{Lemma}~\ref{lem:cutdistcomp}, \textrm{Proposition}~\ref{prop:boundbycompr}, and \eqref{eq:regular} that
	\begin{equation}\label{eq:apprstep1}
		\deltaL(L,L^d_{\mathcal{P}})\leq\deltaL(L,L^d)+\deltaL(L^d,L^d_{\mathcal{P}})\leq \tfrac{1}{2^{d-1}}+2^d\cdot\sqrt{\tfrac{2^{d+1}}{\log(r)}}<\eps/2\;.
	\end{equation}
	
	We now fix arbitrarily for every $j\in[m]$ another partition $\mathcal{Q}_j=\{Q_{j,1},\dots,Q_{j,r}\}$ of 
	$\interval{m}{j}$
	such that $Q_{j,t}$ is an interval and $\lambda(Q_{j,t})=\mu(P_{j,t})$ for every $t\in[r]$. Note that this is possible, as $\sum_{t\in[r]}\mu(P_{j,t})=\lambda(
	\interval{m}{j}
	)$ for $j\in[m]$. Then $\mathcal{Q}=\{Q_{j,t}\mid j\in[m],t\in[r]\}$ is a partition of $[0,1]$ consisting of intervals. We now use this partition to define a step Latinon $(U,id)$ on $[0,1]$ by setting for every $x\in Q_{i,s},y\in Q_{j,t}$ and $I\subseteq[0,1]$ that $U(x,y)(I)\coloneqq W'(x',y')(I)$ for some $x'\in P_{i,s}$ and $y'\in P_{j,t}$. Note that since $L^d_{\mathcal{P}}$ is a step Latinon this definition does not depend on the choice of $x',y'$ and also consequently $(U,id)$ is a step Latinon. 
	
	Next, we observe that $(U,id)$ is close to $L^d_{\mathcal{P}}$ in the generalised cut distance $\deltaL^*$ (see Section~\ref{sec:diffgroundspaces}) which can compare Latinons defined on different ground spaces.
	\begin{claim}\label{clm:dist01} We have
		$$\deltaL^*(L^d_{\mathcal{P}},(U,id))<\eps/2.$$
	\end{claim}
	
	\begin{proof}[Proof of Claim~\ref{clm:dist01}]
		For every $j\in[m],t\in[r]$ let $\vphi_{j,t}:Q_{j,t}\rightarrow P_{j,t}$ be an arbitrary measure preserving bijection. These bijections altogether naturally define a measure preserving bijection $\vphi:[0,1]\rightarrow\Omega$. Note that then $(W')^{\vphi,\vphi}(x,y)=U(x,y)$ for every $x,y\in[0,1]$ by definition of $U$. Furthermore, we also have $\vphi(
		\interval{m}{k}
		)=f^{-1}(
		\interval{m}{k}
		)$ and therefore $(f\circ\vphi)(
		\interval{m}{k}
		)=
		\interval{m}{k}
		$ for every $k\in[m]$. Thus
		\[\deltaL^*((U,id),L^d_{\mathcal{P}})\leq 2\cutn{O-O^{f\circ\vphi }}+\cutnL{U-(W')^{\vphi,\vphi}}\leq 2\sum_{k\in[m]}\lambda(
		\interval{m}{k}
		)^2= \tfrac{2}{m}<\eps/2\;.\]
	\end{proof}
	
	We therefore have by \eqref{eq:apprstep1}, Claim~\ref{clm:dist01}, and the generalised triangle inequality from Section~\ref{sec:diffgroundspaces} that
	$$\deltaL^*(L,(U,id))\leq\deltaL(L,L^d_{\mathcal{P}})+\deltaL^*(L^d_{\mathcal{P}},(U,id))<\eps/2+\eps/2=\eps\;.$$

\end{proof}

\subsection{Constructing Latin squares}\label{sec:constrlatinsquares}
It is more convenient to represent the sought Latin square $L_n$ as in Lemma~\ref{lem:LSapprox} as a $K_3$-decomposition of the complete tripartite graph $K_{n,n,n}$. Indeed, given any Latin square $L_n$ of order $n$, we view the three parts of $K_{n,n,n}$ as numbering rows, columns, and values, respectively. Then each entry of $L_n$ corresponds to a triangle in $K_{n,n,n}$, and the collection of all the entries corresponds to a \emph{$K_3$-decomposition} of $K_{n,n,n}$, that is a partition of $E(K_{n,n,n})$ into triangles. In fact, it is easy to see that once the role of the three parts of $K_{n,n,n}$ and their numberings are fixed, there is a one-to-one correspondence between Latin squares of order $n$ and $K_3$-decompositions of $K_{n,n,n}$.

We will need to recall some notation and facts from the theory of designs, mainly following  Keevash's breakthrough results~\cite{Kee:DesignsI,Kee:DesignsII}.
As we explained above, it suffices to recall these results in the special case of $K_3$-decompositions of tripartite graphs. Given a graph $G$, a vertex $v\in V(G)$, and $S,V_1,V_2\subseteq V(G)$ where $V_1\cap V_2=\emptyset$, we write $N_G(v,S)=N_G(v)\cap S$ for the neighbourhood of $v$ in $S$. We also write $\overline{N}_G(S)=\bigcap_{w\in S}N_G(w)$ for the common neighbourhood of $S$ and $d_G(V_1,V_2)=\frac{|E(G[V_1,V_2])|}{|V_1||V_2|}$ for the density of the pair $(V_1,V_2)$.

\begin{defi}
	Let $G$ be a tripartite graph with partite sets $V_1,V_2,V_3$ and densities $d_{\{i,j\}}=\tfrac{|E(G[V_i,V_{j}])|}{|V_i||V_{j}|}$, $\{i,j\}\in\binom{[3]}{2}$.
	\begin{enumerate}[label={(\roman*)}]
		\item We say that $G$ is (partite) \emph{$K_3$-decomposable} if there exists an edge-disjoint collection of copies of $K_3$ in $G$ which covers $E(G)$.
		\item We say that $G$ is \emph{$K_3$-balanced} if for every $\{i,j,k\}=[3]$ and every $v\in V_k$, $\deg(v,V_i)=\deg(v,V_j)$.
		\item Let $c>0$ and $h\in\mathbb N$. We say that $G$ is \emph{$(c,h)$-typical}, if for every pair of distinct $i,j\subseteq[3]$, and every set $S\subseteq V_i\cup V_j$ with $|S|\leq h$ we have that $|\overline{N}_G(S)\cap V_k|=(1\pm c)d_{\{i,k\}}^{|S\cap V_i|}d_{\{j,k\}}^{|S\cap V_j|}|V_k|$, where $k\in[3]\setminus\{i,j\}$.
	\end{enumerate}
\end{defi}

It is easy to see that a $K_3$-decomposable tripartite graph needs to be $K_3$-balanced. Now, the important result from the theory of designs we want to use is that for typical graphs the converse is also true. The first theorem we need is that a typical graph can be approximately decomposed such that the leftover graph has bounded degree. The existence of an approximate decomposition was shown first by R\"odl in~\cite{Rodl:Nibble}. For our purposes we need a slightly stronger version, i.e. the additional degree condition for the leftover graph. This follows from a more general theorom of Alon and Yuster in~\cite{AY05}. For completeness we explain how to derive our version in Appendix~\ref{app:nibble}.

\begin{thm}[\cite{Rodl:Nibble}]\label{thm:LSalmost} 
	Suppose that $\tfrac{1}{n}\ll c\ll\alpha\ll d\leq 1$.
	Let $G$ be a tripartite graph with partite sets $V_1,V_2,V_3$ such that $|V_1|=|V_2|=|V_3|=n$ and $d_G(V_i,V_j)=(1\pm c)d$ for all $\{i,j\}\subseteq[3]$.
	If $G$ is $(c,2^{4500})$-typical, then there exists a subgraph $S\subseteq E(G)$ such that
	\begin{enumerate}
		\item $G-S$ has a partite $K_3$-decomposition, and
		\item $\Delta(S)\leq\alpha n$.
	\end{enumerate}
\end{thm}

The second theorem we need is that typical $K_3$-balanced tripartite graphs are $K_3$-decomposable and is a direct application of Theorem~1.7 from~\cite{Kee:DesignsII}.

\begin{thm}[\cite{Kee:DesignsII}]\label{thm:LScomplete}
	Suppose that $\tfrac{1}{n}\ll c\ll d\leq 1$.
	Let $G$ be a $K_3$-balanced tripartite graph with partite sets $V_1,V_2,V_3$ such that $|V_1|=|V_2|=|V_3|=n$ and $d_G(V_i,V_j)\geq d$ for all $\{i,j\}\in\binom{[3]}{2}$.
	If $G$ is $(c,2^{4500})$-typical, then $G$ has a partite $K_3$-decomposition. 
\end{thm}

We will also make use of the following simple lemma stating that selecting edges independently and uniformly at random from a typical graph results in a typical subgraph with high probability.

\begin{lem}\label{lem:typrandsubgr} Suppose that $1/n\ll c_1\ll c_2 \ll1/h,\rho, d\leq 1$.
	Let $G$ be a $(c_1,h)$-typical tripartite graph with vertex classes $V_1$, $V_2$, $V_3$ with $|V_1|=|V_2|=|V_3|$ and $d_G(V_i,V_j)=(1\pm c_1)d$ for every $\{i,j\}\in\binom{[3]}{2}$. If we select each edge of $G$ with probability $\rho$ independently and uniformly at random for a subgraph $H$, then with probability $1-\exp(-\Theta(n))$ we have
	\begin{enumerate}[label={(\roman*)}]
		\item\label{x:one} $d_H(V_i,V_j)=(1\pm c_2)\rho d$ for every $\{i,j\}\in\binom{[3]}{2}$,
		\item\label{x:two} $H$ is $(c_2,h)$-typical, and
		\item\label{x:three} $G-H$ is $(c_2,h)$-typical.
	\end{enumerate}
\end{lem}

\begin{proof}
	For every $\{i,j\} \in\binom{[3]}{2}$ we have that $\Expect[|H[V_i,V_j]|]=\rho (1\pm c_1)dn^2$ and by Lemma~\ref{lem:chernoff} we then get that 
	\[\Prob\left[\left||H[V_i,V_j]|-\rho (1\pm c_1)dn^2\right|\geq (c_2/2)\rho dn^2\right]\leq 2\exp(-c_2^2\rho dn^2/12)\;.\]
	Therefore by the union bound~\ref{x:one} holds with probability $1-\exp(-\Theta(n))$. Furthermore, for every $\{i,j,k\}=[3]$, and every set $S\subseteq V_i\cup V_j$ with $|S|\leq h$ we have that \[\Expect\left[|\overline{N}_H(S)\cap V_k|\right]=(1\pm c_1)(\rho(1\pm c_1)d)^{|S|}n\;,\] and again by Lemma~\ref{lem:chernoff} we then get that 
	\[\Prob\left[\left||\overline{N}_H(S)|-(1\pm c_1)(\rho(1\pm c_1)d)^{|S|}n\right|\geq (c_2/2)(\rho d)^{|S|}n\right]\leq 2\exp(-c_2^2 (\rho d)^{|S|}n/12)\;.\]
	Again by using the union bound the probability that $H$ is $(c_2,h)$-typical is at least $$1-h\cdot \binom{3n}{h}\cdot \exp(-c_2^2 (\rho d)^{h}n/12)=1-\exp(-\Theta(n))\;.$$ Similarly, one can show that~\ref{x:three} holds with probability $1-\exp(-\Theta(n))$. Therefore with probability $1-\exp(-\Theta(n))$, \ref{x:one}-\ref{x:three} hold altogether.
\end{proof}

With these tools at hand we can prove Lemma~\ref{lem:LSapprox}. As we explained at the beginning of Section~\ref{sec:constrlatinsquares}, we will find a triangle decomposition of $K_{n,n,n}$ which represents the Latin square we want to construct. We partition the three copies of $[n]$, which represent the rows, columns and values, into \emph{clusters} corresponding to the parts of $[0,1]$ defined by the given step Latinon $(U,id)$. The construction itself has three steps. First, we remove a sparse \emph{reservoir} $R\subset E(K_{n,n,n})$ at random. Second, we find an approximate $K_3$-decomposition\footnote{\label{foot:sparse}By an `approximate $K_3$-decomposition' we mean that a small number of edges of $K_{n,n,n}-R$ will not be covered; here the number of edges will be chosen to be much smaller than that of $R$.} $B$ of $K_{n,n,n}-R$ with the key property that between any triple of clusters (traversing the three colour classes of $K_{n,n,n}$) the number of triangles of $B$ spanned by these three clusters is proportional to the value of $U$ on the steps corresponding to these three clusters. Theorem~\ref{thm:LSalmost} is used to this end. Third, we use Theorem~\ref{thm:LScomplete} to find a $K_3$-decomposition in the graph $K_{n,n,n}-B$. Theorem~\ref{thm:LScomplete} indeed applies in this setting as $K_{n,n,n}-B$ is almost identical to $R$ (recall Footnote~\ref{foot:sparse}) which has excellent typicality properties because of Lemma~\ref{lem:typrandsubgr}. 

\begin{proof}[Proof of Lemma~\ref{lem:LSapprox}] First note that by passing to a refinement we may assume that $(U,id)$ has the same corresponding partition $\mathcal{C}=\{C_1,\dots,C_\ell\}$ for each coordinate. Since $(U,id)$ is a step Latinon we therefore have real values $K_{i,j,k}$, $i,j,k\in[\ell]$, such that $U(x,y)(I)=K_{i,j,k}\lambda(I)/\lambda(C_k)$ for every $x\in C_i$, $y\in C_j$, and $I\subseteq C_k$. Let $\eps>0$, set $h=2^{4500}$, and $K\coloneqq\min \{K_{i,j,k} \mid i,j,k\in[\ell], K_{i,j,k}>0\}$. Choose additional constants such that
	\[\tfrac{1}{n}\ll c_1\ll c_2\ll\alpha\ll c\ll\rho\ll\tfrac{1}{M},\tfrac{1}{q}\ll\tfrac{1}{n_0},\tfrac{1}{\ell}, K,\tfrac{1}{h},\eps,\min_i \lambda(C_i)\;,\]
	and such that $m\coloneqq n/(qM)=2^d$ for some $d\in\mathbb N$.
	
	In the following we want to approximate $(U,id)$ with another step Latinon $(U',id)$ whose corresponding partition fulfills some convenient divisibility conditions. 
	Let $(U_{d,1},\dots,U_{d,m})$ be the compression of $U$ of depth $d$, i.e. $U_{d,i}(x,y)=U(x,y)(
	\dyainterval{d}{i}
	)$ for every $x,y\in[0,1]$, $i\in[m]$. We now set $U'\coloneqq\AntiCompr{(U_{d,1})^{\Join\dyapartition{d}},\dots,(U_{d,m})^{\Join\dyapartition{d}}}$ and note that $(U',id)$ is a step Latinon. Set $X\coloneqq\{i\in[m]\mid \dyainterval{d}{i}\not\subseteq C_j\; \textrm{ for all } j\in[\ell]\}$ and $Y\coloneqq\bigcup_{i\in X}\dyainterval{d}{i}$. Since $\dyapartition{d}$ consists of smaller intervals than $\mathcal{C}$ we therefore have $|X|\leq\ell$ and $\lambda(Y)\leq\ell/m$. Then the $3$-dimensional measures $\AntiCompr{U_{d,1},\dots,U_{d,m}}$ and $U'$ agree everywhere except on $Q=(Y\times[0,1]\times[0,1])\cup ([0,1]\times Y\times[0,1])\cup([0,1]\times[0,1]\times Y)$, that is, for each $Z\subset [0,1]^3\setminus Q$ we have
	\[\AntiCompr{U_{d,1},\dots,U_{d,m}}(Z)=U'(Z)\;.\]
	Hence by the uniform marginals we get
	$$\cutnL{\AntiCompr{U_{d,1},\dots,U_{d,m}}-U'}\leq 3\ell/m=\tfrac{3\ell qM}{n}<\eps/6\;.$$
	Also by Proposition~\ref{lem:cutdistcomp}\ref{en:caj1} we have that
	$$\cutnL{U-\AntiCompr{U_{d,1},\dots,U_{d,m}}}\leq \tfrac{2}{m}<\eps/6$$
	and therefore
	\begin{equation}
		\cutnL{U-U'}\leq\cutnL{U-\AntiCompr{U_{d,1},\dots,U_{d,m}}}+\cutnL{\AntiCompr{U_{d,1},\dots,U_{d,m}}-U'}<\eps/3\;.\label{eq:suitintervals}
	\end{equation}
	Since $(U',id)$ is a step Latinon we therefore have real values $M_{i,j,k}$, $i,j,k\in[m]$, such that $U'(x,y)(I)=M_{i,j,k}\lambda(I)/\lambda(\dyainterval{d}{k})$ for every $x\in \dyainterval{d}{i}$, $y\in \dyainterval{d}{j}$, and $I\subseteq \dyainterval{d}{k}$. By the uniform marginals of $(U',id)$ we also have for every $i,j\in[m]$ that
	\begin{equation}
		\sum_{k\in[m]}M_{i,j,k}=\sum_{k\in[m]}M_{k,i,j}=\sum_{k\in[m]}M_{i,k,j}=1\;.\label{eq:umarginals}
	\end{equation}
	From now on we will only work with $(U',id)$ and construct a corresponding $3$-partite graph. For this we equipartition each interval $\dyainterval{d}{i}$, $i\in[m]$, further into consecutive intervals $P_{i,1},\dots,P_{i,M}$. Next, we define a corresponding partition of $[n]$ into sets $A_{i,r}$, $i\in[m]$, $r\in[M]$ by setting $A_{i,r}\coloneqq\{1+((i-1)M+r-1)q,\dots,((i-1)M+r)q\}$. Note that $|A_{i,r}|=q=\lambda(P_{i,r})\cdot n$ for every $i\in[m],r\in[M]$. Consider $K_{n,n,n}$; we want to find a $K_3$-decomposition of $K_{n,n,n}$ whose corresponding Latin square is $\eps$-close to $(U',id)$ in the Latinon cut distance. For this we first include each edge of $K_{n,n,n}$ uniformly at random with probability $\rho$ in a reservoir graph $R$. The following claim follows directly from Lemma~\ref{lem:typrandsubgr}. We denote the partite sets of $K_{n,n,n}$ by $V_1$, $V_2$ and $V_3$.
	
	\begin{claim}\label{clm:reservoir}
		We have with probability $1-o(1)$ that
		
		\begin{enumerate}[label=(\roman*)]
			\item $|R[V_i,V_j]|=(1\pm c_1)\rho n^2$ for every $\{i,j\subseteq[3]\}$,
			
			\item $R$ is $(c_1,h)$-typical, and
			
			\item $G'\coloneqq G-R$ is $(c_1,h)$-typical.
		\end{enumerate}
	\end{claim}
	
	We can therefore choose an outcome of the random selection such that Claim~\ref{clm:reservoir} (i)-(iii) hold. In the following we will denote by $A_{i,r}^{(a)}$, $i\in[m],r\in[M],a\in[3]$ the copy of $A_{i,r}$ in $V_a$. We now make a further random selection; every edge will be independently at random selected for one of the graphs $G_{i,j,k}^{r,s,t}$, $i,j,k\in[m]$, $r,s,t\in[M]$. For every $\{i,j\}\subseteq[m]$ and $\{r,s\}\subseteq [M]$ we add an edge $\{u,v\}\in G'[A_{i,r}^{(1)}, A_{j,s}^{(2)}]$ to one of the graphs $G_{i,j,k}^{r,s,t}$ with the probability to be selected for $G_{i,j,k}^{r,s,t}$ being $M_{i,j,k}/M$, $k\in[m],t\in[M]$. Similarly, for every $\{i,k\}\subseteq[m]$ and $\{r,t\}\subseteq [M]$ we add an edge $\{u,v\}\in G'[A_{i,r}^{(1)}, A_{k,t}^{(3)}]$ with probability $M_{i,j,k}/M$ to the subgraph $G_{i,j,k}^{r,s,t}$, $j\in[m],s\in[M]$. Lastly, for every $\{j,k\}\subseteq[m]$ and $\{s,t\}\subseteq [M]$ we add an edge $\{u,v\}\in G'[A_{j,s}^{(2)}, A_{k,t}^{(3)}]$ with probability $M_{i,j,k}/M$ to the subgraph $G_{i,j,k}^{r,s,t}$, $i\in[m],r\in[M]$. Note that by \eqref{eq:umarginals} the probabilities for each edge indeed sum up to~1 and that we may assume $V(G_{i,j,k}^{r,s,t})=A_{i,r}^{(1)}\cup A_{j,s}^{(2)}\cup A_{k,t}^{(3)}$.
	
	\begin{claim}\label{clm:subgraphs}
		We have with probability $1-o(1)$ that for every $i,j,k\in[m]$ and $r,s,t\in[M]$ 
		\begin{enumerate}[label=(\roman*)]
			\item $|E(G_{i,j,k}^{r,s,t}[V_a,V_b])|=(1\pm c_2)(1-\rho)q^2\cdot M_{i,j,k}/M$ for all $\{a,b\}\subseteq[3]$ and
			\item $G_{i,j,k}^{r,s,t}$ is $(c_2,h)$-typical.
		\end{enumerate}
	\end{claim}
	
	\begin{proof}[Proof of Claim~\ref{clm:subgraphs}]
		For fixed $i,j,k\in[m]$ and $r,s,t\in[M]$ we can model the random selection of $G_{i,j,k}^{r,s,t}$ by a random selection  from the $(c_1,h)$-typical graph $G-R$ according to the distribution $\textrm{Bi}(n,M_{i,j,k}/M)$. Therefore by Lemma~\ref{lem:typrandsubgr} we have that (i) and (ii) hold with probability $1-\exp(-\Theta(n))$ for $G_{i,j,k}^{r,s,t}$. By using the union bound we get that the probability that (i) or (ii) does not hold for all $i,j,k\in[m]$ and $r,s,t\in[M]$ is at most $\sum_{i,j,k\in[m]}\sum_{r,s,t\in[M]}\exp(-\Theta(n))\leq m^3M^3\exp(-\Theta(n))=\exp(-\Theta(n))$ which proves the assertion.
	\end{proof}
	
	For every $i,j,k\in[m]$ and $r,s,t\in[M]$, we can now use Theorem~\ref{thm:LSalmost} on $G_{i,j,k}^{r,s,t}$ to get an almost $K_3$-decomposition $T_{i,j,k}^{r,s,t}$ of the edge set such that the leftover graph has maximum degree at most $\alpha q$. After having done this for every $G_{i,j,k}^{r,s,t}$ we denote the set of leftover edges of $G'$ by $F$ and note that $\Delta(F)<(mM)^2\alpha q=\alpha mMn$. Let us write $B\coloneqq\bigcup_{i,j,k\in[m],r,s,t\in[M]}E(T_{i,j,k}^{r,s,t})$. Now, $G''\coloneqq G-B$ is $K_3$-balanced, as it is the leftover of the $K_3$-balanced graph $K_{n,n,n}$ after removing a set of edge-disjoint triangles. Furthermore, note that $G''=R\cup F$ is $(c,h)$-typical with density at least $\rho$, as for every $S\subseteq V_i\cup V_j$, $\{i,j\}\subseteq[3]$, with $|S|\leq h$ we have that $|\overline{N}_{R\cup F}(S)|=(1\pm c_1)\rho^{|S|}n\pm \alpha mMn=(1\pm c)\rho^{|S|}n$. We can therefore use Theorem~\ref{thm:LScomplete} with $\rho$ in place of $d$ to find a $K_3$-decomposition $D'$ of $G''$ and hence $D\coloneqq D'\cup B$ is a $K_3$-decomposition of $K_{n,n,n}$ corresponding to a Latin square which we denote by $L_n$. Note that we have for the linear $3$-uniform hypergraph $H$ which corresponds to $L_n$ that
	\begin{equation}
		\left|E(H[A_{i,r}^{(1)},A_{j,s}^{(2)},A_{k,t}^{(3)}])\right|=(1\pm 2\rho)M_{i,j,k}q^2/M,\textrm{ for all }i,j,k\in[m],r,s,t\in[M]\;.\label{eq:approxedges}
	\end{equation}
	
	It now only remains to show that $L_n$ is $\eps$-close to $(U,id)$. We write $L \coloneqq (W_{L_n},id)$ for the Latinon-representation of $L_n$ over $[0,1]$, and recalling \eqref{eq:suitintervals} and $\deltaL((U,id),L_n)=\deltaL((U,id),L)$, it suffices to show $\cutnL{U'-W_{L_n}} \leq 2\eps/3$. So let $S,T\subseteq [0,1]$ and $V\subseteq [0,1]$ be an interval. Since $1/M\ll\eps$ we can choose $S',T',V'\subseteq [0,1]$ such that $V'$ is an interval and suitable index sets $Q_i^{S'},Q_j^{T'},Q_k^{V'}\subseteq[M]$, for $i,j,k\in[m]$, such that
	\begin{itemize}
		\item $S'=\bigcup_{i\in[m]}\bigcup_{r\in Q_i^{S'}}P_{i,r}$,\\
		\item $T'=\bigcup_{j\in[m]}\bigcup_{s\in Q_j^{T'}}P_{j,s}$,\\
		\item $V'=\bigcup_{k\in[m]}\bigcup_{t\in Q_k^{V'}}P_{k,t}$, and\\
		\item $\lambda(S\triangle S'),\lambda(T\triangle T'),\lambda(V\triangle V')< \eps/3$.
	\end{itemize}
	
	We then have by the uniform marginals property of $(U',id)$ and $L$ that
	\begin{equation}
		\left|\int_S\int_T (U'(x,y)-W_{L_n}(x,y))(V)\diffsymb x \diffsymb y\right|<\left|\int_{S'}\int_{T'} (U'(x,y)-W_{L_n}(x,y))(V')\diffsymb x \diffsymb y\right|+\eps/3\;.\label{eq:approxint}
	\end{equation}
	
	Furthermore, note that $P_{i,r}\coloneqq [(\min\{A_{i,r}\}-1)/n,\max\{A_{i,r}\}/n]$ and therefore the partitions of $L$ and $U'$ conveniently coincide. We then have that
	
	\begin{align*}
		&\left|\int_{S'}\int_{T'} U'(x,y)(V')-W_{L_n}(x,y)(V')\diffsymb x\diffsymb y\right|\\
		=&\left|\sum_{i,j,k\in[m]}\sum_{r\in Q^{S'}_i,s\in Q^{T'}_j,t\in Q^{V'}_k}\int_{P_{i,r}}\int_{P_{j,s}}U'(x,y)(P_{k,t})-W_{L_n}(x,y)(P_{k,t})\diffsymb x\diffsymb y\right|\\
		=&\left|\sum_{i,j,k\in[m]}\sum_{r\in Q^{S'}_i,s\in Q^{T'}_j,t\in Q^{V'}_k}\int_{P_{i,r}}\int_{P_{j,s}}U'(x,y)(P_{k,t})-|E(H[A_{i,r}^{(1)},A_{j,s}^{(2)},A_{k,t}^{(3)}])|/q^2\diffsymb x\diffsymb y\right|\\
		\overset{\eqref{eq:approxedges}}{\leq}&\sum_{i,j,k\in[m]}\sum_{r\in Q^{S'}_i,s\in Q^{T'}_j,t\in Q^{V'}_k}\left|\lambda(P_{i,s})\lambda(P_{j,t})(M_{i,j,k}\lambda(P_{k,t})/\lambda(Q_{d,k})-(1\pm2\rho)M_{i,j,k}/M)\right|\\
		\leq&\sum_{i,j,k\in[m]} M^3\cdot (q^2/n^2) \cdot 2\rho M_{i,j,k}/M\\
		=&\;(2\rho/m^2)\cdot\sum_{i,j,k\in[m]}M_{i,j,k}=2\rho<\eps/3\;,
	\end{align*}
	and thus together with \eqref{eq:approxint} we obtain $\cutnL{U'-W_{L_n}}\le 2\eps/3$ as required.
\end{proof}

\section{Possible further work}\label{sec:future}
Here, we include some questions which would be a natural continuation of the theory we introduced in this paper.

\subsection{Removal lemma for Latin squares}\label{ssec:Removal}
The removal lemma for graphs~\cite{RuSz:TripleSystems,MR1251840,MR1404036} states that if a large $n$-vertex graph $G$ contains only $o(n^{v(H)})$ copies of $H$ then by removing a suitable set of $o(n^2)$ edges from $G$ we obtain a graph with no copies of $H$ at all. Similar statements were investigated for permutations. Cooper~\cite{Co06:permutonreg} formulated a first removal-type lemma for permutations and Klimo\v sov\'a and Kr\'al' \cite{MR3376446} a strong version using a metric called the \emph{Spearman's footrule distance}.
In~\cite{GaHlKuPe22} the first and the third author together with Kun and Pek\'arkov\'a showed how to derive a non-quantitative version of the removal lemma from~\cite{MR3376446} using the theory of permutation limits. It could be that a similar statement for Latin squares, which we state as a conjecture, can be obtained using the theory we introduce here.

\begin{conj}\label{conj:removal}
	Let $k,\ell\in\mathbb N$ and $A\in\mathcal{R}(k,\ell)$. For every $\eps>0$ there exists $\delta>0$ such that the following holds. Suppose that $L$ is a Latin square of order $n$ with $t(A,L)<\delta$. Then there exists a Latin square $\tilde L$ of order $n$ which is $A$-avoiding and for which we have 
\begin{equation}\label{eq:removalchange}
	\sum_{i,j = 1}^n|L(i,j)-\tilde L(i,j)|<\eps n^3\;.
\end{equation}
\end{conj}

The removal lemma in~\cite{GaHlKuPe22} is deduced from another result in~\cite{GaHlKuPe22} which says that permutons avoiding a fixed patterns must have a `one-dimensional' structure. (Note that a finitary version of this one-dimensionality result is Theorem 3 in~\cite{Co06:permutonreg}). It would be interesting to investigate counterparts of this result for Latinons.

%

\subsection{Entropy and counting of Latin squares}\label{ssec:Entropy}
It is a basic calculation in the theory of random graphs that for any $\eps>0$ and any graph $H$, asymptotically almost every (as $n\rightarrow\infty$) graph $G$ of order $n$ satisfies that the density of $H$ in $G$ is $(1\pm \eps)(\frac{1}{2})^{e(H)}$. Rephrased in the limit language: if $(G_n)_n$ is a sequence of uniformly random graphs of order $n$, then almost surely, $(G_n)_n$ converges to a graphon $W\equiv \frac12$.

Here, we prove a Latinon counterpart of this result.
\begin{thm}\label{thm:LimitofUniformLS}
	Given any $\eps>0$, $k,\ell\in\mathbb N$ and $A\in\mathcal{R}(k,\ell)$, asymptotically almost every (as $n\rightarrow \infty$) Latin square $L$ of order $n$ satisfies that $t(A,L)=\tfrac{1}{(k\ell)!}\pm\eps$. Rephrased in the limit language: if $(L_n)_n$ is a sequence of uniformly random Latin squares of order $n$, then almost surely, $(L_n)_n$ converges to a Latinon $W\equiv \lambda_{\restriction [0,1]}$.
\end{thm}
To prove Theorem~\ref{thm:LimitofUniformLS} we use the following nontrivial result of Kwan and Sudakov.
\begin{thm}[{Theorem 3 in~\cite{KwSu10}}]\label{thm:KwSu}
	Asymptotically almost every (as $n\rightarrow \infty$) Latin square $L$ of order $n$ satisfies that for each set $R\subset [n]$ of row indices, for each set $C\subset [n]$ of column indices, and each set $V\subset [n]$ of values, we have
	$$\left|\{(i,j)\in R\times C \mid L_{ij}\in V\}-\frac{|R|\:|C|\:|V|}{n}\right|=O\left(n^{3/2}\log n\right)\;.$$
\end{thm}
With this result, Theorem~\ref{thm:LimitofUniformLS} follows simply by computing the cut distance between Latinons. 
\begin{proof}[Proof of Theorem~\ref{thm:LimitofUniformLS}]
	That is, we show asymptotically almost surely $\deltaL(W,L_n)$ tends to zero as $n$ tends to infinity, where $W \equiv \mathrm{Lebesgue}[0,1]$ and $(L_n,id)$ is the Latinon-representation of a uniformly random Latin square of order $n$. Since we have $$\deltaL(W,L_n) \leq \cutn{W-L_n}= \sup_{\substack{R,C,V\subseteq[0,1],\\V\textrm{ interval}}} \left| \int_{x \in R} \int_{y \in C} (W-L_n)(x,y)(V) \diffsymb x \diffsymb y \right|$$ we may assume that the $R,C,V$ used in the supremum are at least $\Omega(1)$ in measure;
	then Theorem~\ref{thm:KwSu} implies that $\int_{x \in R} \int_{y \in C} L_n(x,y)(V) \diffsymb x \diffsymb y = \lambda(R) \lambda(C) \lambda(V) + O(\frac{\log n}{\sqrt{n}})$ for all such $R,C,V$. Then $\deltaL(W,L_n)=O(\frac{\log n}{\sqrt{n}})$ which tends to zero as $n$ tends to infinity, as required.
\end{proof}

We believe that the concept of Latinons records counts of finite Latin squares. To this end, we introduce the entropy of a Latinon, thus paralleling previous work on graphons~\cite{MR2825532} and for permutons~\cite{KeKrRaWi:Permutations}. Suppose that $L=(W,f)$ is a Latinon. Suppose that for almost every $(x,y)\in\Omega^2$, the measure $W(x,y)$ has a Radon--Nikodym derivative $g(x,y,\cdot)$ with respect to the Lebesgue measure. Then define the \emph{entropy} of $L$ by $\Entropy(L)\coloneqq \int_{x\in\Omega}\int_{y\in\Omega}\int_{z\in[0,1]}g(x,y,z)\log g(x,y,z) \diffsymb x\diffsymb y \diffsymb z$.

The number of Latin squares of order $n$ is $((1\pm o(1))\frac{n}{e^2})^{n^2}$, \cite{MR1871828}. We believe that the entropy tells us how many Latin squares there are close to a given Latinon. More precisely, we believe that there are $((1\pm o(1))\cdot e^{\Entropy(L)}\cdot\frac{n}{e^2})^{n^2}$ Latin squares of order $n$ close to a given Latinon $L$. The easiest way to state this formally is the following.
\begin{conj}
	Suppose that $(A_i\in\bigcup_{k,\ell} \mathcal{R}(k,\ell))_{i\in I}$ is a finite collection of patterns and $(a_i)_{i\in I}$ are reals in $[0,1]$. Let $\alpha\coloneqq \sup \Entropy(L)$, where $L$ ranges through all Latinons with $t(A_i,L)=a_i$ for each $i\in I$. Then for each $\eps>0$ there exists $\delta>0$ and $n_0$ such that for each $n>n_0$ we have that for the number $N$ of Latin squares $K$ of order $n$ satisfying $t(A_i,K)=a_i\pm \delta$ for each $i\in I$ we have
	$$\frac{1}{n^2}\cdot\log\left(\frac{N}{(\frac{n}{e^2})^{n^2}}\right)=\alpha\pm\eps\;.$$
\end{conj}

\subsection{Higher dimensional permutations}\label{ssec:higherdimPerm}
It has been a long programme of Linial to identify higher-dimensional counterparts to common combinatorial objects. Let us introduce his concept of higher dimensional permutations here. For $n\in\mathbb{N}$, any $\{0,1\}$-sequence with $n$ entries which contains $n-1$ zeros and $1$ one is called a \emph{$0$-dimensional permutation}. Now, inductively we say that a $(d+1)$-dimensional $\{0,1\}$-array whose each coordinate is indexed by $[n]$ is a \emph{$d$-dimensional permutation} of order $n$ if restricting this array by fixing arbitrarily one coordinate we get a $(d-1)$-dimensional permutation. So, there is a one-to-one correspondence between 1-dimensional permutations and ordinary permutations on $[n]$ (represented as the corresponding permutation matrix), and 2-dimensional permutations correspond to Latin squares. It is very likely that with the same approach that we developed in this paper, one could create a theory of limits of permutations of any fixed dimension. However, there is a substantial technical challenge. Indeed, our Latinons can be encoded as compressions of graphons. Going up, say, by one dimension, we would have to work with compressions of 3-uniform hypergraphons. It is likely that many steps in our proof would actually need a substantial revision in that setting. For example, there are substantial difficulties with any `cut distance' for 3-uniform hypergraphs (see Section~3 of \cite{Zh15}).

\subsection{Quasirandomness for Latinons, now solved}\label{ssec:quasirandomness}
The concept of pseudorandomness is omnipresent in contemporary mathematics. In graph theory, several related concepts appeared in the late 1980's, \cite{ThomasonPseudorandom,Chung1989,Rodl:Universality}; all these concepts stemmed from properties of the then recent regularity lemma of Szemer\'edi~\cite{Sze:ReguLemma}. It turns out that the (now well established) concept of quasirandomness is best formalised in the language of graphons. In that language, the most comprehensive of the above results, the famous Chung--Graham--Wilson Theorem~\cite{Chung1989} states that, among others, the following properties are equivalent for a graphon $W$,
\begin{enumerate}[label={(\alph*)}]
	\item\label{en:CGW1} $W$ is a constant-$p$ function, almost everywhere;
	\item\label{en:CGW2} $t(H,W)=p^{e(H)}$ for every graph $H$;
	\item\label{en:CGW3} $t(H,W)=p^{e(H)}$ for $H=K_2$ and for $H=C_4$.
\end{enumerate}
Counterparts of the Chung--Graham--Wilson Theorem were obtained for other combinatorial structures, including uniform hypergraphs~\cite{Towsner:QuasirandomHyper,ACHPS:QuasirandomHyper}, classes of oriented graphs~\cite{ChGr:Quasirandomtournaments,Griffiths:QuasirandomOriented}, and permutations~\cite{KrPi}. Of course, the direction \ref{en:CGW1}$\Rightarrow$\ref{en:CGW2}$\Rightarrow$\ref{en:CGW3} is trivial in the graph case as well as in counterparts for other structures. After we published this preprint on the arXiv, Cooper, Kr\'al', Lamaison and Mohr~\cite{CKLM:QuasirandomLatin} characterized constant Latinons, thus proving one of the conjectures posted in the preprint.
\begin{thm}[\cite{CKLM:QuasirandomLatin}]\label{thm:quasiranomLS}
	Suppose that $L=(W,f)$ is a Latinon. Then the following are equivalent.
	\begin{enumerate}[label={(\alph*)}]
		\item\label{en:QL1} $W(x,y)$ is the Lebesgue measure on $[0,1]$ for almost every $(x,y)$,
		\item\label{en:QL2} $t(A,L)=1/(k\ell)!$ for each $k,\ell\in \mathbb{N}$ and each $A\in \mathcal{R}(k,\ell)$,
		\item\label{en:QL3} $t(A,L)=1/6!$ for each $A\in \mathcal{R}(3,2)$.
	\end{enumerate} 
\end{thm}
Moreover, Cooper et al.~show that~\ref{en:QL3} cannot be strengthened by using patterns $\mathcal{R}(2,2)$ or $\mathcal{R}(k,1)$ for any $k\in \mathbb{N}$.

Of course, Theorem~\ref{thm:quasiranomLS} can be equivalently formulated for sequences of Latin squares of growing orders.

\section*{Acknowledgements}
We thank Stefan Glock for helpful discussions about Section~\ref{sec:ApproximatingLatinons}, Leonardo Nagami Coregliano, Dan Král' and Alexander Razborov for discussions about Section~\ref{ssec:NewDiscussion}, and Jan Greb\'ik and Martin Dole\v zal for helping us with Section~\ref{ssec:fromStandardToLocGlob}.
We thank the anonymous referees for their helpful comments. In particular for pointing out that the removal lemmas, suggested in Section~\ref{ssec:Removal} in an earlier version of the paper, were not formulated in a meaningful way.

\appendix

\section{Proof of Proposition~\ref{prop:transformIntoLocalGlobal}}\label{app:localglobal}

\begin{proof}[Proof of Proposition~\ref{prop:transformIntoLocalGlobal}]
	
	We replace $\Omega$ by the unit interval $\mathcal{I}=[0,1]$ with the Lebesgue measure (which we still denote by $\mu$). Again, this can be assumed without loss of generality, since every separable atomless probability space is measure isomorphic to $\mathcal{I}$. Note that the unit interval now appears in the modified setting of Proposition~\ref{prop:transformIntoLocalGlobal} in two different contexts; however the numerical values of elements in the space now replacing $\Omega$ are irrelevant. So, we shall use the symbol $\mathcal{I}$ to denote just this replacement of $\Omega$. For example, we view $f$ as a function from $\mathcal{I}$ to $[0,1]$. In this setting, we want to find $h:[0,1]\times[0,1]\rightarrow \mathcal{I}$ satisfying that $f(h(a,x))=a$.
	
	Let $\Lambda:[0,1]\rightarrow \mathcal{B}$ be a disintegration of $\mu$ through $f$, that is, $\Lambda$ satisfies
	\begin{align}
		\Lambda(r)(f^{-1}(r))&=1 \quad\mbox{for almost every $r\in[0,1]$, and}\\
		\int_{x\in \mathcal{I}} \Lambda(f(x))(A)d\mu&=\mu(A)\quad\mbox{for every $A\subset \mathcal{I}$.}
	\end{align}
	Let $\{C_n\}_{n\in\mathbb{N}}$ be an enumeration of all non-degenerate closed intervals with rational endpoints in $\mathcal{I}$. 
	
	Let $\mathbbm{i}:\Omega\rightarrow\{0,1\}^\mathbb{N}$ be a map defined on its $n$-th coordinate by $(\mathbbm{i}(\omega))_n\coloneqq \mathbbm{1}_{\omega\in C_n}$. Since $\{C_n\}_{n\in\mathbb{N}}$ separates points, we have that $\mathbbm{i}$ is injective. 
	Next, we use $\mathbbm{i}$ to construct a linear order $\prec^*$ on $\mathcal{I}$ as follows. Suppose that $\omega_1,\omega_2\in\mathcal{I}$ are two distinct elements. If $f(\omega_i)<f(\omega_{3-i})$ for some $i\in\{1,2\}$ then put $\omega_i\prec^*\omega_{3-i}$. It remains to order $\omega_1$ and $\omega_2$ in case that $f(\omega_1)=f(\omega_2)$. In that case, we order $\omega_1$ and $\omega_2$ in $\prec^*$ in a way that is consistent with the lexicographic ordering of $\mathbbm{i}(\omega_1)$ and $\mathbbm{i}(\omega_2)$. 
	
	Let $B\subset \mathcal{I}$ be the set of elements $\omega\in\mathcal{I}$ that are atoms for the measure $\Lambda(f(\omega))$. We claim that $B$ is a Borel set. Indeed, this follows upon writing
	$$B=\left\{\omega\in\mathcal{I}\mid \Lambda(f(\omega))(\omega)>0\right\}\;,$$ 
	and recalling that $\Lambda$ is a Borel map.
	
	Let $\mathcal{I}'\coloneqq \mathcal{I}\setminus B$, $f'\coloneqq f\restriction\mathcal{I}'$, and let $\mu'$ be the restriction of $\mu$ on $\mathcal{I}'$. Also, for $r\in[0,1]$, let $\Lambda'(r)$ be the restriction of $\Lambda(r)$ on $\mathcal{I}'$. We then have
	\begin{align}
		\Lambda'(r)(f^{-1}(r))&=1-\Lambda(r)(B) \quad\mbox{for almost every $r\in[0,1]$, and}\\
		\int_{x\in\mathcal{I}'} \Lambda'(f(x))(A)d\mu'&=\mu'(A)\quad\mbox{for every $A\subset \mathcal{I}$.}
	\end{align}
	
	For each $\omega\in\mathcal{I}$, write $\llbracket\omega\rrbracket\coloneqq \{x\in f^{-1}(f(\omega))\mid x\preceq^*\omega\}$.
	
	For the time being, it is convenient to work with a slight modification of the Lebesgue measure on $[0,1]^2$. Define a measure $\gamma$ on $[0,1]^2$ by setting $$\gamma(A)\coloneqq \int_x \Lambda'(x)[0,1] \;\cdot\;\lambda(\{y\in[0,1]\mid(x,y)\in A\}) \diffsymb \lambda(x)\;.$$
	
	We are now going to define a key map $m:\mathcal{I}'\rightarrow [0,1]^2$. Given $\omega\in\mathcal{I}'$, set $$m(\omega)\coloneqq \left(\;
	f(\omega)
	\;,\;
	\frac{1}{1-\Lambda(f(\omega))(B)}\cdot\Lambda'(f(\omega)) \llbracket\omega\rrbracket\;
	\right)\;.$$
	
	The following claim is obvious.
	\begin{claim}The map $m$ is measure preserving with respect to the measures $\mu'$ and $\gamma$.
	\end{claim}
	The following claim is crucial.
	\begin{claim}There exists $\mathcal{I}^*\subset \mathcal{I}'$ with $\mu(\mathcal{I}^*)=\mu(\mathcal{I}')$ so that $m$ is injective on $\mathcal{I}^*$.
	\end{claim}
	\begin{proof}
		Let us introduce a piece of notation. Suppose that $d\in\mathbb{N}$ and $\mathbf{p}\in\{0,1\}^d$. Then define
		\begin{equation*}
			\bigwedge_{\mathbf{p}}(C_n)_n\coloneqq \mathcal{I}'\cap\bigcap_{i\in[d]:\mathbf{p}_i=1}C_i\setminus \left(\bigcup_{i\in[d]:\mathbf{p}_i=0}C_i\right)\;.
		\end{equation*}
		In other words, $\bigwedge_{\mathbf{p}}(C_n)_n$ are the elements of $\mathcal{I}'$ whose encoding using indicator functions of the sets $(C_n)_n$ agrees with $\mathbf{p}$ on the first $d$ coordinates.
		
		For each $t\in [0,1]$ and each $\mathbf{p}\in \bigcup_{d=1}^\infty\{0,1\}^d$, let $w(t,\mathbf{p})\coloneqq f^{-1}(t)\left(\bigwedge_{\mathbf{p}}(C_n)_n\right)$. For a given $\mathbf{p}\in \bigcup_{d=1}^\infty\{0,1\}^d$, let $T_\mathbf{p}=\{t\in[0,1] \mid w(t,\mathbf{p})=0\}$. Observe that $T_\mathbf{p}$ is a Borel subset of $[0,1]$. 
		
		Let $L\subset\mathcal{I}'$ be defined as follows. Take $\omega\in \mathcal{I}'$. We put $\omega$ in $L$ if and only if there exists $k$ such that for each $\ell>k$ either $(\mathbbm{i}(\omega))_\ell=0$, or $w(f(\omega),\mathbf{p})=0$, where $\mathbf{p}\in\{0,1\}^\ell$ has the first $\ell-1$ coordinates as in $\mathbbm{i}(\omega)$ and the $\ell$-th set to $0$. $L$ is a Borel set. Furthermore, observe that for each $t\in[0,1]$, $f^{-1}(t)\cap L$ is at most countable. Since $\Lambda(t)$ is atomless on $\mathcal{I}'$, we conclude that $\Lambda(t)(L)=0$. Integrating, we get that $\mu(L)=0$.
		Set $\mathcal{I}^*\coloneqq \mathcal{I}' \setminus L$. By the above, $\mathcal{I}^*$ has full measure in $\mathcal{I}'$.
		
		We claim that $m$ is injective on $\mathcal{I}^*$. Indeed, suppose that $\omega_1,\omega_2\in\mathcal{I}^*$ are two distinct elements. If $f(\omega_1)\neq f(\omega_2)$ then clearly $m(\omega_1)\neq m(\omega_2)$. So, it remains to assume that $f(\omega_1)= f(\omega_2)=t$ for some $t\in[0,1]$. Without loss of generality, assume that $\mathbbm{i}(\omega_1)$ is smaller in the lexicographic order than $\mathbbm{i}(\omega_2)$, and that the first coordinate on which $\mathbbm{i}(\omega_1)$ and $\mathbbm{i}(\omega_2)$ differ is, say, $h$. Since $\omega_2\not\in L$, we know that there exists $\ell>h$ such that the $\ell$-th coordinate of $\mathbbm{i}(\omega_2)$ is 1, and $w(f(\omega),\mathbf{p})>0$, where $\mathbf{p}\in\{0,1\}^\ell$ has the first $\ell-1$ coordinates equal to $\mathbbm{i}(\omega_2)$ and the $\ell$-th is $0$. We have
		$\llbracket\omega_1\rrbracket\subset \llbracket\omega_2\rrbracket$ with $\llbracket\omega_1\rrbracket\cap f^{-1}(t)\cap\bigwedge_{\mathbf{p}}(C_n)_n=\emptyset$, yet $f^{-1}(t)\cap\bigwedge_{\mathbf{p}}(C_n)_n\subset \llbracket\omega_2\rrbracket$.  Since $\Lambda(t)\left(\bigwedge_{\mathbf{p}}(C_n)_n\right)=w(t,\mathbf{p})>0$, we have that 
		$$\Lambda(t)\left(\llbracket\omega_2\rrbracket\right)\ge \Lambda(t)\left(\llbracket\omega_1\rrbracket\right)+\Lambda(t)\left(\bigwedge_{\mathbf{p}}(C_n)_n\right)>\Lambda(t)\left(\llbracket\omega_1\rrbracket\right)\;.$$
		In particular, we conclude that $m(\omega_2)$ and $m(\omega_1)$ are different on the second coordinate, as was needed.
	\end{proof}
	Now, let $J\subset [0,1]^2$ be the image of $\mathcal{I}^*$ under $m$. As $m$ is measure preserving and $\mathcal{I}^*$ is conull, we conclude that $J$ is conull in $[0,1]^2$ with respect to $\gamma$. Let $g:[0,1]^2\rightarrow \mathcal{I}^*$ be defined as the inverse map to $m$ on $J$, and arbitrarily on the rest (which is null).
	
	We are now ready to define the map $h$. Suppose that $(a,x)\in [0,1]^2$ is given. If $x< 1-\Lambda(a)(B)$ then define $h(a,x)\coloneqq g(a,\frac{x}{1-\Lambda(a)(B)})$. Otherwise, let $h(a,x)$ be the atom $\zeta\in B\cap f^{-1}(a)$ for which $$1-\Lambda(a)(B)+\Lambda(a)\left(\llbracket \zeta\rrbracket\right)> x >1-\Lambda(a)(B)+\Lambda(a)\left(\llbracket \zeta\rrbracket\setminus \{\zeta\}\right)\;.$$
	Such an atom exists for almost every choice of $x\in[0,1]$; the exceptions are at most countably many values for which the down-interval from that atom has exactly that measure.
	
	Obviously, $h$ has all the properties we want.
\end{proof}

\section{Proof of Theorem~\ref{thm:LSalmost}}\label{app:nibble}

We will derive Theorem~\ref{thm:LSalmost} from the following theorem of Alon and Yuster~\cite{AY05}. Let $H=(V,E)$ be an $r$-uniform hypergraph. Recall that the \emph{degree} of a vertex $v$ is the number of edges in $H$ that contain $v$, and that the \emph{codegree} of two distinct vertices $u$ and $v$ is the number of edges in $H$ that contain $u$ and $v$. We write $\delta(H)$ for the minimum degree, $\Delta(H)$ for the maximum degree, and $\Delta_2(H)$ for the maximum codegree. Let $\mathcal{F}\subseteq 2^V$ and let $0\leq \alpha\leq1$. A matching $M$ in $H$ is \emph{$(\alpha,\mathcal{F})$-perfect} if for each $F\in\mathcal{F}$, at least $\alpha|F|$ vertices of $F$ are covered by $M$. Also we write $s(\mathcal{F})=\min_{F\in\mathcal{F}}\{|F|\}$ and $g(H)=\Delta(H)/\Delta_2(H)$.

\begin{thm}[Theorem 1.2~\cite{AY05}]\label{thm:AY}
	For an integer $r\geq 2$, a real $C>1$ and a real $\eps>0$ there exists a real $\mu=\mu(r, C, \eps)$ and a real $K=K(r, C, \eps)$ so that the following holds:  If the $r$-uniform hypergraph $H=(V,E)$ on $N$ vertices satisfies:
	\begin{enumerate}
		\item $\delta(H)\geq(1-\mu)\Delta(H)$,
		\item $g(H)>\max\{1/\mu ,  K(\log(N))^6\}$,
	\end{enumerate}
	then for every $\mathcal{F}\subseteq 2^V$ with $|\mathcal{F}|\leq C^{g(H)^{1/(3r-3)}}$ and with $s(\mathcal{F})\geq 5g(H)^{1/(3r-3)}\log(|\mathcal{F}|g(H))$ there is a $(1-\eps,\mathcal{F})$-perfect matching in $H$.
\end{thm}

\begin{proof}[Proof of Theorem~\ref{thm:LSalmost}]
	Recall that we have for the given constants that
	\[1/n\ll c\ll \alpha \ll d\;.\]
	Let $\mu$ and $K$ be the constants given by Theorem~\ref{thm:AY} applied with $r=3,C=2$ and $\eps=\alpha$. Therefore we may also assume that
	\[1/n\ll c\ll \mu,1/K\;.\]
	We define $H=(V,E)$ to be the $3$-uniform hypegraph with vertex set $V=E(G)$ and edge set $E=\{\{\{x,y\},\{y,z\},\{x,z\}\}\subseteq V\mid \{x,y\},\{y,z\},\{x,z\}\in E(G)\}$. In other words, the vertices of $H$ are the edges of $G$ and the hyperedges of $H$ are the triangles of $G$. Furthermore, for every $x\in V(G)$ we define $F_x\coloneqq \{\{x,y\}\mid y\in N_G(x)\}$, and set $\mathcal{F}\coloneqq\{F_x\mid x\in V(G)\}$. Note that $\mathcal{F}\subseteq 2^V$, $|\mathcal{F}|=3n$,
	and that by the $(c,h)$-typicality of $G$,
	\[s(\mathcal{F})=\min_{F\in\mathcal{F}}\{|F|\}=2(1\pm 2c)dn\;.\]
	Furthermore, also by the $(c,h)$-typicality of $G$, we have that $|V|=|E(G)|=3(1\pm c)dn^2$, $\Delta(H)=(1\pm c)d^2n$, $\delta(H)=(1\pm c)d^2n$, and $\Delta_2(H)=1$. Hence
	\[g(H)=\Delta(H)/\Delta_2(H)=(1\pm c)d^2n/1=(1\pm c)d^2n\;.\]
	We check the conditions for applying Theorem~\ref{thm:AY}.
	\[\delta(H)\geq (1-c)d^2n\geq (1-\mu)(1+c)d^2n\geq(1-\mu)\Delta(H)\;.\]
	Also
	\[g(H)\geq (1-c)d^2n\geq\max\{1/\mu,K(\log(3(1\pm c)dn^2))^6\}\;.\]
	Also note that
	\[|\mathcal{F}|\leq 3n\leq C^{((1+ c)d^2n)^{1/6}}\;,\]
	and
	\[s(\mathcal{F})\geq2(1- 2c)dn\geq n^{1/3}\geq5((1+ c)d^2n)^{1/6}\log(3n(1+ c)d^2n)\geq 5g(H)^{1/6}\log(|\mathcal{F}|g(H))\;.\]
	Hence there exists a $(1-\alpha,\mathcal{F})$-perfect matching $T$ for $H$. Note that elements of $T$ correspond to edge-disjoint triangles in $G$ and that each vertex $x\in V(G)$ lies in at least $(1-\alpha)|F_x|$ triangles of $T$. Therefore
	\[\Delta(G-T)\leq\alpha\Delta(G)\leq\alpha n\;.\]
\end{proof}

\bibliographystyle{amsplain}
\bibliography{LLS}

\begin{dajauthors}
\begin{authorinfo}[fg]
  Frederik Garbe\\
  Masaryk University\\
  Brno, Czech Republic\\
  garbe\imageat{}fi\imagedot{}muni\imagedot{}cz \\
  \url{https://www.fi.muni.cz/~garbe/}
\end{authorinfo}
\begin{authorinfo}[rh]
  Robert Hancock\\
  Institut f\"ur Informatik \\ 
  Heidelberg University\\ 
  Heidelberg, Germany \\ 
  Previous affiliation: \\
  Masaryk University\\
  Brno, Czech Republic\\
  hancock\imageat{}informatik\imagedot{}uni-heidelberg\imagedot{}de \\
  \url{https://sites.google.com/view/robert-hancock/} 
\end{authorinfo}
\begin{authorinfo}[jh]
  Jan Hladk\'y\\
  Researcher\\
  Institute of Computer Sciences of the Czech Academy of Sciences\\
  Prague, Czechia\\
  hladky\imageat{}cs\imagedot{}cas\imagedot{}cz\\
  \url{https://www.cs.cas.cz/~hladky/}
\end{authorinfo}
\begin{authorinfo}[ms]
  Maryam Sharifzadeh\\
  Institutionen för matematik och matematisk statistik\\
  Umeå universitet\\
  Umeå, Sweden\\
  maryam\imagedot{}sharifzadeh\imageat{}umu\imagedot{}se \\
  \url{https://www.umu.se/personal/maryam-sharifzadeh/}
\end{authorinfo}
\end{dajauthors}

\end{document}